\documentclass[11pt]{amsart}

\usepackage{amsmath, amsthm, amsfonts, amscd, flafter, graphicx, verbatim, rotating, pkmacros, url}

\newcommand\kh{\widetilde{\textit{Kh}}}
\newcommand\ckh{\widetilde{\textit{CKh}}}

\newcommand\hf{\widehat{HF}}

\newcommand\zz{\mathbb{Z}}

\DeclareMathOperator{\rk}{rk}

\newcommand\ftwo{\mathbb{F}_2}
\newcommand\grt{\text{gr}^{(2)}}
\newcommand\grq{\text{gr}^\mathbb{Q}}

\newcommand\bg{\boldsymbol\gamma}
\newcommand\bbg{\breve{\boldsymbol\gamma}}
\newcommand\bbd{\breve{\boldsymbol\delta}}

\newcommand\cptb{\overline{\mathbb{CP}}^2}

\newcommand\cm{\check{C}}
\newcommand\cmhat{\tilde{C}}
\newcommand\hm{\Hto}
\newcommand\hmb{\Hto_\bullet}
\newcommand\hmfrom{\widehat{\textit{HM}}_\bullet}
\newcommand\hmbar{\overline{\textit{HM}}_\bullet}
\newcommand\hmhat{\widetilde{\textit{HM}}_\bullet}

\newcommand\mw{\check{m}(W)}

\newcommand\qq{\check{Q}}
\newcommand\Aa{\check{A}}

\newcommand\Ii{\mathit{0}}
\newcommand\If{\mathit{1}}

\newcommand\wij{W_{IJ}}

\newcommand\pij{P_{IJ}}
\newcommand\qij{Q_{IJ}}

\newcommand\yi{Y_I}
\newcommand\yj{Y_J}

\newcommand\ddd{\check{\partial}}

\newcommand\dd{\check{D}}
\newcommand\dij{\check{D}^I_J}
\newcommand\dkj{\check{D}^K_J}
\newcommand\dik{\check{D}^I_K}
\newcommand\dii{\check{D}^I_I}
\newcommand\djj{\check{D}^J_J}

\newcommand\dooij{D^o_o ({}^I_J)}
\newcommand\dosij{D^o_s ({}^I_J)}
\newcommand\dusij{D^u_s ({}^I_J)}
\newcommand\duoij{D^u_o ({}^I_J)}
\newcommand\dbssij{\bar{D}^s_s ({}^I_J)}
\newcommand\dbsuij{\bar{D}^s_u ({}^I_J)}
\newcommand\dbusij{\bar{D}^u_s ({}^I_J)}
\newcommand\dbuuij{\bar{D}^u_u ({}^I_J)}

\newcommand\dooik{D^o_o ({}^I_K)}
\newcommand\dosik{D^o_s ({}^I_K)}

\newcommand\dbssik{\bar{D}^s_s ({}^I_K)}
\newcommand\dbsuik{\bar{D}^s_u ({}^I_K)}

\newcommand\dookj{D^o_o ({}^K_J)}
\newcommand\doskj{D^o_s ({}^K_J)}
\newcommand\duskj{D^u_s ({}^K_J)}
\newcommand\duokj{D^u_o ({}^K_J)}
\newcommand\dbsskj{\bar{D}^s_s ({}^K_J)}

\newcommand\dbsuil{\bar{D}^s_u ({}^I_L)}

\newcommand\duslk{D^u_s ({}^L_K)}
\newcommand\duolk{D^u_o ({}^L_K)}

\newcommand\dbsukm{\bar{D}^s_u ({}^K_M)}

\newcommand\dusmj{D^u_s ({}^M_J)}
\newcommand\duomj{D^u_o ({}^M_J)}

\newcommand\duoktj{D^u_o({}^{K_2}_{J})}
\newcommand\dbuukokt{\bar{D}^u_u({}^{K_1}_{K_2})}
\newcommand\dbsuiko{\bar{D}^s_u({}^{I}_{K_1})}

\newcommand\doo{\partial^o_o}
\newcommand\dos{\partial^o_s}
\newcommand\duo{\partial^u_o}
\newcommand\dus{\partial^u_s}

\newcommand\ess{\bar{\partial}^s_s}
\newcommand\esu{\bar{\partial}^s_u}
\newcommand\eus{\bar{\partial}^u_s}
\newcommand\euu{\bar{\partial}^u_u}

\newcommand\mma{\check{m}({}^{00}_{01})}
\newcommand\mmb{\check{m}({}^{01}_{11})}
\newcommand\mmc{\check{m}({}^{00}_{10})}
\newcommand\mmd{\check{m}({}^{10}_{11})}

\newcommand\hh{\check{H}}

\newcommand\aij{\check{A}^I_J}

\newcommand\aoo{A^o_o}
\newcommand\aos{A^o_s}
\newcommand\auo{A^u_o}
\newcommand\aus{A^u_s}

\newcommand\abss{\bar{A}^s_s}
\newcommand\absu{\bar{A}^s_u}
\newcommand\abus{\bar{A}^u_s}
\newcommand\abuu{\bar{A}^u_u}

\newcommand\aooij{\aoo({}^I_J)}
\newcommand\aosij{\aos({}^I_J)}
\newcommand\auoij{\auo({}^I_J)}
\newcommand\ausij{\aus({}^I_J)}
\newcommand\abssij{\abss({}^I_J)}
\newcommand\absuij{\absu({}^I_J)}
\newcommand\abusij{\abus({}^I_J)}
\newcommand\abuuij{\abuu({}^I_J)}

\newcommand\absuik{\absu({}^I_K)}

\newcommand\auokj{\auo({}^K_J)}
\newcommand\auskj{\aus({}^K_J)}

\newcommand\auokoj{\auo({}^{K_1}_J)}
\newcommand\absuikt{\absu({}^I_{K_2})}

\newcommand\boo{B^o_o}
\newcommand\bos{B^o_s}
\newcommand\buo{B^u_o}
\newcommand\bus{B^u_s}

\newcommand\bbss{\bar{B}^s_s}
\newcommand\bbsu{\bar{B}^s_u}
\newcommand\bbus{\bar{B}^u_s}
\newcommand\bbuu{\bar{B}^u_u}

\newcommand\bij{\check{B}^I_J}

\newcommand\bbsua{\bar{B}^s_u({}^{00}_{01})}

\newcommand\buob{B^u_o({}^{01}_{11})}
\newcommand\busb{B^u_s({}^{01}_{11})}

\newcommand\bbsub{\bar{B}^s_u({}^{01}_{11})}

\newcommand\buoc{B^u_o({}^{00}_{10})}

\newcommand\bbsuc{\bar{B}^s_u({}^{00}_{10})}

\newcommand\buod{B^u_o({}^{10}_{11})}
\newcommand\busd{B^u_s({}^{10}_{11})}

\newcommand\Eoo{E^o_o}
\newcommand\Eos{E^o_s}
\newcommand\Euo{E^u_o}
\newcommand\Eus{E^u_s}

\newcommand\Ebss{\bar{E}^s_s}
\newcommand\Ebsu{\bar{E}^s_u}
\newcommand\Ebus{\bar{E}^u_s}
\newcommand\Ebuu{\bar{E}^u_u}

\newcommand\Foo{F^o_o}
\newcommand\Fos{F^o_s}
\newcommand\Fuo{F^u_o}
\newcommand\Fus{F^u_s}

\newcommand\Fbss{\bar{F}^s_s}
\newcommand\Fbsu{\bar{F}^s_u}
\newcommand\Fbus{\bar{F}^u_s}
\newcommand\Fbuu{\bar{F}^u_u}

\newcommand\hoo{H^o_o}
\newcommand\hos{H^o_s}
\newcommand\huo{H^u_o}
\newcommand\hus{H^u_s}

\newcommand\kss{\bar{H}^s_s}
\newcommand\ksu{\bar{H}^s_u}
\newcommand\kus{\bar{H}^u_s}
\newcommand\kuu{\bar{H}^u_u}

\newcommand\qoo{Q^o_o}
\newcommand\qos{Q^o_s}
\newcommand\quo{Q^u_o}
\newcommand\qus{Q^u_s}

\newcommand\rss{\bar{Q}^s_s}
\newcommand\rsu{\bar{Q}^s_u}
\newcommand\rus{\bar{Q}^u_s}
\newcommand\ruu{\bar{Q}^u_u}

\newcommand\mm{\check{m}}

\newcommand\moo{m^o_o}
\newcommand\mos{m^o_s}
\newcommand\muo{m^u_o}
\newcommand\mus{m^u_s}

\newcommand\nss{\bar{m}^s_s}
\newcommand\nsu{\bar{m}^s_u}
\newcommand\nus{\bar{m}^u_s}
\newcommand\nuu{\bar{m}^u_u}

\newcommand\mook{m^o_o({}^{0}_{1})}
\newcommand\mosk{m^o_s({}^{0}_{1})}
\newcommand\muok{m^u_o({}^{0}_{1})}
\newcommand\musk{m^u_s({}^{0}_{1})}

\newcommand\nssk{\bar{m}^s_s({}^{0}_{1})}
\newcommand\nsuk{\bar{m}^s_u({}^{0}_{1})}
\newcommand\nusk{\bar{m}^u_s({}^{0}_{1})}
\newcommand\nuuk{\bar{m}^u_u({}^{0}_{1})}

\newcommand\mool{m^o_o({}^{1}_{\infty})}
\newcommand\mosl{m^o_s({}^{1}_{\infty})}
\newcommand\muol{m^u_o({}^{1}_{\infty})}
\newcommand\musl{m^u_s({}^{1}_{\infty})}

\newcommand\nssl{\bar{m}^s_s({}^{1}_{\infty})}
\newcommand\nsul{\bar{m}^s_u({}^{1}_{\infty})}
\newcommand\nusl{\bar{m}^u_s({}^{1}_{\infty})}
\newcommand\nuul{\bar{m}^u_u({}^{1}_{\infty})}

\newcommand\mooa{m^o_o({}^{00}_{01})}
\newcommand\moob{m^o_o({}^{01}_{11})}
\newcommand\mooc{m^o_o({}^{00}_{10})}
\newcommand\mood{m^o_o({}^{10}_{11})}

\newcommand\mosa{m^o_s({}^{00}_{01})}
\newcommand\mosb{m^o_s({}^{01}_{11})}
\newcommand\mosc{m^o_s({}^{00}_{10})}
\newcommand\mosd{m^o_s({}^{10}_{11})}

\newcommand\muoa{m^u_o({}^{00}_{01})}
\newcommand\muob{m^u_o({}^{01}_{11})}
\newcommand\muoc{m^u_o({}^{00}_{10})}
\newcommand\muod{m^u_o({}^{10}_{11})}

\newcommand\musa{m^u_s({}^{00}_{01})}
\newcommand\musb{m^u_s({}^{01}_{11})}
\newcommand\musc{m^u_s({}^{00}_{10})}
\newcommand\musd{m^u_s({}^{10}_{11})}

\newcommand\nssa{\bar{m}^s_s({}^{00}_{01})}
\newcommand\nssb{\bar{m}^s_s({}^{01}_{11})}
\newcommand\nssc{\bar{m}^s_s({}^{00}_{10})}
\newcommand\nssd{\bar{m}^s_s({}^{10}_{11})}

\newcommand\nsua{\bar{m}^s_u({}^{00}_{01})}
\newcommand\nsub{\bar{m}^s_u({}^{01}_{11})}
\newcommand\nsuc{\bar{m}^s_u({}^{00}_{10})}
\newcommand\nsud{\bar{m}^s_u({}^{10}_{11})}

\newcommand\nusa{\bar{m}^u_s({}^{00}_{01})}
\newcommand\nusb{\bar{m}^u_s({}^{01}_{11})}
\newcommand\nusc{\bar{m}^u_s({}^{00}_{10})}
\newcommand\nusd{\bar{m}^u_s({}^{10}_{11})}

\newcommand\nuua{\bar{m}^u_u({}^{00}_{01})}
\newcommand\nuub{\bar{m}^u_u({}^{01}_{11})}
\newcommand\nuuc{\bar{m}^u_u({}^{00}_{10})}
\newcommand\nuud{\bar{m}^u_u({}^{10}_{11})}

\newcommand\mmooa{m^o_o({}^{000}_{001})}
\newcommand\mmoob{m^o_o({}^{000}_{010})}
\newcommand\mmooc{m^o_o({}^{000}_{100})}

\newcommand\mmooj{m^o_o({}^{011}_{111})}
\newcommand\mmook{m^o_o({}^{101}_{111})}
\newcommand\mmool{m^o_o({}^{110}_{111})}

\newcommand\mmosa{m^o_s({}^{000}_{001})}
\newcommand\mmosb{m^o_s({}^{000}_{010})}
\newcommand\mmosc{m^o_s({}^{000}_{100})}

\newcommand\mmosj{m^o_s({}^{011}_{111})}
\newcommand\mmosk{m^o_s({}^{101}_{111})}
\newcommand\mmosl{m^o_s({}^{110}_{111})}

\newcommand\mmuoa{m^u_o({}^{000}_{001})}
\newcommand\mmuob{m^u_o({}^{000}_{010})}
\newcommand\mmuoc{m^u_o({}^{000}_{100})}

\newcommand\mmuoj{m^u_o({}^{011}_{111})}
\newcommand\mmuok{m^u_o({}^{101}_{111})}
\newcommand\mmuol{m^u_o({}^{110}_{111})}

\newcommand\mmusa{m^u_s({}^{000}_{001})}
\newcommand\mmusb{m^u_s({}^{000}_{010})}
\newcommand\mmusc{m^u_s({}^{000}_{100})}

\newcommand\mmusj{m^u_s({}^{011}_{111})}
\newcommand\mmusk{m^u_s({}^{101}_{111})}
\newcommand\mmusl{m^u_s({}^{110}_{111})}

\newcommand\nnssa{\bar{m}^s_s({}^{000}_{001})}
\newcommand\nnssb{\bar{m}^s_s({}^{000}_{010})}
\newcommand\nnssc{\bar{m}^s_s({}^{000}_{100})}

\newcommand\nnssj{\bar{m}^s_s({}^{011}_{111})}
\newcommand\nnssk{\bar{m}^s_s({}^{101}_{111})}
\newcommand\nnssl{\bar{m}^s_s({}^{110}_{111})}

\newcommand\nnsua{\bar{m}^s_u({}^{000}_{001})}
\newcommand\nnsub{\bar{m}^s_u({}^{000}_{010})}
\newcommand\nnsuc{\bar{m}^s_u({}^{000}_{100})}
\newcommand\nnsud{\bar{m}^s_u({}^{001}_{011})}
\newcommand\nnsue{\bar{m}^s_u({}^{001}_{101})}
\newcommand\nnsuf{\bar{m}^s_u({}^{010}_{011})}
\newcommand\nnsug{\bar{m}^s_u({}^{010}_{110})}
\newcommand\nnsuh{\bar{m}^s_u({}^{100}_{101})}
\newcommand\nnsui{\bar{m}^s_u({}^{100}_{110})}
\newcommand\nnsuj{\bar{m}^s_u({}^{011}_{111})}
\newcommand\nnsuk{\bar{m}^s_u({}^{101}_{111})}
\newcommand\nnsul{\bar{m}^s_u({}^{110}_{111})}

\newcommand\nnusa{\bar{m}^u_s({}^{000}_{001})}
\newcommand\nnusb{\bar{m}^u_s({}^{000}_{010})}
\newcommand\nnusc{\bar{m}^u_s({}^{000}_{100})}

\newcommand\nnusj{\bar{m}^u_s({}^{011}_{111})}
\newcommand\nnusk{\bar{m}^u_s({}^{101}_{111})}
\newcommand\nnusl{\bar{m}^u_s({}^{110}_{111})}

\newcommand\nnuua{\bar{m}^u_u({}^{000}_{001})}
\newcommand\nnuub{\bar{m}^u_u({}^{000}_{010})}
\newcommand\nnuuc{\bar{m}^u_u({}^{000}_{100})}

\newcommand\nnuuj{\bar{m}^u_u({}^{011}_{111})}
\newcommand\nnuuk{\bar{m}^u_u({}^{101}_{111})}
\newcommand\nnuul{\bar{m}^u_u({}^{110}_{111})}

\newcommand\hhooa{H^o_o({}^{000}_{011})}
\newcommand\hhoob{H^o_o({}^{000}_{101})}
\newcommand\hhooc{H^o_o({}^{000}_{110})}
\newcommand\hhood{H^o_o({}^{001}_{111})}
\newcommand\hhooe{H^o_o({}^{010}_{111})}
\newcommand\hhoof{H^o_o({}^{100}_{111})}

\newcommand\hhosa{H^o_s({}^{000}_{011})}
\newcommand\hhosb{H^o_s({}^{000}_{101})}
\newcommand\hhosc{H^o_s({}^{000}_{110})}
\newcommand\hhosd{H^o_s({}^{001}_{111})}
\newcommand\hhose{H^o_s({}^{010}_{111})}
\newcommand\hhosf{H^o_s({}^{100}_{111})}

\newcommand\hhuoa{H^u_o({}^{000}_{011})}
\newcommand\hhuob{H^u_o({}^{000}_{101})}
\newcommand\hhuoc{H^u_o({}^{000}_{110})}
\newcommand\hhuod{H^u_o({}^{001}_{111})}
\newcommand\hhuoe{H^u_o({}^{010}_{111})}
\newcommand\hhuof{H^u_o({}^{100}_{111})}

\newcommand\hhusa{H^u_s({}^{000}_{011})}
\newcommand\hhusb{H^u_s({}^{000}_{101})}
\newcommand\hhusc{H^u_s({}^{000}_{110})}
\newcommand\hhusd{H^u_s({}^{001}_{111})}
\newcommand\hhuse{H^u_s({}^{010}_{111})}
\newcommand\hhusf{H^u_s({}^{100}_{111})}

\newcommand\kkssa{\bar{H}^s_s({}^{000}_{011})}
\newcommand\kkssb{\bar{H}^s_s({}^{000}_{101})}
\newcommand\kkssc{\bar{H}^s_s({}^{000}_{110})}
\newcommand\kkssd{\bar{H}^s_s({}^{001}_{111})}
\newcommand\kksse{\bar{H}^s_s({}^{010}_{111})}
\newcommand\kkssf{\bar{H}^s_s({}^{100}_{111})}

\newcommand\kksua{\bar{H}^s_u({}^{000}_{011})}
\newcommand\kksub{\bar{H}^s_u({}^{000}_{101})}
\newcommand\kksuc{\bar{H}^s_u({}^{000}_{110})}
\newcommand\kksud{\bar{H}^s_u({}^{001}_{111})}
\newcommand\kksue{\bar{H}^s_u({}^{010}_{111})}
\newcommand\kksuf{\bar{H}^s_u({}^{100}_{111})}

\newcommand\kkusa{\bar{H}^u_s({}^{000}_{011})}
\newcommand\kkusb{\bar{H}^u_s({}^{000}_{101})}
\newcommand\kkusc{\bar{H}^u_s({}^{000}_{110})}
\newcommand\kkusd{\bar{H}^u_s({}^{001}_{111})}
\newcommand\kkuse{\bar{H}^u_s({}^{010}_{111})}
\newcommand\kkusf{\bar{H}^u_s({}^{100}_{111})}

\newcommand\kkuua{\bar{H}^u_u({}^{000}_{011})}
\newcommand\kkuub{\bar{H}^u_u({}^{000}_{101})}
\newcommand\kkuuc{\bar{H}^u_u({}^{000}_{110})}
\newcommand\kkuud{\bar{H}^u_u({}^{001}_{111})}
\newcommand\kkuue{\bar{H}^u_u({}^{010}_{111})}
\newcommand\kkuuf{\bar{H}^u_u({}^{100}_{111})}

\newcommand\gc{\check{G}}

\newcommand\mmma{\check{m}({}^{000}_{001})}
\newcommand\mmmb{\check{m}({}^{000}_{010})}
\newcommand\mmmc{\check{m}({}^{000}_{100})}
\newcommand\mmmj{\check{m}({}^{011}_{111})}
\newcommand\mmmk{\check{m}({}^{101}_{111})}
\newcommand\mmml{\check{m}({}^{110}_{111})}

\newcommand\hhha{\check{H}({}^{000}_{011})}
\newcommand\hhhb{\check{H}({}^{000}_{101})}
\newcommand\hhhc{\check{H}({}^{000}_{110})}
\newcommand\hhhd{\check{H}({}^{001}_{111})}
\newcommand\hhhe{\check{H}({}^{010}_{111})}
\newcommand\hhhf{\check{H}({}^{100}_{111})}

\newcommand\goo{G^o_o}
\newcommand\gos{G^o_s}
\newcommand\guo{G^u_o}
\newcommand\gus{G^u_s}

\newcommand\gbss{\bar{G}^s_s}
\newcommand\gbsu{\bar{G}^s_u}
\newcommand\gbus{\bar{G}^u_s}
\newcommand\gbuu{\bar{G}^u_u}

\newcommand\amf{\mathfrak{a}}
\newcommand\bmf{\mathfrak{b}}
\newcommand\cmf{\mathfrak{c}}

\newcommand\Mod{M(\amf,W^*,\bmf)}

\newcommand\Modzplus{M_z^+(\amf,W^*,\bmf)}
\newcommand\Modzijplus{M_z^+(\amf,W_{IJ}^*,\bmf)}

\newcommand\Modplus{M^+(\amf,W^*,\bmf)}

\newcommand\Modz{M_z(\amf,W^*,\bmf)}

\newcommand\Modzij{M_z(\amf,\wij^*,\bmf)}
\newcommand\Modzrij{M_z^{\text{red}}(\amf,\wij^*,\bmf)}

\newcommand\Modzprijplus{M_z^{\text{red}+}(\amf,\wij^*,\bmf)_{\pij}}

\newcommand\Modzpij{\Modzij_{\pij}}
\newcommand\Modzprij{\Modzrij_{\pij}}
\newcommand\Modzqijplus{M_z^+(\amf,\wij^*,\bmf)_{Q_{IJ}}}

\newcommand\qqij{\check{Q}^I_J}

\newcommand\Modzpijplus{\Modzijplus_{\pij}}

\newcommand\modzij{m_z(\amf, \wij^*,\bmf)}
\newcommand\modzrij{\bar{m}_z(\amf, \wij^*,\bmf)}

\newcommand\Cmf{\mathfrak{C}}
\newcommand\Comf{\mathfrak{C}^o}
\newcommand\Csmf{\mathfrak{C}^s}
\newcommand\Cumf{\mathfrak{C}^u}

\newcommand\bswij{\mathcal{B}^\sigma(\wij)}

\providecommand{\abs}[1]{\lvert#1\rvert}

\newtheorem{theorem}{Theorem}[section]
\newtheorem{lemma}[theorem]{Lemma}

\newtheorem{conjecture}[theorem]{Conjecture}
\newtheorem{corollary}[theorem]{Corollary}
\newtheorem{proposition}[theorem]{Proposition}

\theoremstyle{definition}

\newtheorem{remark}[theorem]{Remark}

\newtheorem{example}[theorem]{Example}

\title[A Link Surgery Spectral Sequence]{A Link Surgery Spectral Sequence in \\ Monopole {F}loer Homology}

\author{Jonathan M.\ Bloom}

\begin{document}

\begin{abstract}
To a link $L \subset S^3$, we associate a spectral sequence whose $E^2$ page is the reduced Khovanov homology of $L$ and which converges to a version of the monopole Floer homology of the branched double cover.  The pages $E^k$ for $k \geq 2$ depend only on the mutation equivalence class of $L$.  We define a mod 2 grading on the spectral sequence which interpolates between the $\delta$-grading on Khovanov homology and the mod 2 grading on Floer homology.  We also derive a new formula for link signature that is well-adapted to Khovanov homology.

More generally, we construct new bigraded invariants of a framed link in a 3-manifold as the pages of a spectral sequence modeled on the surgery exact triangle.  The differentials count monopoles over families of metrics parameterized by permutohedra.  We utilize a connection between the topology of link surgeries and the combinatorics of graph associahedra.  This also yields simple realizations of permutohedra and associahedra, as refinements of hypercubes.
\end{abstract}

\maketitle

\thispagestyle{empty}

\tableofcontents

\listoffigures


\noindent
\textit{Keywords:} monopole Floer homology, gauge theory, Khovanov homology, branched double cover, framed link, surgery, permutohedron, graph associahedron. \\ \\
The author was supported by NSF grant DMS-0739392. \\ \\
{\small \textsc{Department of Mathematics, Columbia University, New York, NY 10027, USA}} \\
\textit{E-mail address:} \texttt{jbloom@math.columbia.edu} \\
\textit{URL:} \url{http://math.columbia.edu/~jbloom/}
\newpage


\section{Introduction}

Monopole Floer homology is a gauge-theoretic invariant defined via Morse theory on the Chern-Simons-Dirac functional.  As such, the underlying chain complex is generated by Seiberg-Witten monopoles over a 3-manifold, and the differential counts monopoles over the product of the 3-manifold with $\mathbb{R}$.  In \cite{kmos}, a surgery exact triangle is associated to a triple of surgeries on a knot in a 3-manifold (for a precursor in instanton Floer homology, see \cite{bd}, \cite{fl}).

This paper has two main objectives. The first is to construct a link surgery spectral sequence in monopole Floer homology, generalizing the exact triangle.  This is a spectral sequence which starts at the monopole Floer homology of a hypercube of surgeries on $Y$ along $L$, and converges to the monopole Floer homology of $Y$ itself. The differentials count monopoles on 2-handle cobordisms equipped with families of metrics parameterized by polytopes called permutohedra.  Those metrics parameterized by the boundary of the permutohedra are stretched to infinity along collections of hypersurfaces, representing surgered 3-manifolds.  The monopole counts satisfy identities obtained by viewing the map associated to each polytope as a null-homotopy for the map associated to its boundary.  Note that this can be seen as analogue of Ozsv\'ath and Szab\'o's link surgery spectral sequence for Heegaard Floer homology \cite{osz12}. There, the differentials count pseudo-holomorphic polygons in Heegaard multi-diagrams, and they satisfy $A_\infty$ relations which encode degenerations of conformal structures on polygons.  In Section \ref{sec:hfproof}, we develop a dictionary between these gauge-theoretic and symplectic perspectives.

Our construction introduces a number of techniques that we hope will be of more general use.  In Sections \ref{sec:surface} and \ref{sec:graph}, we couple the topology of 2-handle cobordisms arising from link surgeries to the combinatorics of polytopes called graph associahedra \cite{carr}.  For the chain-level Floer maps induced by 2-handle cobordisms, these polytopes encode a mixture of commutativity and associativity up to homotopy.  We hope this coupling, and its relationship to finite product lattices, will be of independent interest to algebraists and combinatorists.  As one application, we obtain simple, recursive realizations of permutohedra as refinements of associahedra, which in turn refine hypercubes (see Figures \ref{fig:assocu} through \ref{fig:assocr}).  Curiously, these realizations are predicted by the ``sliding-the-point'' proof of the naturality of the $U_\dagger$ action in Floer theory.

Our construction of polytopes of metrics was inspired by the pentagon of metrics in the proof of the surgery exact triangle \cite{kmos}.    However, to make use of these polytopes, we must deal effectively with the notorious mix of interior, boundary-stable, and boundary-unstable critical points involved in the construction of the monopole Floer complex.  To this end, we systematize the construction of maps associated to cobordisms equipped with polytopes of metrics, as well as the vanishing operators which count ends of 1-dimensional moduli spaces.  This includes the construction of the usual monopole Floer differentials, cobordism maps, and homotopies as special cases, as well as the operators used in the proof of the surgery exact triangle, which we reorganize in Section \ref{sec:surgexact}.  More generally, we prove that the homotopy type of the link surgery spectral sequence is independent of analytic choices, which may be viewed as a gauge-theoretic analogue of the invariance of $A_\infty$ homotopy type in symplectic geometry \cite{seid}.  In particular, the higher pages are themselves invariants of a framed link in a 3-manifold.

Khovanov homology is a powerful new invariant of classical links in the 3-sphere, arising from representation theory \cite{kh1}.  It is defined combinatorially and categorifies the Jones polynomial.  Our second main objective is to construct a spectral sequence from reduced Khovanov homology to a version of the monopole Floer homology of the branched double cover.  While here the strategy in the Heegaard Floer case may be translated fairly directly, we instead present a more global identification of the $E^1$ page with the reduced Khovanov complex, based on our ``thriftier'' construction of the reduced odd Khovanov complex \cite{b1}.  The intermediate pages are link invariants as well.

Beyond these objectives, we refine the link surgery spectral sequence, as well as its specialization to Khovanov homology, in ways that were not previously known for the Heegaard Floer version, but now follow by parallel arguments.  In particular, we equip the spectral sequence with a mod 2 grading which interpolates between a shift of the $\delta$-grading on Khovanov homology and the mod 2 grading on monopole Floer homology, thereby refining the known rank inequality.
We also derive a new formula for the signature of a link that is well-adapted to Khovanov homology and may be of independent interest.

We have recently learned that  Kronheimer and Mrowka have established a similar connection between Khovanov homology and a version of instanton Floer homology.  It would be interesting to understand the relationship between our approaches. \\

\subsection{Statement of results.}
All monopole Floer homology and Khovanov homology groups are considered over the 2-element field $\ftwo$.  Our notation is consistent with the definitive reference \cite{km}.  In particular, $\cm(Y)$ and $\hmb(Y)$ denote the ``to'' version of the complex and homology group associated to $Y$, while $\mm(W)$ and $\hmb(W)$ denote the chain-level and homology-level homomorphisms associated to a cobordism $W$.

In order to motivate the statement of the link surgery spectral sequence, we first recall the surgery exact triangle.  Let $Y$ be a closed, oriented 3-manifold, equipped with a knot $K$ with framing $\lambda$ and meridian $\mu$.   Orient $\lambda$ and $\mu$ as simple closed curves on the torus boundary of the complement of a neighborhood of $K$, so that the algebraic intersection numbers of the triple $(\lambda, \lambda + \mu, \mu)$ satisfy
$$(\lambda \cdot (\lambda + \mu)) = ((\lambda + \mu) \cdot \mu) = (\mu \cdot \lambda) = -1.$$
Let $Y(0)$ and $Y(1)$ denote the result of surgery on $K$ with respect to $\lambda$ and $\lambda + \mu$, respectively.  In \cite{kmos}, Kronheimer, Mrowka, Ozsv\'ath, and Szab\'o prove that the mapping cone
$$\cm(Y(0)) \xrightarrow{\mm(W(01))} \cm(Y(1))$$
is quasi-isomorphic to the monopole Floer complex $\cm(Y)$, where $\mm(W(01))$ is the chain map induced by the elementary 2-handle cobordism $W(01)$ from $Y(0)$ to $Y(1)$.  The associated long exact sequence on homology is known as the surgery exact triangle.  However, we can also frame the result in another way.  As in \cite{osz12}, if we filter by the index $I$ in $Y(I)$, then the mapping cone induces a spectral sequence with
$$E^1 = \hm(Y(0)) \textstyle\bigoplus \hm(Y(1))$$
and
$$d^1 = \hmb(W(01)),$$
which converges by the $E^2$ page to $\hm(Y)$.

The link surgery spectral sequence generalizes this interpretation of the exact triangle to the case of an $l$-component framed link $L \subset Y$.  For each $I = (m_1, \dots, m_l)$ in the hypercube $\{0,1\}^l$, let $Y(I)$ denote the result of performing $m_i$-surgery on the component $K_i$.  For $I<J$, let $W(IJ)$ denote the associated cobordism, composed of $(w(J)-w(I))$ 2-handles.  The (iterated) mapping cone now takes the form of a hypercube complex
\begin{align*}
X = \bigoplus_{I \in \{0,1\}^l} \cm(Y(I))
\end{align*}
with differential $\dd$ given by the sum of components $\dij : \cm(Y(I)) \to \cm(Y(J))$ for all $I \leq J$.  We filter $X$ by vertex weight $w(I)$, defined as the sum of the coordinates of $I$.  The component $\dd^I_I$ is the usual differential on $\cm(Y(I))$, whereas for $I < J$, the component $\dij$ counts monopoles on $W(IJ)$ over a family of metrics parametrized by the permutohedron of dimension $w(J) - w(I) - 1$.  We define this family in Section \ref{sec:surface} and construct $(X,\dd)$ in Section \ref{sec:linksurg}.  In Section \ref{sec:linksurg2}, we complete the proof of:

\begin{theorem}
\label{thm:linksurg1}
The filtered complex $(X,\dd)$ induces a spectral sequence with $E^1$ page given by
\begin{align*}
E^1 = \bigoplus_{I \in \{0,1\}^l} \hmb(Y(I))
\end{align*}
and $d^1$ differential given by
\begin{align*} 
d^1 = \bigoplus_{\substack{I < J \in \{0,1\}^l \\ w(J) - w(I) = 1}} \hmb(W(IJ)).
\end{align*}
The spectral sequence converges by the $E^{l+1}$ page to $\hmb(Y)$ and comes equipped with an absolute mod 2 grading $\check \delta$ which coincides on $E^\infty$ with that of $\hmb(Y)$.  In addition, each page has an integer grading $\check{t}$ induced by vertex weight.  The differential $d^k$ shifts $\check \delta$ by one and increases $\check{t}$ by $k$.
\end{theorem}

The complex $(X,\dd)$ depends on a choice of regular metric and perturbation of the monopole equations on the full cobordism from $Y(\{0\}^l)$ to $Y(\{1\}^l)$.  For any two such choices, we produce a homotopy equivalence which induces a graded isomorphism between the associated $E^1$ pages.

\begin{theorem}
\label{thm:inv}
For each $i \geq 1$, the $( \check t, \check \delta)$-graded vector space $E^i$ is an invariant of the framed link $L \subset Y$.
\end{theorem}

In fact, reduced Khovanov homology over $\ftwo$ arises as such an invariant.  To frame this properly, in Section \ref{sec:umap}, we introduce another version of monopole Floer homology, pronounced ``H-M-hat'' and denoted $\hmhat$.  By analogy with $\hf$ in Heegaard Floer homology, we define $\hmhat(Y)$ as the homology of the mapping cone of $U_\dagger: \cm(Y) \to \cm(Y)[1]$, where  $U_\dagger$ induces the usual even endomorphism on $\hmb(Y)$.  It follows that $\hmhat(Y)$ inherits a mod 2 grading, and we prove a version of Theorem \ref{thm:linksurg1} for $\hmhat$ as well.  Note that $\hmhat(Y,\mathfrak{s})$ should agree with the sutured monopole Floer homology group $SFH(Y-B^3,\mathfrak{s})$ relative to the equatorial suture, as defined in \cite{km2}.  The latter is given by $\hmb(Y \# (S^1\times F), \mathfrak{s} \# \mathfrak{s}_c)$, where $F$ is an orientable surface of genus $g > 2$ and $\mathfrak{s}_c$ is the canonical spin$^c$-structure with $\left< c_1(\mathfrak{s_c}),[F] \right> = 2g - 2$.  In particular, the rank of $\hmhat(Y)$ is finite over $\ftwo$ and invariant under orientation reversal.

For an oriented link $L \subset S^3$, let $\kh(L)$ denote the reduced Khovanov homology of $L$.  To a diagram of $L$, we will associate a framed link $\mathbb{L}$ in the branched double cover with reversed orientation, denoted $-\Sigma(L)$.  Applying the $\hmhat$ version of Theorem \ref{thm:linksurg1}, we prove:

\begin{theorem}
\label{thm:1}
The link surgery spectral sequence for $\mathbb{L} \subset -\Sigma(L)$ has $E^2$ page isomorphic to $\kh(L)$ and converges by the $E^{l+1}$ page to $\hmhat(-\Sigma(L))$.
\end{theorem}

While the construction of this spectral sequence depends on a choice of diagram for $L$, as well as analytic data, Theorem \ref{thm:1} implies that the $E^2$ and $E^{\infty}$ pages are actually link invariants.  These pages are also insensitive to Conway mutation, since this is true of Khovanov homology over $\ftwo$ as well as branched double covers.  More generally, we prove:

\begin{theorem}
\label{thm:linkinv}
For each $k \geq 2$, the $(\check t, \check \delta)$-graded vector space $E^k$ depends only on the mutation equivalence class of $L$.
\end{theorem}

The analytic invariance described in Theorem \ref{thm:inv} is crucial here.  As explained in Section 9.2, Reidemeister invariance is then an immediate consequence of Baldwin's proof in the Heegaard Floer case \cite{b}, whereas mutation invariance follows from our proof in the Heegaard Floer case \cite{b1}.  Note that both Heegaard Floer proofs, in turn, depend on Roberts' work on invariance with respect to isotopy, handleslides, and stabilization in Heegaard multi-diagrams \cite{lrob3}, and Baldwin's work on invariance with respect to almost complex data \cite{b}.

Recall that Khovanov homology is graded by two integers, the homological grading $t$ and the quantum grading $q$.  We may repackage this as a rational $(t,\delta)$-bigrading, where $$\delta = q/2 - t$$ marks the diagonals of slope two in the $(t,q)$-plane.  On the other hand, monopole Floer homology has a canonical mod 2 grading and decomposes over the set of spin$^c$ structures.  Using the $\check \delta$ grading on the spectral sequence, we derive the first result relating these finer features of monopole or Heegaard Floer homology to those of Khovanov homology, leading to a refinement of the rank inequality
$$\text{rk} \ \kh(L) \geq \text{rk} \ \hmhat({-\Sigma(L))} \geq \det(L).$$
Let $\hmhat^0(Y)$ and $\hmhat^1(Y)$ denote the even and odd graded pieces of $\hmhat(Y)$, respectively.  Let $\kh^0(L)$ and $\kh^1(L)$ denote the even and odd graded pieces of $\kh(L)$ with respect to the integer grading $\delta -(\sigma(L) + \nu(L))/2$.  The terms $\sigma(L)$ and $\nu(L)$ refer to the signature and nullity of $L$, respectively.  Our convention is that the signature of the right-handed trefoil is +2 (that is, minus the signature of a Seifert matrix).  Recall that $\nu(L) = b_1(\Sigma(L))$.
\begin{theorem}
\label{thm:grading}
The $\check \delta$ grading on the spectral sequence coincides with
$$\delta - \frac{1}{2}(\sigma(L) + \nu(L)) \mod 2$$ on the $E^2$ page.  Thus, the rank inequality may be refined to
\begin{align*}
\rk \, \kh^0(L) &\geq \rk \, \hmhat^0({-\Sigma(L))} \geq \det(L)\\
\rk \, \kh^1(L) &\geq \rk \, \hmhat^1({-\Sigma(L))}.
\end{align*}
Furthermore, the $\check \delta$ Euler characteristic of each page is given by $\det(L)$.
\end{theorem}

In particular, all the differentials on the spectral sequence shift $\delta + 2\zz$ by one.  We conclude:
\begin{corollary}
\label{cor:thin}
If $\kh(L)$ is supported on a single diagonal, then the spectral sequence collapses at the $E^2$ page.  In particular, $\hmhat(-\Sigma(L))$ is supported in even grading and has rank $\det(L)$, with one generator in each spin$^c$ structure.
\end{corollary}
In fact, $\kh(L)$ is supported on the single diagonal $\delta = \sigma(L)/2$ whenever $L$ is quasi-alternating \cite{mo}.  This is consistent with Theorem \ref{thm:grading}, since quasi-alternating links have non-zero determinant, and therefore vanishing nullity.

Finally, we present a new formula for $\sigma(L)$.  It follows quickly from the proof of Theorem \ref{thm:grading}, which in turn invokes the Gordon-Litherland signature formula (see \cite{gl}).

\begin{figure}[htp]
\centering
\includegraphics[width=150mm]{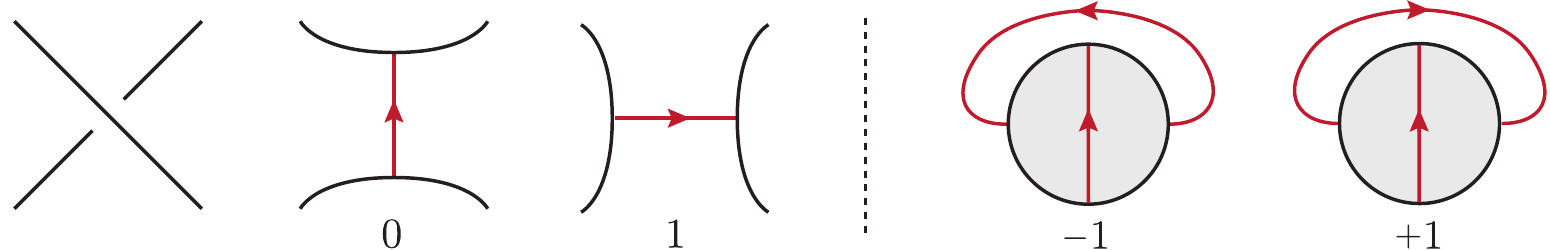}
\caption{Resolution and arc-linking conventions for the signature formula.}
\label{fig:arcs}
\end{figure}

As in \cite{b1}, we first assign a symmetric $l \times l$ matrix $A = (a_{ij})$ to an oriented, connected diagram $\mathcal{D}$ with $l$ numbered crossings as follows.  Fix a vertex $I^* = (m_1^*,\dots,m_l^*)$ such that the resolution $\mathcal{D}(I^*)$ consists of one circle (such a resolution may be obtained by resolving along a spanning tree of the black graph).  Now place a small arbitrarily-oriented arc $x_i$ across each resolved crossing, as shown at left in Figure \ref{fig:arcs}.  Two arcs are {\em linked} if their endpoints are interleaved around the circle.  For each pair of linked arcs $\{x_i, x_j\}$, set $a_{ij} = a_{ji} = \pm 1$ according to the convention at right in Figure \ref{fig:arcs}.  Here we are viewing $\mathcal{D}(I^*)$ on the sphere, so that the outside arc may be pulled to the bottom.  Let $a_{ii} = (-1)^{m^*_i}$ and set all remaining entries to zero.  Let $n_-$ denote the number of negative crossings in $\mathcal{D}$.

\begin{proposition}
\label{prop:newsig}
With the above conventions:
$$\sigma(L) = \sigma(A) + w(I^*) - n_- \qquad \quad \det(L) = \left|\det(A)\right| \qquad \quad  \nu(L) = \nu(A)$$
\end{proposition}

The proof and an example are given at the end of Section \ref{sec:khgrade}.  Unlike the Goeritz matrix, the non-zero entries of $A$ are all $\pm1$.  Remarkably, a deep structure theorem in graph theory due to W. H. Cunningham implies that $A$ alone determines the mutation equivalence class of $L$, the framed isotopy type of $\mathbb{L}$, and therefore $E^i$ for all $i \geq 1$ (see \cite{b1}, \cite{chm}, and Remark \ref{rem:encode}).

\subsection{Philosophy and future directions.}
In outline, the identity of the $E^2$ page in Theorem \ref{thm:1} may be established as follows (our detailed proof in Section \ref{sec:khovanov} is along slightly different lines).  To a diagram of a link $L \subset S^3$, we associate a framed link $\mathbb{L} \subset -\Sigma(L)$.  With respect to $\mathbb{L}$, the link surgery hypercube of 3-manifolds $Y(I)$ and 4-dimensional 2-handle cobordisms $W(IJ)$ is precisely the branched double cover of the Khovanov hypercube of 1-manifolds $\mathcal{D}(I) \subset S^3$ and 2-dimensional 1-handle cobordisms $F_{IJ} \subset S^3 \times [0,1]$, as illustrated for the trefoil knot in Figures \ref{fig:trefoil} and \ref{fig:trefcube}.  Furthermore, the functor $\hmhat$ and the functor $\mathit{CKh}$ underlying Khovanov's {\em unreduced} theory over $\ftwo$ fit into a commutative square of functors.

\begin{figure}[htp]
\centering
\includegraphics[width=110mm]{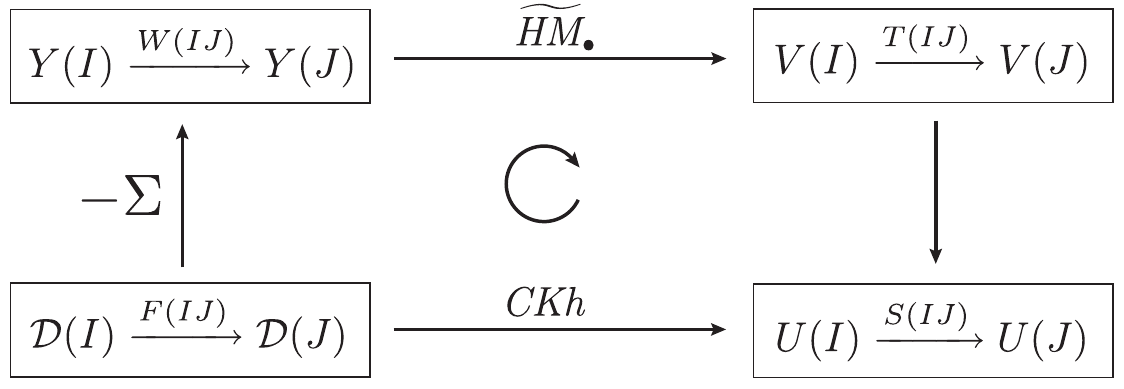}
\label{fig:squarea}
\end{figure}

\noindent Here $S(IJ):U(I) \to U(J)$ and $T(IJ):V(I) \to V(J)$ represent the induced maps of $\ftwo$-vector spaces with respect to each theory.  If we replace $\mathit{CKh}$ with the {\em reduced} Khovanov functor $\ckh$ over $\ftwo$, then the vertical arrow at right induces an equivariant isomorphism of vector spaces.  Consequently, we may identify the complex $(E^1, d^1)$ with $\ckh(L)$, and hence $E^2$ with $\kh(L)$.

In fact, the entire commutative diagram admits a more elementary and unified description, which illuminates {\em why} the functors $\hmhat$ and $\mathit{CKh}$ are connected in the first place.   Both horizontal arrows may be regarded as an instance of a TQFT described by Donaldson in \cite{don}.  The algebraic basis for his construction is as follows.  To an $\ftwo$-vector space $U$, we associate the exterior algebra $\Lambda^* U$.
To a linear map $i: \Gamma \to U_0 \oplus U_1$, we associate a map $| \Gamma | : \Lambda^* U_0 \to \Lambda^* U_1$ defined as follows.  Let $k$ and $n$ be the dimensions of $\Gamma$ and $U_0$, respectively.  By taking the exterior product of the images of the elements in any basis of $\Gamma$, we obtain an element of $\Lambda^k(U_0 \oplus U_1)$, which may be regarded as a map via the series of isomorphisms
\begin{align*}
\Lambda^k(U_0 \oplus U_1) &\cong \textstyle\bigoplus_{i = 0}^k \, \Lambda^i U_0 \textstyle\otimes   \Lambda^{k-i} U_1 \\
&\cong  \textstyle\bigoplus_{i = 0}^k \, (\Lambda^{n-i} U_0)^* \textstyle\otimes   \Lambda^{k-i} U_1 \\
&\cong \textstyle\bigoplus_{i = 0}^k \, \text{Hom}(\Lambda^{n-i} U_0, \Lambda^{k-i} U_1).
\end{align*}
A composition law holds provided that a certain transversality condition is met.

To a manifold $M$, Donaldson associates the exterior algebra $\Lambda^* H^1(M)$.  To a cobordism $N: M_0 \to M_1$, he associates the map
$$| H^1(N) | : \Lambda^* H^1(M_0) \to \Lambda^* H^1(M_1),$$
obtained from the restriction $H^1(N) \to H^1(\partial N) \cong H^1(M_0) \oplus H^1(M_1)$.  If we denote his TQFT by $\Lambda^* H^1$, then the above commutative diagram may be replaced with:

\begin{figure}[htp]
\centering
\includegraphics[width=110mm]{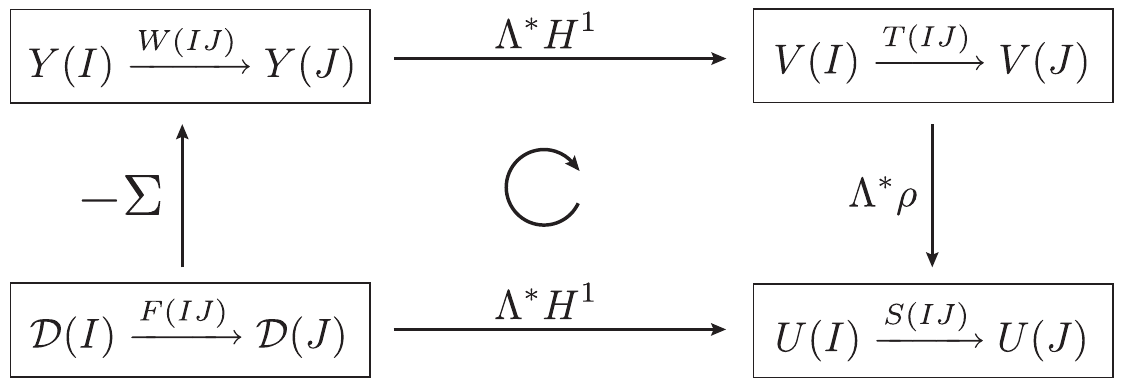}
\label{fig:squareb}
\end{figure}

\noindent  Regarding the vertical arrow at right, note that for any link $L$ with a basepoint, there is a natural map $\widetilde{H}_0(L) \cong H_1(S^3, L) \to H_1(\Sigma(L))$ which takes a relative 1-cycle to its preimage.  Dually, there is a map $$\rho : H^1(\Sigma(L)) \to H^1(S^3, L) \cong \widetilde{H}^0(L) \subset H^0(L) \cong H^1(L),$$ which induces this arrow.  Note that when $F(IJ)$ has positive genus, both $S(IJ)$ and $T(IJ)$ vanish since the restriction map from $H^1(F(IJ))$ has non-trivial kernal.

The equivalence of the two commutative diagrams may be understood as follows.  The manifold $Y(I)$ admits a metric of positive scalar curvature, so it follows from Proposition 36.1.3 of \cite{km} that $\hmhat(Y(I))$ is the cohomology of the Picard torus $H^1(Y(I),\mathbb{R}) / H^1(Y(I),\zz)$, parameterizing flat $U(1)$-connections on $Y(I)$ modulo gauge.  The cobordism $W(IJ)$ also admits a metric of positive scalar curvature, and indeed, it follows from Corollary \ref{cor:s1s2} herein that the induced map $T(IJ)$ coincides with the map on cohomology induced by the correspondence between Picard tori defined by flat connections over $\wij$.  As Donaldson observes, the map on cohomology induced by such a correspondence is encoded in the above TQFT.  Along the bottom row, this TQFT is determined by its Frobenius algebra, which one may easily check is the same as the one underlying Khovanov homology over $\ftwo$.

A version of the spectral sequence with $\mathbb{Z}$ coefficients is work in progress.  Indeed, Donaldson's TQFT lifts to $\zz$ by equipping cobordisms with homology orientations, and we then recover the monopole Floer and {\em odd} Khovanov functors in our commutative diagram.   We also expect that the $(\check t, \check \delta)$ bigrading can be lifted and shifted to an invariant rational $(t,\delta)$ bigrading on the higher pages (compare with Conjecture 1.1 of \cite{b}).
For such links, we then obtain a ``higher'' Khovanov homology and Jones polynomial on each page $E^k$ for $k > 2$.  These are conjectured for a family of torus knots in Section \ref{sec:torus}.  Watson has shown that Khovanov homology is not an invariant of the branched double cover \cite{w}.   One wonders whether the same is true of the bigrading on the $E^\infty$ page, and what this bigrading encodes.

Finally, we strongly suspect that the entire construction of the link surgery spectral sequence is functorial.  Broadly speaking, to a framed surface in a cobordism of 3-manifolds, we would like to associate a map between the spectral sequences associated to the framed links in the 3-manifold ends.  In the Khovanov specialization, the branched double cover of a classical link cobordism provides a 4-manifold $W$ with boundary, and on the $E^\infty$ page, we expect to see the associated monopole Floer map.  Perhaps one could then associate a framed surface in $W$ to a combinatorial description of the link cobordism, so that Khovanov's combinatorially-defined link cobordism maps appear on the $E^1$ page.  This would provide a Floer-theoretic realization of the functoriality of Khovanov's theory. \\

\noindent \textbf{Acknowledgements.} It is my pleasure to thank: my advisor Peter Ozsv\'ath, for recommending this engaging topic and providing invaluable guidance throughout.  Tim Perutz, for first suggesting the relevance of Donaldson's TQFT, and for sharing his expertise in all applicable fields.  Adam Knapp and John Baldwin, for their endless enthusiasm and many insightful discussions.  Satyan Devadoss, for a wonderful visit and for \cite{d}.  Dave Bayer, for his interest in my work and expert programming assistance.  I am also indebted to the following people for their wisdom and support: Eli Grigsby, Robert Lipshitz, Max Lipyanskiy, Maria Lissitsyna, James Stasheff, Dylan Thurston, and Rumen Zarev.

\newpage


\section{Hypercubes and permutohedra}
\label{sec:surface}

This section involves no Floer homology whatsoever, but rather surgery theory and Kirby calculus as described in Part 2 of \cite{gs}.  In particular, with respect to a 2-handle $D^2 \times D^2$, the terms {\em core}, {\em cocore}, and {\em attaching region} will refer to the subsets $D^2 \times \{0\}$, $\{0\} \times D^2$, and $\partial D^2 \times D^2$, respectively.

Let $Y$ be a closed, oriented 3-manifold, equipped with an $l$-component, framed link $L = K_1 \bigcup \cdots \bigcup K_l$, and let $Y^\prime$ denote the result of (integral) surgery on $L$.  There is a standard oriented cobordism $W: Y \to Y^\prime$, built by thickening $Y$ to $[0,1] \times Y$ and attaching 2-handles $h_i$ to $\{1\} \times Y$ by identifying the attaching region of $h_i$ with a tubular neighborhood $\nu(K_i)$ in accordance with the framing.  The diffeomorphism type of $W$ is insensitive to whether the handles are attached simultaneously as above, or instead in a succession of batches which express $W$ as a composite cobordism.  Our goal in this section is to construct a family of metrics on $W$, parameterized by the permutohedron $P_l$, which smoothly interpolates between all ways of expressing $W$ as a composite cobordism.

In order to keep track of the $l!$ ways to build up $W$ one handle at a time, we introduce the hypercube poset $\{0,1\}^l$, with $I = (m_1,\dots,m_l) \leq J = (m^\prime_1,\dots,m^\prime_l)$ if and only if $m_i \leq m^\prime_i$ for all $1\leq i \leq l$.  $J$ is called an {\em immediate successor} of $I$ if there is a $k$ such that $m_k = 0$, $m^\prime_k = 1$, and $m_i = m^\prime_i$ for all $i \neq k$.   We define a {\em path} of length $k$ from $I$ to $J$ to be a sequence of immediate successors $I= I_0 < I_1 < \cdots < I_k = J$.  The {\em weight} of a vertex $I$ is given by $w(I) = \sum_{i=1}^l m_i$.  We use $\Ii$ and $\If$ as shorthand for the initial and terminal vertices of $\{0,1\}^l$, which we call {\em external}.  The other $2^l - 2$ vertices will be called {\em internal}.  A totally ordered subset of a poset is called a {\em chain}.  A chain is {\em maximal} if it is not properly contained in any other chain.  In $\{0,1\}^l$, the maximal chains are precisely the paths from $\Ii$ to $\If$, with each such path determined by its internal vertices.

To each vertex $I$, we associate the 3-manifold $\yi$ obtained by surgery on the framed sublink $$L(I) = \bigcup_{\{i \,|\, m_i = 1 \}} K_i$$ in $Y$.  Note that the remaining components of $L$ constitute a framed link in $\yi$.

\begin{remark}
The 3-manifold denoted $Y(I)$ in the introduction and in \cite{osz12} is obtained from $\yi$ by shifting forward one frame in the surgery exact triangle for each component of $L$.  We will use $\yi$ throughout and address this discrepancy in Remark \ref{rem:bdcview}.
\end{remark}

We regard $\{\yi \, | \, I \in \{0,1\}^l \}$ as a poset isomorphic to $\{0,1\}^l$, with $Y_\Ii$ and $Y_\If$ external and the rest internal.  To a pair of vertices $(I,J)$ with $I < J$, we associate the 2-handle cobordism $$\wij = \yi \times [0,1] \bigcup_{\{i \,|\, m_i = 0, \, m^\prime_i = 1\}} h_i$$ from $\yi$ to $\yj$.  In particular, if $J$ is an immediate successor of $I$, then $\wij$ is an {\em elementary cobordism}, given by a single 2-handle addition.  More generally, $\wij$ will be the composition of $w(J)-w(I)$ elementary cobordisms.

In order to quantify how far two vertices are from being ordered, we define a symmetric function $\rho$ on pairs of vertices by
$$\rho(I,J) = \min \left\{ \left| \{i \,|\, m_i > m^\prime_i \} \right|, \left| \{i \,|\, m^\prime_i > m_i \} \right| \right\}.$$
Note that $\rho(I,J) = 0$ if and only if $I$ and $J$ are ordered.  In this case, $Y_I$ and $Y_J$ are disjoint:

\begin{lemma}
\label{lem:hyper}
The full set of $2^l - 2$ internal hypersurfaces $\yi$ can be simultaneously embedded in the interior of the cobordism $W$ so that the following conditions hold:

\renewcommand{\labelenumi}{(\roman{enumi})}
\begin{enumerate}
\item the hypersurfaces in any subset are pairwise disjoint if and only if they form a chain.  In this case, cutting on $Y_{I_1} < Y_{I_2} < ... < Y_{I_k}$ breaks $W$ into the disjoint union $$W_{\Ii I_1} \  \textstyle\coprod \ W_{I_1I_2} \ \coprod  \ \cdots \  \coprod \ W_{I_k \If}.$$

\item distinct hypersurfaces $\yi$ and $\yj$ intersect in exactly $\rho(I,J)$ disjoint tori.
\end{enumerate}
\end{lemma}

\begin{figure}[htp]
\centering
\includegraphics[width=150mm]{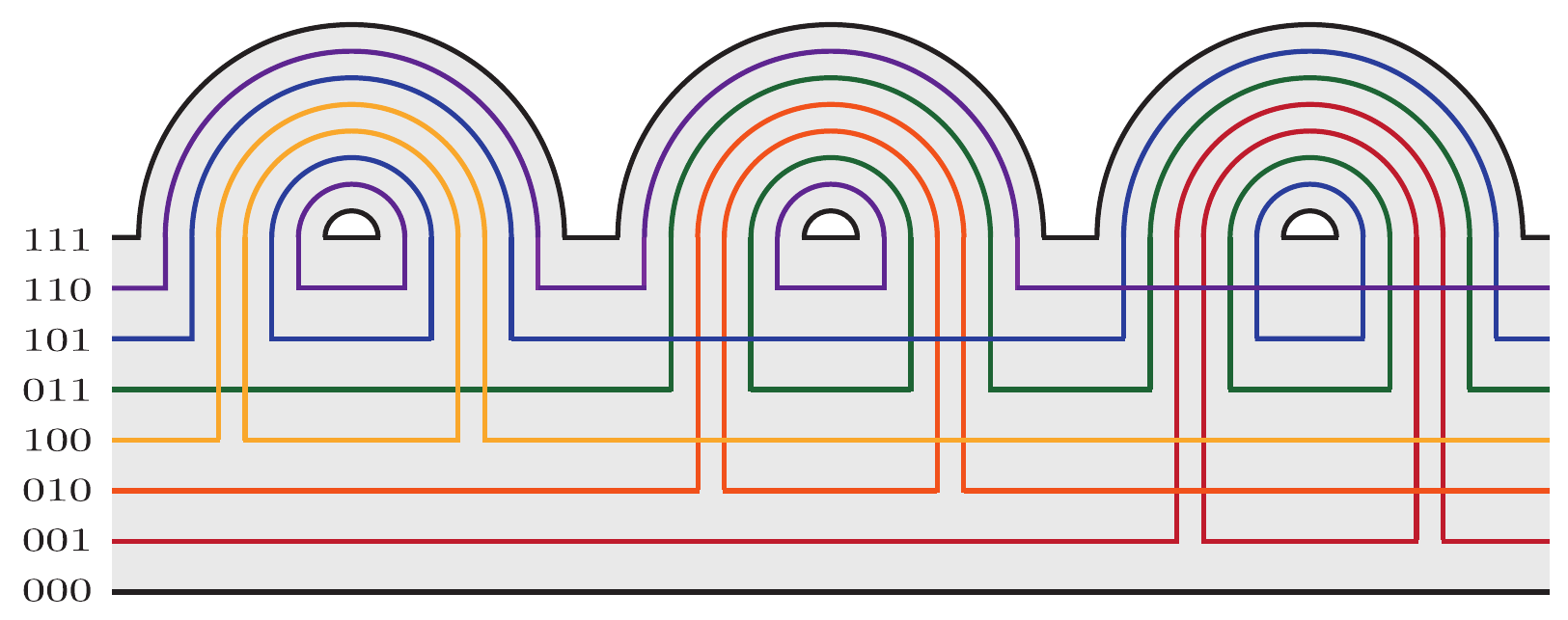}
\caption[The cobordism for the hypercube $\{0,1\}^3$.]{Half-dimensional diagram of the cobordism $W$ for the hypercube $\{0,1\}^3$.}
\label{fig:intersections1}
\end{figure}

\begin{remark} The reader who is convinced by Figure \ref{fig:intersections1}  may safely skip the proof.
\end{remark}

\begin{proof}
We list all of the vertices as $I_0, I_1, ..., I_{2^l-1}$, first in order of increasing weight and then numerically within each weight class.  We express the full cobordism as  $$W = [0,2^l-1]  \times Y \ \ \bigcup_{i=1}^l \ \ h_i$$ and embed $Y_\Ii$ and $Y_\If$ as the boundary.  We then embed the interior hypersurfaces as follows.  For $1 \leq q\leq 2^l-2$, define a slimmer 2-handle $h^q_i$ as the image of $D^2 \times D_q^2$ in $h_i$, where $D_q^2$ is the disk of radius $\frac{q}{2^l}$.  Let $\nu_q(K_i)$ be the region to which $h_i^q$ is attached, considered as a subset of $Y$.  Then we may regard $$\tilde{h}_i^q = [q,2^l-1] \times \nu_q(K_i)  \bigcup_{\{2^l-1\} \times \nu_q(K_i)}  h^q_i$$ as a longer 2-handle which tunnels through $[q,2^l-1] \times Y$ in order to attach to $ [0,q] \times Y$ along $\{q\} \times v_q(K_i)$.  In this way, we embed $W_{\Ii I_q}$ in $W$ as $$ W_{\Ii I_q} = [0,q] \times Y \bigcup_{\{i \,|\, m_i = 1\}} \tilde{h}_i^q$$
and $Y_{I_q}$ as a component of the boundary.

Now consider two vertices $I_q = (m_1,...,m_l)$ and $I_{q^\prime} = (m^\prime_1,...,m^\prime_l)$ and assume without loss of generality that $q < q^\prime$.  By construction, $Y_{I_q} \bigcap Y_{I_{q^\prime}}$ is confined to the union of the thickened attaching regions $[q,q^\prime] \times \nu(K_i)$ in $ [q,q^\prime] \times Y$ with $m_i = 1$.  If $m^\prime_i = 1$ as well, then $\tilde{h}_i^q$ is contained in the interior of $W_{\Ii q^\prime}$.  On the other hand, if $m_i > m^\prime_i$ then $\tilde{h}_i^q$ and $\partial W_{\Ii q^\prime}$ intersect in the solid torus $\{q^\prime\} \times \nu_q(K_i)$.  It follows that $Y_{I_q}$ and $Y_{I_{q^\prime}}$ intersect in one torus for each $i$ such that $m_i > m^\prime_i$.  With $q < q^\prime$, the number of such $i$ is exactly $\rho(I_q, I_q^{\prime})$, verifying (ii).  The first part of (i) immediately follows, since a subset of $\{0,1\}^l$ forms a chain if and only if $\rho$ vanishes on every pair of vertices in the subset.  In this case, $W$ decomposes as claimed by construction.
\end{proof}

We are now ready to build a special family of metrics on the cobordism $W$, starting from an initial Riemannian metric $g_0$ which is cylindrical near every $Y_I$.  Fix a path $\gamma$ from $\Ii$ to $\If$.   By Lemma \ref{lem:hyper}, $\gamma$ corresponds to a maximal subset of {\em disjoint} internal hypersurfaces $Y_{I_1} < Y_{I_2} < ... < Y_{I_{l-1}}$ in $W$.  So for each point $(T_1,\dots,T_{l-1})\in [0,\infty)^{l-1}$, we may insert necks to express $W$ as the Riemannian cobordism $W_\gamma(T_1,\dots, T_{l-1})$ given by
\begin{align}
\label{eqn:stretch}
W_{\Ii I_1} \ \bigcup_{Y_{I_1}}  \ \left([-T_1,T_1] \times Y_{I_1} \right) \  \bigcup_{Y_{I_1}} \ W_{I_1I_2}  \ \bigcup_{Y_{I_2}}  \ \cdots \ \bigcup_{Y_{I_{l-1}}}  \ \left( [-T_{l-1},T_{l-1}] \times Y_{I_{l-1}} \right)  \ \bigcup_{Y_{I_{l-1}}}  \ W_{I_{l-1}\If}.
\end{align}
We then extend this family to the cube $[0,\infty]^{l-1}$ by degenerating the metric on $Y_j$ when $T_j = \infty$.  As in the proof of the composition law, when $T_j$ grows, the $Y_{I_j}$-neck stretches, and when $T_j = \infty$, it breaks.  In particular, $W_\gamma(0,\dots,0)$ has the metric $g_0$, while $W_\gamma(\infty,\dots, \infty)$ is the disjoint union of $l$ elementary cobordisms which compose to give $W$ with metric $g_0$.

In this way, we obtain $l!$ families of metrics on $W$, each parameterized by a cube $C_\gamma$.  The facets of each cube fall evenly into two types.  A facet is {\em interior} if it is specified by fixing a coordinate at 0, and {\em exterior} if it is specified by fixing a coordinate at $\infty$.  Note that each ``almost maximal'' chain $Y_{I_1} < \dots < \widehat{Y_{I_j}} < \dots < Y_{I_{l-1}}$ can be completed to a maximal chain in exactly two ways.  It follows that each internal facet has a twin on another cube, in the sense that the twins parameterize identical families of metrics on $W$.  By gluing the cubes together along twin facets, we can build a single family of metrics which interpolates between the various ways of expressing $W$ as a composite cobordism.  In fact, this construction realizes the cubical subdivision of the following ubiquitous convex polytope.

\begin{figure}[htp]
\centering
\includegraphics[width=70mm]{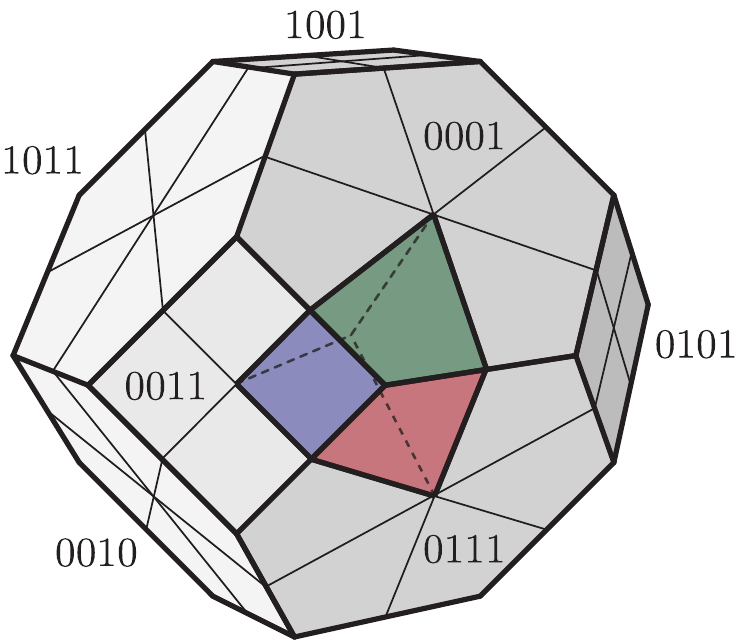}
\caption[The permutohedron of metrics for $\{0,1\}^4$.]{The cubical subdivision of the permutohedron $P_4$ consists of 24 cubes, corresponding to the $4!$ paths from $0000$ to $1111$ in $\{0,1\}^4$.  Above, the cube corresponding to the path $0000 < 0001 < 0011 < 0111 < 1111$ is shown with its exterior faces in translucent color.  Each cube shares one vertex with $P_4$ and has one vertex at the center.}
\label{fig:permutohedron}
\end{figure}

The permutohedron $P_l$ of order $l$ arises as the convex hull of all points in $\mathbb{R}^l$ whose coordinates are a permutation of $(1,2,3,\dots, l)$.   These points lie in general position in the hyperplane $x_1 + \cdots + x_l = \frac{l(l-1)}{2}$, so $P_l$ has dimension $l-1$.   The first four permutohedra are the point, interval, hexagon, and truncated octohedron (see Figure \ref{fig:permutohedron}).  The 1-skeleton of $P_l$ is the Cayley graph of the standard presentation of the symmetric group on $l$ letters:
$$S_l = \langle\sigma_1,\cdots,\sigma_{l-1} \, | \, \sigma_i^2 = 1, \, \sigma_i\sigma_{i+1}\sigma_i = \sigma_{i+1}\sigma_i \sigma_{i+1}, \, \sigma_i\sigma_j = \sigma_j \sigma_i \text{ for } \abs{i-j}>1\rangle.$$
More generally, the $(l-d)$-dimensional faces of $P_l$ correspond to partitions of the set $\{1,\dots, l\}$ into an ordered $d$-tuple of subsets $(A_1,\dots,A_d)$.  Inclusion of faces corresponds to merging of neighboring $A_j$.

The connection between the permutohedron and the hypercube rests on a simple observation: the face poset of $P_l$ is dual to the poset of chains of internal vertices in the hypercube $\{0,1\}^l$.  Namely, to each face $(A_1,\dots,A_d)$, we assign the chain $I_1 < \cdots < I_{d-1}$, where $I_j$ has $i^\text{th}$ coordinate 1 if and only if $i \in A_1 \bigcup \cdots \bigcup A_{j-1}$.  For example, in the case of the edges of the hexagon $P_3$, the correspondence is given by:
\begin{align*}
\begin{array}{cccccc}
(\{3\}, \{1,2\}) & (\{2,3\}, \{1\}) & (\{2\}, \{1,3\}) & (\{1, 2\}, \{3\}) & (\{1\}, \{2,3\}) & (\{1,3\}, \{2\}) \\
001 & 011 & 010 & 110 & 100 & 101
\end{array}
\end{align*}
In particular, each path $\gamma$ from $\Ii$ to $\If$ corresponds to a vertex $V_\gamma$ of $P_l$.

\begin{figure}[htp]
\centering
\includegraphics[width=155mm]{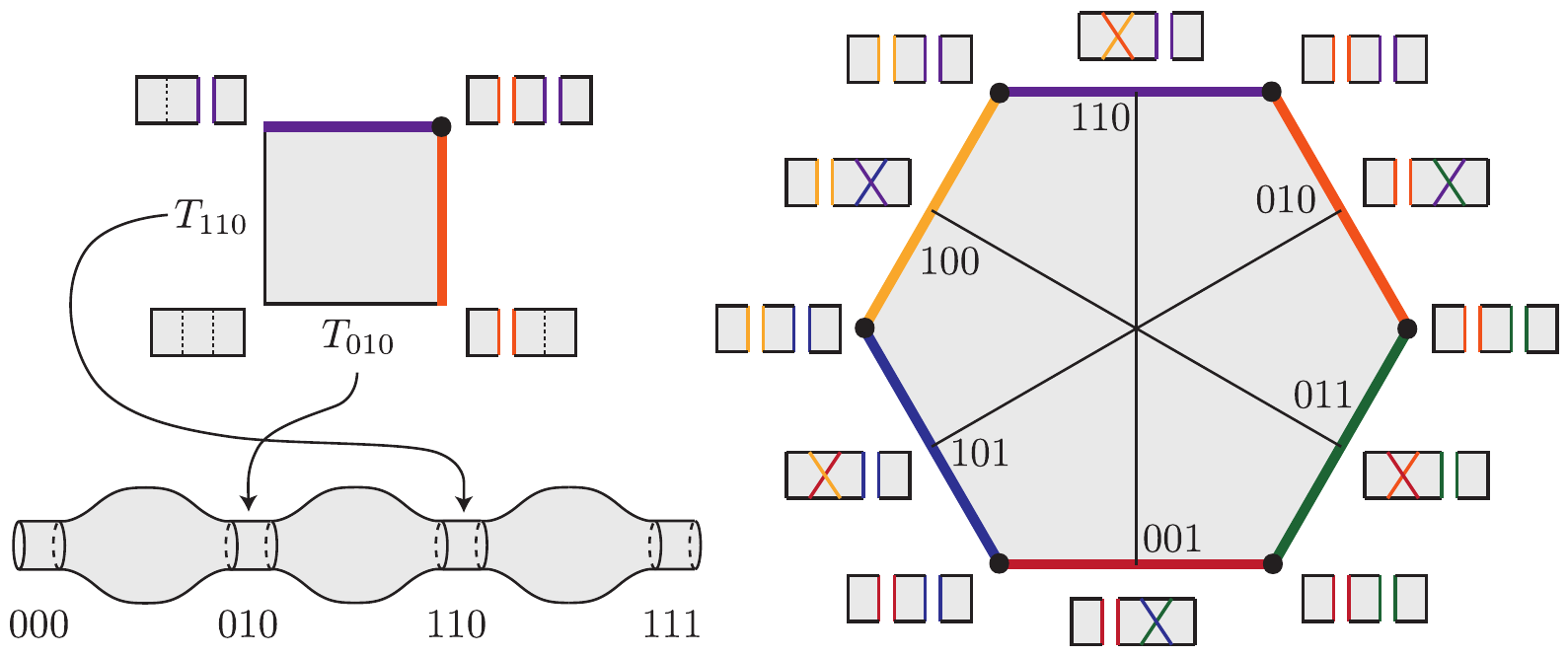}
\caption[The hexagon of metrics for $\{0,1\}^3$.]{At left, we consider the path $\gamma$ given by $000 < 010 < 110 < 111$ in $\{0,1\}^3$.  The corresponding square $C_\gamma$ with coordinates $(T_{010}, T_{011})$ parameterizes a family of metrics on the cobordism $W^*$ which stretches at $Y_{010}$ and $Y_{110}$.  We have one square for each non-intersecting pair of hypersurfaces in Figure \ref{fig:intersections1}.  These six squares fit together to form the hexagon $P_3$ at right.  The small figures at the vertices and edges illustrate the metric degenerations on $W$, read as composite cobordisms from left to right.}
\label{fig:hexagon}
\end{figure}

Now in the cubical subdivision of $P_l$, we may identify the cube containing $V_\gamma$ with the cube of metrics $C_\gamma$ so that twin interior facets are identified (see Figures \ref{fig:hexagon} and \ref{fig:permutohedron}).  In this way, the interior of $P_l$ parameterizes a family of non-degenerate metrics on $W$, while the boundary parameterizes a family of degenerate metrics.  The parameterization can be made smooth on the interior by a slight adjustment of the rate of stretching.  We summarize these observations in the following proposition.

\begin{proposition}
\label{prop:perm}
The face poset of the permutohedron $P_l$ is dual to the poset of chains of internal hypersurfaces in $W$.  In particular, the facets of $P_l$ correspond to the ways of breaking $W$ into a composite cobordism along a single interior hypersurface.  The interior of $P_l$ smoothly parameterizes a family of non-degenerate metrics on $W$, which extends naturally to the boundary in such a way that the interior of each face parameterizes those metrics which are degenerate on precisely the corresponding chain.
\end{proposition}

\begin{remark}
\label{rem:dualgraph}
We describe an alternative view of the above construction which is not essential, but will be helpful in Section \ref{sec:graph} when we consider more general lattices than the hypercube.   Consider the directed graph $\Gamma$ associated to $\{0,1\}^l$, with an edge from $I$ to $J$ whenever $J$ is an immediate successor of $I$.  Let $\bar{\Gamma}$ be the transitive closure of $\Gamma - \{\Ii, \If\}$.  The nodes of $\bar{\Gamma}$ correspond to internal hypersurfaces, and by Lemma \ref{lem:hyper}, two nodes are joined by an edge if and only if the corresponding internal hypersurfaces are disjoint.  In fact, $\bar{\Gamma}$ is the 1-sleleton of a simplical complex $\mathcal{C}_l$, whose face poset is isomorphic to the poset of non-empty cliques in $\bar{\Gamma}_l$ under inclusion.  Then $\mathcal{C}_l$ is dual to the boundary of $P_{l}$.
\end{remark}


\section{The composition law}
\label{sec:monopole}

To set notation and motivate the constructions in Section \ref{sec:linksurg}, we recall the formal properties of monopole Floer theory and the proof of the composition law, following \cite{km} (see also \cite{kmos} for an efficient survey).  We will always work over the 2-element field $\ftwo$.  Let COB be the category whose objects are compact, connected, oriented 3-manifolds and whose morphisms are isomorphism classes of connected cobordisms.  Then the monopole Floer homology groups define covariant functors from the oriented cobordism category COB to the category MOD${}_\dagger$ of modules for the ring $\ftwo[[U_\dagger]]$:
\begin{align*}
\hmb: & \, \text{COB} \to \text{MOD}_\dagger \\
\hmfrom: & \, \text{COB} \to \text{MOD}_\dagger \\
\hmbar: & \, \text{COB} \to \text{MOD}_\dagger \, .
\end{align*}
The module structure may be extended over the exterior algebra $\Lambda^*(H_1(Y) / \text{torsion})$.  These modules have a canonical mod 2 grading, and fit into a long exact sequence
\begin{align*}
\cdots \longrightarrow \hmb(Y) \longrightarrow \hmfrom(Y) \longrightarrow \hmbar(Y) \longrightarrow \cdots
\end{align*}
which is natural with respect to the maps induced by cobordisms.  The map $$\hmb(W): \hmb(Y_0) \to \hmb(Y_2)$$ induced by a cobordism $W: Y_0 \to Y_2$ satisfies the composition law
\begin{align}
\label{eqn:complaw}
\hmb(W) = \hmb(W_2) \circ \hmb(W_1)
\end{align}
whenever $W$ is the composition of cobordisms $W_1$ and $W_2$.  The composition law follows from a ``stretching the neck'' argument, as do many of the results in this paper, so we now take a moment to review the proof (see Proposition 4.16 of \cite{kmos} for details over $\ftwo$, and Proposition 26.1.2 of \cite{km} for details over $\mathbb{Z}$).

We refer the reader to \cite{km} for the full construction of the monopole Floer groups.  We first summarize the construction of the chain map $\mw: \cm(Y_0) \to \cm(Y_1)$ which induces $\hmb(W)$.  Here the monopole Floer complex $\cm(Y_i)$ is the $\ftwo$-vector space over the basis $e_\amf$ indexed by (irreducible or boundary stable) monopoles $\amf \in \check{\Cmf}(Y_i) = \Comf(Y_i) \bigcup \Csmf(Y_i)$.  Given a cobordism $W : Y_0 \to Y_1$ equipped with a metric and perturbation which are cylindrical near the boundary, we denote by $W^*$ the Riemannian manifold built by attaching the infinite cylinders $\mathbb{R} \times Y_i$ to each end of $W$. For monopoles $\amf \in \check{\Cmf}(\yi)$ and $\bmf \in \check{\Cmf}(\yj)$, and a relative homotopy class $z$ from $\amf$ to $\bmf$ in the configuration space $\mathcal{B}^\sigma(W)$, we consider the moduli space $\Modz$ of trajectories (mod gauge) on $W^*$ asymptotic to $\amf$ and $\bmf$ and in class $z$.  The map $\mw$ is defined to count isolated trajectories in such moduli spaces.  In particular, when $\amf$ is irreducible, the coefficient of $e_\bmf$ in $\mw(e_\amf)$ is the number of trajectories in $\Modz$, summed over all $z$ such that $\Modz$ is 0-dimensional.  When $\Modz$ is 1-dimensional, it has a compactification $\Modzplus$ formed by considering broken trajectories as well.  The composite maps $\ddd \mw$ and $\mw \ddd$ then count the (even) number of boundary points, so
\begin{align*}
\ddd \mw + \mw \ddd = 0,
\end{align*}
and we conclude that $\mw$ is a chain map.

More generally, suppose we have a family of metrics on $W$, smoothly parameterized by a closed manifold $P$.  The map $\mw_P: \cm(Y_0) \to \cm(Y_1)$ is defined to count isolated trajectories in the union
\begin{align}
\label{eqn:fiberprod}
\Mod_P = \bigcup_z \Modz_P
\end{align}
of fiber products
\begin{align}
\Modz_P = \bigcup_{p\in P} M_z(\amf, W(p)^*, \bmf),
\end{align}
where $W(p)$ denotes $W$ with the metric over $p$.   The compact fiber product $\Modzplus_P$ is defined similarly.  By counting boundary points of $\Modz_P$, we again conclude
\begin{align*}
\ddd \mw_P + \mw_P \ddd = 0.
\end{align*}

On the other hand, if $P$ is a compact manifold with boundary $Q$, then $\mw_P$ is no longer a chain map, because the boundary of $\Modz_P$ now includes the fibers over $Q$.  Including these contributions, we have
\begin{align}
\label{eqn:qnull}
\ddd \mw_P + \mw_P \ddd = \mw_Q.
\end{align}
Thus, $\mw_Q$ is null-homotopic and $\mw_P$ provides the chain homotopy.  That $\hmb(Y)$ is independent of the choice of metric and perturbation follows by letting $P$ be the interval $[0,1]$ parameterizing a path between two such choices.

Now let $W: Y_0 \to Y_2$ be a composite cobordism $$W: Y_0 \xrightarrow{W_1} Y_1 \xrightarrow{W_2} Y_2$$ and fix a metric on $W$ which is cylindrical near each $Y_i$.  For each $T \geq 0$, we construct a new Riemannian cobordism $W(T)$ by cutting $W$ along $Y_1$ and splicing in the cylinder $[-T,T] \times Y_1$ with the cylindrical metric.
We also define $W(\infty)$ as the disjoint union $W_1 \coprod W_2$.  In this way, $P = [0,\infty]$ parameterizes a family of metrics on $W$, where the metric degenerates on $Y_1$ at infinity.  In other words, as $T$ increases, the cylindrical neck stretches, and when $T = \infty$, it breaks.  

We again define $\mw_P$ to count isolated trajectories in the fiber products $\Modz_P$ of \eqref{eqn:fiberprod}, where now
\begin{align}
\label{eqn:inftyfiber}
M_z(\amf, W(\infty)^*, \bmf) = \bigcup_{\cmf \in \cm(Y_1)}  \bigcup_{z_1, z_2} M_{z_1}(\amf, W_1^*, \cmf) \times M_{z_2}(\cmf, W_2^*, \bmf),
\end{align}
and the inner union is taken over homotopy classes $z_1$ and $z_2$ which concatenate to give $z$.  The compact fiber product $\Modplus_P$ is defined similarly.  By counting boundary points, we conclude
\begin{align}
\label{eqn:complaw2}
\ddd \mw_P + \mw_P \ddd =  \mw + \check{m}(W_2) \check{m}(W_1).
\end{align}
Here  $\mw$ and $\check{m}(W_2) \check{m}(W_1)$ count trajectories in the fibers over $0$ and $\infty$, respectively.  Viewing $\mw_P$ as a chain homotopy, the composition law now follows.  Note that, while formally similar, \eqref{eqn:qnull} does not imply \eqref{eqn:complaw2} because the latter involves a degenerate metric.  The key technical machinery behind this generalization consists of compactness and gluing theorems for moduli spaces on cobordisms with cylindrical ends, as developed in \cite{km} and \cite{kmos}.  Our workhorse version is Lemma \ref{lem:ends} in the following section.


\section{The link surgery spectral sequence: construction}
\label{sec:linksurg}

Let $W$ be the cobordism associated to surgery on a framed link $L \subset Y$.  In Section \ref{sec:surface}, we constructed a family of metrics on $W$, parameterized by a permutohedron $P_l$ and degenerate on the boundary $Q_l$.  We now use such families to define maps between summands in a hypercube complex $X$ associated to the framed link.  That these maps define a differential will follow from a generalization of \eqref{eqn:qnull} similar in spirit to \eqref{eqn:complaw2}.  The link surgery spectral sequence is then induced by the filtration on the hypercube complex given by vertex weight.

Fix a regular metric and perturbation on the cobordism $W$ which are cylindrical near every hypersurface $\yi$.  Let $X$ be the direct sum of the monopole Floer complexes of the hypersurfaces, considered as a vector space over $\ftwo$:
$$X = \bigoplus_{I\in \{0,1\}^l} \cm(\yi)$$
We will define a differential $\dd: X \to X$ as the sum of maps $\dij: \cm(\yi) \to \cm(\yj)$ over all $I \leq J$, with $\dii$ the differential on the monopole Floer complex $\cm(\yi)$.  We now construct the maps $\dij$ when $I < J$.

Fix vertices $I < J$ and let $k = w(J) - w(I)$.  Regarding $\wij$ as the cobordism arising by surgery on a k-component, framed link in $Y_I$, with initial metric induced by $W$, we apply Proposition \ref{prop:perm} to obtain a family of metrics on $\wij$ parameterized by the permutohedron $\pij$ of dimension $k-1$.  Consider a pair of critical points $\amf \in \Cmf(\yi)$ and $\bmf \in \Cmf(\yj)$, and a relative homotopy class $z$ from $\amf$ to $\bmf$ in the configuration space $\bswij$.  As in \eqref{eqn:inftyfiber}, we must extend the definition of $M_z(\amf,\wij(p)^*,\bmf)$ to the degenerate metrics over the boundary of $\pij$.  If $p$ is in the interior of the face $I_1 < I_2 < \cdots < I_{q-1}$, then an element $\boldsymbol{\gamma}$ of $M_z(\amf,\wij(p)^*,\bmf)$ is a $q$-tuple
$$(\gamma_{01}, \gamma_{12}, \dots, \gamma_{q-1\,q})$$
where
\begin{align*}
\gamma_{j\,j+1} & \in M(\amf_j, W_{I_j I_{j+1}}^*(p),\amf_{j+1}) \\
\amf_0 &= \amf \\
\amf_q &= \bmf
\end{align*}
and the homotopy classes of these elements compose to give $z$.  Here, the metric on $W_{I_j I_{j+1}}(p)$ is the restriction of the metric on $W(p)$.  We then define $\Modzpij$ as the fiber product $$\Modzpij = \bigcup_{p\in P} \{p\} \times M_z(\amf,\wij(p)^*,\bmf).$$

In order to count the points in this moduli space, we define two elements of $\ftwo$ by
\begin{align}
\label{eqn:m}
\modzij &= \left\{ \begin{array}{l}
         \abs{\Modzpij} \text{ mod 2}, \hspace{10mm}  \text{if dim } \Modzpij = 0\\
        0, \hspace{50mm} \text{otherwise,}\end{array} \right. \\
\label{eqn:mb}
\modzrij &= \left\{ \begin{array}{l}
         \abs{\Modzprij} \text{ mod 2}, \hspace{6.9mm} \text{if dim } \Modzprij = 0\\
        0, \hspace{50mm} \text{otherwise.}\end{array} \right.
\end{align}

\begin{remark}
When $I = J$, we replace $\Modzpij$ in \eqref{eqn:m} by the moduli space $\breve{M}_z(\amf, \bmf)$ of unparameterized trajectories on the cylinder $\mathbb{R} \times Y$ (see the definition below).  We similarly replace $\Modzprij$ in \eqref{eqn:mb} by $\breve{M}^\text{red}_z(\amf, \bmf)$.
\end{remark}

Recall that $C^o(Y)$, $C^s(Y)$,  and $C^u(Y)$ are vector spaces over $\ftwo$, with bases $e_\amf$ indexed by the monopoles $\amf$ in $\Comf(Y)$, $\Csmf(Y)$, and $\Cumf(Y)$, respectively.  We use the above counts to construct eight linear maps $\dooij$, $\dosij$, $\duoij$, $\dusij$, $\dbssij$, $\dbsuij$, $\dbusij$, $\dbuuij$, where for example,
\begin{align}
\label{eqn:thed}
\left.\begin{array}{cc}{\displaystyle \dusij : C^u_\bullet (\yi) \to C^s_\bullet (\yj)} & {\displaystyle \hspace{1cm} \dusij e_\amf = \sum_{\bmf \in \Cumf(Y_J)} \sum_z \modzij e_\bmf \, ;} \\{\displaystyle \dbusij : C^u_\bullet (\yi) \to C^s_\bullet (\yj)} & {\displaystyle \hspace{1cm} \dbusij e_\amf = \sum_{\bmf \in \Cumf(Y_J)} \sum_z \modzrij e_\bmf \, .} \end{array}\right.
\end{align}
Note that the above two maps are distinct.  We then define $\dij : \cm(\yi) \to \cm(\yj)$ by the matrix
\begin{align}
\label{eqn:dij}
\dij = \left[\begin{array}{rr}\dooij & \sum_{I \leq K \leq J} \duokj \dbsuik \\\dosij & \dbssij + \sum_{I \leq K \leq J} \duskj \dbsuik\end{array}\right],
\end{align}
with respect to the decomposition $\cm(Y) = C^o(Y) \bigoplus C^s(Y)$.  The motivation behind this definition is explained in Appendix I, and we write our this map in full for $0 \leq l \leq 3$ in Appendix II.  Finally, as promised, we let $\dd : X \to X$ be the sum $$\dd = \sum_{I \leq J} \dij.$$

We now turn to proving that $\dd$ is a differential.  As in the proof of the composition law, the argument proceeds by constructing an appropriate compactification of $\Modzpij$ and counting boundary points.  We first consider the compactification of the space of unparameterized trajectories on $Y$, repeating nearly verbatim the definitions given in Section 16.1 of \cite{km}.  A trajectory $\gamma$ belonging to $M_z(\amf, \bmf)$ is {\em non-trivial} if it is not invariant under the action of $\mathbb{R}$ by translation on the cylinder $\mathbb{R} \times Y$.  An {\em unparameterized trajectory} is an equivalence class of non-trivial trajectories in $M_z(\amf, \bmf)$.  We write $\breve{M}_z(\amf, \bmf)$ for the space of unparameterized trajectories.  An {\em unparameterized broken trajectory} joining $\amf$ to $\bmf$ consists of the following data:
\renewcommand{\labelenumi}{$\bullet$}
\begin{enumerate}
\item an integer $n \geq 0$, the {\em number of components};
\item an $(n+1)$-tuple of critical points $\amf_0, \dots, \amf_n$ with $\amf_0 = \amf$ and $\amf_n = \bmf$, the {\em restpoints};
\item for each $i$ with $1 \leq i \leq n$, an unparameterized trajectory $\breve{\gamma}_i$ in $\breve{M}_z(\amf_{i-1}, \amf_i)$, the $i$th component of the broken trajectory.
\end{enumerate}
The {\em homotopy class} of the broken trajectory is the class of the path obtained by concatenating representatives of the classes $z_i$, or the constant path at $\amf$ if $n = 0$.  We write $\breve{M}^+_z(\amf, \bmf)$ for the space of unparameterized broken trajectories in the homotopy class $z$, and denote a typical element by $\bbg = (\bg_1,\dots,\bg_n)$.  This space is compact for the appropriate topology (see \cite{km}, Section 24.6).  Note that if $z$ is the class of the constant path at $\amf$, then $\breve{M}_z(\amf, \amf)$ is empty, while $\breve{M}_z^+(\amf, \amf)$ is a single point, a broken trajectory with no components.

We are now ready to define the compactification $M_z^+(\amf,\wij(p)^*,\bmf)$.  If $p$ is in the interior of the face $I_1 < I_2 < \cdots < I_{q-1}$, then an element $\bbg$ of $M_z^+(\amf,\wij(p)^*,\bmf)$ is a $(2q + 1)$-tuple
$$(\bbg_0, \gamma_{01}, \bbg_1, \gamma_{12}, \dots, \bbg_{q-1}, \gamma_{q-1\,q}, \bbg_q)$$
where
\begin{align*}
\bbg_j & \in \breve{M}^+(\amf_j, \amf_j) \\
\gamma_{j\,j+1} & \in M(\amf_j, W_{I_j I_{j+1}}^*(p),\amf_{j+1}) \\
\amf_0 &= \amf \\
\amf_q &= \bmf
\end{align*}
and $\bbg$ is in the homotopy class $z$.  The fiber product
\begin{align*}
\Modzpijplus = \bigcup_{p \in P} \{p\} \times M_z^+(\amf,W_{IJ}(p)^*,\bmf)
\end{align*}
is compact for the appropriate topology (see \cite{km}, Section 26.1).  We also write $\Modzqijplus$ for the restriction of $\Modzpijplus$ to the fibers over the boundary $\qij$ .  We can similarly define a compactification $\Modzprijplus$ of $\Modzprij$ by only considering reducible trajectories.  Fix a regular choice of metric and perturbation.

\begin{remark}
The intuition behind the following classification of ends comes from the model case of Morse homology for manifolds with boundary.  We encourage the interested reader to see Appendix I at this time.
\end{remark}

\begin{lemma}
\label{lem:ends}
If $\Modzpij$ is $0$-dimensional, then it is compact.  If $\Modzpij$ is $1$-dimensional and contains irreducibles, then $\Modzpijplus$ is a compact, $1$-dimensional space stratified by manifolds.  The $1$-dimensional stratum is the irreducible part of $\Modzpij$, while the $0$-dimensional stratum (the boundary) has an even number of points and consists of:
\renewcommand{\labelenumi}{(\Alph{enumi})}
\begin{enumerate}
\item trajectories with two or three components.  In the case of three components, the middle one is boundary-obstructed.
\item the reducibles locus $\Modzprij$ in the case that the moduli space contains reducibles as well (which requires $\amf$ to be boundary-unstable and $\bmf$ to be boundary-stable).
\end{enumerate}
If $\Modzprij$ is $0$-dimensional, then it is compact. If $\Modzprij$ is $1$-dimensional, then $\Modzprijplus$ is a compact, $1$-dimensional $C^0$-manifold with boundary.
The boundary has an even number of points and consists of:
\renewcommand{\labelenumi}{(C)}
\begin{enumerate}
\item trajectories with exactly two components.
\end{enumerate}
\end{lemma}
\begin{proof}
This is essentially Lemma 4.15 of \cite{kmos}, which in turn is a generalization of the gluing theorems in \cite{km} leading up to the proof of the composition law (see Corollary 21.3.2, Theorem 24.7.2, and Propositions 24.6.10, 25.1.1, and 26.1.6).
\end{proof}
\begin{remark}
\label{rem:dcase}
When $I=J$, Lemma \ref{lem:ends} holds with $\Modzpij$, $\Modzpijplus$, $\Modzprij$, and $\Modzprijplus$ replaced by $\breve{M}_z(\amf, \bmf)$, $\breve{M}^+_z(\amf, \bmf)$, $\breve{M}^\text{red}_z(\amf, \bmf)$, and $\breve{M}^{\text{red}+}_z(\amf, \bmf)$, respectively.
\end{remark}

We obtain a number of identities from the fact that these moduli spaces have an even number of boundary points.  We now bundle these identities into a single operator $\aij$, constructed by analogy with $\dij$.  Fix a pair of critical points $\amf \in \Cmf(\yi)$ and $\bmf \in \Cmf(\yj)$, and a relative homotopy class $z$ from $\amf$ to $\bmf$ in the configuration space $\bswij$.  We define two elements of $\ftwo$ by
\begin{align*}
n_z(\amf, W_{IJ}^*,\bmf)_{P_{IJ}} &= \left\{ \begin{array}{ll}
         \abs{\text{\{trajectories in (A) or (B)\}}} \text{ mod 2}, \hspace{3.7mm} \text{if dim } \Modzpij = 1\\
        0, \hspace{62mm} \text{otherwise,}\end{array} \right. \\
\bar{n}_z(\amf, W_{IJ}^*,\bmf)_{P_{IJ}} &= \left\{ \begin{array}{ll}
         \abs{\text{\{trajectories in (C)\}}} \text{ mod 2}, \hspace{15.7mm}  \text{if dim } \Modzprij = 1\\
        0,  \hspace{62mm} \text{otherwise.}\end{array} \right.
\end{align*}
\begin{remark}
When $I = J$, we again replace $\Modzpij$ and $\Modzprij$ by $\breve{M}_z(\amf, \bmf)$ and $\breve{M}^\text{red}_z(\amf, \bmf)$, respectively.
\end{remark}
\begin{remark}
Trajectories of type (A) necessarily have at least one irreducible component.  It follows that if $\Modzpij$ is 1-dimensional and does not contain irreducibles, then it can only have boundary points in strata of type (C).  So the condition ``if dim $\Modzpij = 1$'' is equivalent to the usual condition ``if dim $\Modzpij = 1$ and $\Modzpij$ contains irreducibles.''  A similar remark holds for the definition of $\modzij$.
\end{remark}

By Lemma \ref{lem:ends} and the above remark, $n_z(\amf, W_{IJ}^*,\bmf)_{P_{IJ}}$ counts the boundary points of $\Modzpij$ when it is 1-dimensional and contains irreducibles, and is zero otherwise.  Similarly, $\bar{n}_z(\amf, W_{IJ}^*,\bmf)_{P_{IJ}}$ counts the boundary points of $\Modzprij$ when it is 1-dimensional, and is zero otherwise.  Since the number of boundary points is even, we conclude:
\begin{align}
\label{eqn:ends}
n_z(\amf, W_{IJ}^*,\bmf)_{P_{IJ}} \textit{ and } \bar{n}_z(\amf, W_{IJ}^*,\bmf)_{P_{IJ}} \textit{ vanish for all choices of } \amf, \bmf\textit{, and } z.
\end{align}

We proceed by analogy with $\dij$, using $n_z(\amf, W_{IJ}^*,\bmf)_{P_{IJ}}$ to define linear maps $\aooij$, $\aosij$, $\auoij$, and $\ausij$, and $\bar{n}_z(\amf, W_{IJ}^*,\bmf)_{P_{IJ}}$ to define linear maps $\abssij$ and $\absuij$ (we will not need $\abusij$ or $\abuuij$).    Again, these maps all vanish identically by \eqref{eqn:ends}.  Each of these maps can be expressed as a sum of terms which are themselves compositions of the component maps of $\dij$.  Finally, we define the map $\aij : \cm(\yi) \to \cm(\yj)$ by the matrix
 \begin{align}
 \label{eqn:aij}
 \aij =  \left[\begin{array}{rr}
\aooij & \sum_{I \leq K \leq J} \left(\auokj\dbsuik + \duokj\absuik \right)\\
\aosij & \abssij + \sum_{I \leq K \leq J} \left(\auskj\dbsuik + \duskj\absuik \right)
\end{array}\right].
\end{align}
It follows that $\aij$ vanishes identically as well.  Note that the motivation behind the definition of $\aij$ is explained in Appendix I, with the cases $0 \leq l \leq 3$ written out in full in Appendix II.

\begin{lemma} 
\label{lem:aij}
$\aij$ is equal to the component of $\dd^2$ from $\cm(\yi)$ to $\cm(\yj)$:
$$\aij = \sum_{I \leq K \leq J} \dkj\dik.$$
\end{lemma}

\begin{proof}
We must show that corresponding matrix entries are equal, that is
\begin{align*}
\aooij &= \sum_{I\leq K \leq J} \dookj\dooik \\
&+ \sum_{I \leq K \leq M \leq J} \duomj\dbsukm\dosik \\
\aosij &= \sum_{I \leq K \leq J} \doskj\dooik \\
&+ \sum_{I \leq K \leq J} \dbsskj\dosik \\
&+ \sum_{I \leq K \leq M \leq J} \dusmj\dbsukm\dosik \\
\sum_{I \leq K \leq J} \left(\auokj\dbsuik + \duokj\absuik\right) &= \sum_{I\leq L \leq K \leq J} \dookj\duolk\dbsuil \\ 
&+ \sum_{I \leq K \leq M \leq J} \duomj\dbsukm\dbssik \\
&+  \sum_{I \leq L \leq K \leq M \leq J} \duomj\dbsukm\duslk\dbsuil \\
\abssij + \sum_{I \leq K \leq J} \left(\auskj\dbsuik + \duskj\absuik\right) &= \sum_{I \leq L \leq K \leq J} \doskj\duolk\dbsuil \\
&+ \sum_{I\leq K \leq J} \dbsskj\dbssik \\
&+ \sum_{I \leq L \leq K \leq J} \dbsskj\duslk\dbsuil \\ 
&+ \sum_{I \leq K \leq M \leq J} \dusmj\dbsukm\dbssik \\
&+ \sum_{I \leq L \leq K \leq M \leq J} \dusmj\dbsukm\duslk\dbsuil.
\end{align*}
After expanding out the $A^*_*$ and distributing, all terms on the right appear exactly once on the left by Lemma \ref{lem:ends} (the terms with four components appear only once since $\bar{D}^s_u D^u_s\bar{D}^s_u$ is not a term of $A^s_u$).  All other terms on the left are of the form $D^u_o \bar{D}^u_u \bar{D}^s_u$, $D^u_s \bar{D}^u_u \bar{D}^s_u$, or $\bar{D}^u_s\bar{D}^s_u$.   In the first case, $\duoktj\dbuukokt\dbsuiko$ is a term of both $\auokoj\dbsuiko$ and $\duoktj\absuikt$.   Similarly, $D^u_s \bar{D}^u_u \bar{D}^s_u$ occurs in $A^u_s\bar{D}^s_u$ and $D^u_s \bar{A}^s_u$, and $\bar{D}^u_s\bar{D}^s_u$ occurs in $A^u_s\bar{D}^s_u$ and $A^s_s$.  Therefore, each of the extra terms occurs twice and we have equality over $\ftwo$.
\end{proof}

\begin{remark}
\label{rem:badbreak}
An internal restpoint of $\bbg$ is called a {\em break}.  A break is {\em good} if the corresponding monopole is irreducible or boundary-stable.  A trajectory $\bbg \in M^+_z(\amf_0, W^*, \bmf_0)$ occurs in the {\em extended boundary} of a 1-dimensional stratum if $\bbg$ can be obtained by appending (possibly zero) additional components to either end of a boundary point of a 1-dimensinal moduli space $\Modzpij$ or $\Modzprij$.  In these terms, we have shown that among the trajectories counted by $\aij$, those with no good break each occur in the extended boundary of exactly two 1-dimensional strata.  The remaining trajectories each have one good break and occur in the extended boundary of exactly one 1-dimensional stratum.  In particular, $\dkj\dik$ counts those isolated trajectories which break well on $Y_K$.  This remark may also be understood from the perspective of path algebras, as explained in Appendix I.
\end{remark}

\begin{remark}
\label{rem:originalform}
A break of $\bbg = (\bbg_0, \gamma_{01}, \dots, \bbg_q)$ is {\em central} if it is {\em not} internal to $\bbg_0$ or $\bbg_q$.  Note that $\bbg$ has a central break if and only if it lies over a boundary fiber.  So we can express $\aij$ as the sum of similarly defined maps $\bij$ and $\qqij$, which count boundary points with and without a central break, respectively.   It follows from Remark \ref{rem:badbreak} that
\begin{align*}
\bij &= \dij\dii + \djj\dij \, ; \\
\qqij &= \sum_{I<K<J} \dkj\dik \, .
\end{align*}
$\bij$ may be thought of (imprecisely) as an operator associated to the interior of $\pij$, while $\qqij$ is (precisely) the operator associated to the boundary $\qij$ (in the case $l = 3$ in Figure \ref{fig:hexagon}, $\check{Q}^{000}_{111}$ is the sum of six composite operators, one for each edge of the hexagon).  We can then express $\aij = 0$ as $$\bij = \qqij,$$ which has the form $$\dij\dii + \djj\dij = \qqij.$$  This is the sense in which Lemma \ref{lem:aij} should be viewed as a generalization of \eqref{eqn:qnull}.  As in that case, $\qqij$ is null-homotopic and $\dij$ provides the chain homotopy.  In Appendix II, we have written out the operator $\qq$ in full in the case $l = 2$.
\end{remark}

We now conclude:

\begin{proposition}
\label{prop:dd}
$(X, \dd, F)$ is a filtered chain complex, where $F$ is the filtration induced by weight, namely
\begin{align*}
F^i X = \bigoplus_{\substack{I\in \{0,1\}^l \\ w(I) \geq i}} \cm(\yi).
\end{align*}
\end{proposition}

\begin{proof}
The equation $\dd^2 = 0$ holds by Lemma \ref{lem:aij} and the fact that the operators $\aij$ all vanish identically. The differential $\dd$ respects the filtration, as $I \leq J$ implies $w(I) \leq w(J)$.
\end{proof}

In order to describe $H_*(X,\dd)$, we recall some topology.  Let $Y_0$ be a closed, oriented 3-manifold, equipped with an oriented, framed knot $K_0$, and let $Y_1$ be the result of surgery on $K_0$ (this surgery is insensitive to the orientation of $K_0$).  $Y_1$ comes equipped with a canonical oriented, framed knot $K_1$, obtained as the boundary of the cocore of the 2-handle in the associated elementary cobordism, and given the -1 framing with respect to the cocore (see Section 42.1 of \cite{km} for details).  So we may iterate this surgery process, yielding a sequence of pairs $\{(Y_n, K_n)\}_{n \geq 0}$.  It is well-known that this sequence is 3-periodic, in the sense that for each $i \geq 0$, there is an orientation-preserving diffeomorphism $$(Y_{i+3},K_{i+3}) \xrightarrow{\cong} (Y_i,K_i)$$ which carries the oriented, framed knot $K_{i+3}$ to $K_i$.  Applying this construction to each component of the link $L \subset Y$, we may extend our collection of surgered 3-manifolds $\yi$ from the hypercube $\{0,1\}^l$ to the lattice $\{0,1,\infty\}^l$.  We may now state the 2-handle version of the link surgery spectral sequence, which computes $H_*(X,\dd)$ in stages.

\begin{theorem}
\label{thm:linksurg}
Let $Y$ be a closed, oriented 3-manifold, equipped with an $l$-component framed link $L$.  Then the filtered complex $(X,\dd,F)$ induces a spectral sequence with $E^1$-term given by
$$E^1 = \bigoplus_{I \in \{0,1\}^l} \hmb(\yi)$$
 and $d^1$ differential given by
 $$d^1 = \bigoplus_{w(J) - w(I) = 1} \hmb(\wij).$$  The link surgery spectral sequence collapses by stage $l+1$ to $\hm(Y_\infty)$.  Each page has an integer grading $\check{t}$ induced by vertex weight, which the differential $d^k$ increases by $k$.
\end{theorem}

\begin{remark}
\label{rem:bdcview}
The above statement uses different notation than that given in Theorem \ref{thm:linksurg1} in the introduction and in Theorem 4.1 of \cite{osz12}, emphasizing 2-handle addition over surgery.  To reconcile the two forms, we describe the 3-periodicity above in the case of a knot $K_0 \subset Y$ from the surgery perspective (see Section 42.1 of \cite{km}).  The complements $Y- \nu(K_n)$ are all diffeomorphic, so we may view each of the surgered manifolds $Y_n$ as obtained by gluing a solid torus to the fixed complement $Y_0 - \nu(K_0)$.  If we denote the meridian and framing of $K_n$ by $\mu_n$ and $\lambda_n$, respectively, thought of as curves on the torus $\partial\nu(K_1)$, then we have the relations
\begin{align*}
\mu_{n+1} & = \lambda_n \\
\lambda_{n+1} & = -\mu_n - \lambda_n
\end{align*}
which correspond to the matrix
$$\left[\begin{array}{cc}0 & -1 \\1 & -1\end{array}\right]$$
of order 3.   Since the framing is insensitive to the orientation of the curve, we can regard $K_0$, $K_1$, and $K_2 = K_\infty$ as having the framings $\lambda_0$, $\lambda_0 + \mu_0$, and $\mu_0$, respectively.  Therefore, $Y_I$ is shifted one step from $Y(I)$, i.e. $Y_1 =  Y(0)$, $Y_\infty = Y(1)$, and $Y_0 = Y(\infty)$.  So Theorem \ref{thm:linksurg1} is simply Theorem \ref{thm:linksurg} applied to $K_1 \subset Y_1$.  In the case of a link, the same shift in the 3-periodic sequence occurs in each component.
\end{remark}

All but the final claim of Theorem \ref{thm:linksurg} follow immediately from the usual construction of the spectral sequence associated to a filtered complex, in this case $(X, \dd)$.  The $\check t$ grading is well-defined since each differential $d^k$ is homogenous with respect to vertex weight.  We complete the proof in two stages.  First, in Section \ref{sec:graph}, we will define a complex $(\widetilde{X},\dd)$, modeled on the lattice $\{0,1,\infty\}^l$, in which $(X,\dd)$ sits as a subcomplex.  Then, in Section \ref{sec:surgexact}, we use the surgery exact triangle to conclude that $\widetilde{X}$ is null-homotopic.  The identity of the $E^\infty$ term quickly follows.


\section{Product lattices and graph associahedra}
\label{sec:graph}

Consider the lattice $\{0,1,\infty\}^l$, with the product order induced by the convention $0 < 1 < \infty$.  An $\infty$ digit contributes two to the weight.  We will sometimes also use $\infty$ to denote the final vertex $\{\infty\}^l$, with the meaning clear from context.  Consider the full cobordism $W$ from $Y_\Ii$ to $Y_\infty$, the result of attaching two rounds of $l$ 2-handles:
 $$W = \left([0,1] \times Y \ \ \bigcup_{i=1}^l h_i \right) \bigcup_{j=1}^l g_j.$$
Here $h_i$ is attached to the component $K_i$ of $L \subset Y$, and  $g_j$ is attached to $K'_j \subset Y_\If$, where $K'_j$ denotes the boundary of the co-core of $h_j$ with -1 framing.  A valid order of attachment corresponds to a maximal chain in $\{0,1,\infty\}^l$, or equivalently to a path in $\Gamma$ from $\Ii$ to $\infty$, of which there are $\frac{(2l)!}{2^l}$.  For each vertex $I = (m_1, \dots, m_l)$, we have the hypersurface $Y_I$, diffeomorphic to a boundary component of $$W_{\Ii I} = [0,1] \times Y \bigcup_{\{i \,|\, m_i \geq 1\}} h_i \bigcup_{\{i \,|\, m_i = \infty\}} g_i.$$  An $\infty$ digit corresponds to attaching a stack of two 2-handles to a component of $L \subset Y$.

As in Section \ref{sec:surface}, we will construct a polytope of metrics $\pij$ on the cobordism $\wij$ for all pairs of vertices $I<J$.  The simplest new case occurs when $l=1$, $I = 0$, and $J = \infty$.  Since $w(\infty) - w(0) = 2$, the polytope $P_{0\infty}$ should be a closed interval with degenerate metrics over the two boundary points.  However, we now have only one interior hypersurface, $Y_1$, on which to degenerate the metric.  The solution, as in \cite{kmos}, is to construct an auxiliary hypersurface $S_1$ as follows.  Let $E_1$ be the 2-sphere formed by gluing the cocore of $h_1$ to the core of $g_1$ along their common boundary $K'_1$.  Due to the -1-framing on $K'_1$, $\nu(E_1)$ is a $D^2$-bundle of Euler class -1, with $E_1$ embedded as the zero-section.  It follows that $$\nu(E_1) \cong \overline{\mathbb{CP}}^2 - \text{int}(D^4)$$ and we define the hypersurface $S_1$ to be the bounding 3-sphere $\partial \nu(K_i)$.  $P_{0\infty}$ is then identified with the interval $[-\infty, \infty]$, with the metric degenerating on $S_1$ at $-\infty$ and $Y_1$ at $\infty$.

For the lattice $\{0,1,\infty\}^l$, we will embed $l$ auxiliary 3-spheres $S_1,\dots, S_l$ in addition to the $3^l - 2$ internal hypersurfaces.  We must then construct a family of metrics which interpolates between the $\sum_{i=1}^l \binom{l}{i} \frac{(2l-i)!}{2^{l-i}}$ ways to decompose $W$ along $2l-1$ pairwise-disjoint hypersurfaces.  As a first step, we generalize Lemma \ref{lem:hyper}.  The following proposition is motivated by a half-dimensional diagram in the spirit of Figures \ref{fig:intersections1} and Figure \ref{fig:intersections2}.

\begin{proposition}
\label{prop:hyper2}
The full set of $3^l - 2$ internal hypersurfaces $\yi$ and $l$ spheres $S_i$ can be simultaneously embedded in the interior of $W$ so that the following conditions hold:
\renewcommand{\labelenumi}{(\roman{enumi})}
\begin{enumerate}

\item The internal hypersurfaces in any subset are pairwise disjoint as submanifolds of $W$ if and only if they form a chain.  In this case, cutting along $Y_{I_1} < Y_{I_2} < ... < Y_{I_k}$ breaks $W$ into the disjoint union $$W_{\Ii I_1} \ \textstyle\coprod \ W_{I_1I_2} \ \coprod \ \cdots \ \coprod \ W_{I_k\infty}.$$

\item Distinct $\yi$ and $\yj$ intersect in exactly $\rho(I,J)$ disjoint tori.

\item $\yi$ and $S_i$ intersect if and only if $m_i = 1$, where $I = (m_1,\dots,m_l)$.  In this case, they intersect in a torus.

\item The $S_i$ are pairwise disjoint.

\end{enumerate}
\end{proposition}

\begin{proof}
List the vertices as $I_0, I_1, ..., I_{3^l}$, first in order of increasing weight and then numerically within each weight class.  We express the full cobordism as $$W =  [0,3^l] \times Y \ \ \bigcup_{i=1}^l \ h_i \ \bigcup_{i=1}^l \ g_i$$ and embed $Y_\Ii$ and $Y_\infty$ as the boundary.  As in the proof of Lemma \ref{lem:hyper}, for each $1 \leq q \leq 3^l-1$, we have slimmer 2-handles $h^q_i$ and $g^q_i$ as the images of $D^2 \times D_q^2$ in $h_i$ and $g_i$, respectively, where $D_q^2$ is the disk of radius $\frac{q}{3^n}$.  Again, we think of $$\tilde{h}_i^q = [q,3^l] \times \nu_q(K_i) \bigcup_{\{3^l\} \times \nu_q(K_i)}  h^q_i$$ as a longer 2-handle which tunnels through $[q,3^l] \times Y$ in order to attach to $[0,q] \times Y$ along $\{q\} \times \nu_q(K_i)$.  Let $K_i^\prime$ be the boundary of the cocore of $h_i$, so that $\nu_q(K^\prime_i) = D^2_q \times \partial D^2$ is the region of $h_i$ to which $g_i^q$ attaches.  Let $A_q^i$ be the annulus given in polar coordinates $(r, \theta)$ by $[\frac{q}{3^l},1] \times S^1$, thought of as sitting in the cocore of $h_i$.  The boundary of $D_q^2 \times A^q_i$ consists of $\nu_q(K^\prime_i)$ and a radial contraction of $\nu_q(K^\prime_i)$ into the interior of $h_i$, denoted $\tilde{\nu}_q(K^\prime_i)$.  So we may regard $$\tilde{g}_i^q = D_q^2 \times A^q_i \bigcup_{\nu_q(K^\prime_i)} g^q_i$$ as a longer 2-handle which tunnels through $D^2 \times A_q^i \subset h_i$ in order to attach to $\tilde{h}^q_i$ along $\tilde{\nu}_q(K^\prime_i) \subset \partial \tilde{h}^q_i$.  In this way, we embed $W_{\Ii I_q}$ in $W$ as $$W_{\Ii I_q} = Y \times [0,q] \bigcup_{\{i \,|\, m_i \geq 1\}} \tilde{h}_i^q\bigcup_{\{i \,|\, m_i = \infty\}} \tilde{g}^q_i$$
and $Y_{I_q}$ as a component of its boundary.  Here $I_q = (m_1, ..., m_l)$.  Next, let the 2-sphere $E_i$ be the result of gluing the cocore of $h_i$ and the core of $g_i$ along their common boundary $K_i^\prime$, and let $\nu(E_i)$ be the result of gluing together the corresponding trivial $D^2$-bundles of radius $\frac{1}{2 \cdot 3^l}$.  Then $\nu(E_i)$ is a $D^2$-bundle of Euler class -1, and we embed the 3-sphere $S_i$ as its boundary.

Conditions (i) and (ii) now follow from a straightforward generalization of the proof of Lemma \ref{lem:hyper}.   For (iii), note that if $m_i = 1$, then the intersection of $Y_I$ and  $S_i$ is the boundary of the restriction of the $D^2$-bundle $\nu(E_i)$ to $K_i^\prime$.  Finally, the $S_i$ are pairwise disjoint because they live in different pairs of handles.
\end{proof}

For fixed $I < J$, the interval $\{K \, |\; I \leq K \leq J \}$ takes the form $\{0,1,\infty\}^m \times \{0,1\}^k$ for some pair of non-negative integers $(m,k)$ with $m + k = l$.  In order to define the maps $\dij$ in general, we need to construct a polytope $P_{m,k}$ of dimension $2m+k-1$ for each pair $(m,k)$.  We define $P_{m,k}$ abstractly to have a face of co-dimension $d$ for every subset of $d$ mutually disjoint hypersurfaces in the interior of $W$, with inclusion of faces dual to inclusion of subsets.  Our definition is justified by Theorem \ref{thm:graph}, which realizes $P_{m,k}$ concretely as a convex polytope.

\begin{figure}[htp]
\centering
\includegraphics[width=155mm]{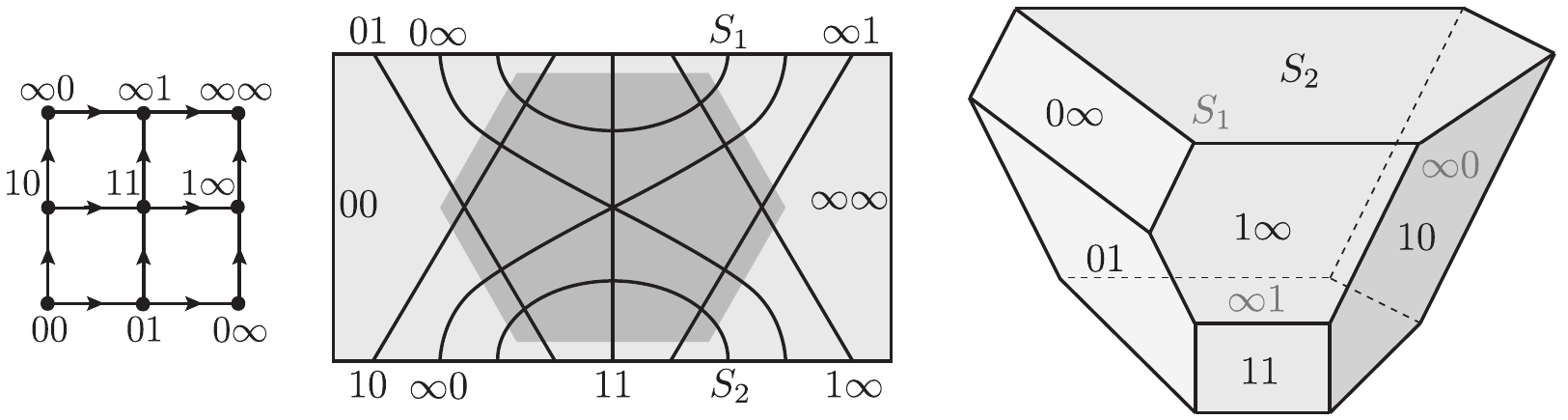}
\caption[The associahedron of metrics for $\{0,1,\infty\}^2$.]{Consider the full cobordism $W$ corresponding to the lattice $\{0,1,\infty\}^2$ at left.  The seven interior hypersurfaces and two auxiliary 3-spheres are embedded in $W$ in such a way that the diagram at center accurately depicts which pairs intersect (although the triple intersection point is an artifact).  The nine internal arcs in the diagram are arranged so that by stretching normal to disjoint subsets, we obtain a parameterization of the space of conformal structures on the hexagon, which is known to compactify to the associahedron $K_5$ at right.  This connection between associahedra and conformal structures on polygons leads to a dictionary between the techniques in this paper and their counterparts in Heegaard Floer homology, as explained in Section \ref{sec:hfproof}.}
\label{fig:assoc}
\end{figure}

In order to motivate this theorem, we first construct those $P_{m,k}$ of dimension three or less by hand.  The polytopes $P_{0,1}$, $P_{0,2}$, $P_{0,3}$, and $P_{0,4}$ are the first few permutohedra of Proposition \ref{prop:perm}, namely a point, an interval, a hexagon, and a truncated octahedron (recall Figure \ref{fig:permutohedron}).   We saw that $P_{1,0}$ is an interval, and it is easy to see that $P_{1,1}$ is the associahedron $K_4$, otherwise known as the pentagon.  $P_{2,0}$ is more interesting.  In Figure \ref{fig:assoc}, we use a trick to establish that it is $K_5$, also known as Stasheff's polytope \cite{stash}.  For a fun and informal introduction to associahedra, see \cite{cass}.  Note that $K_n$ has dimension $n-2$, while $P_n$ has dimension $n-1$.

\begin{figure}[htp]
\centering
\includegraphics[width=155mm]{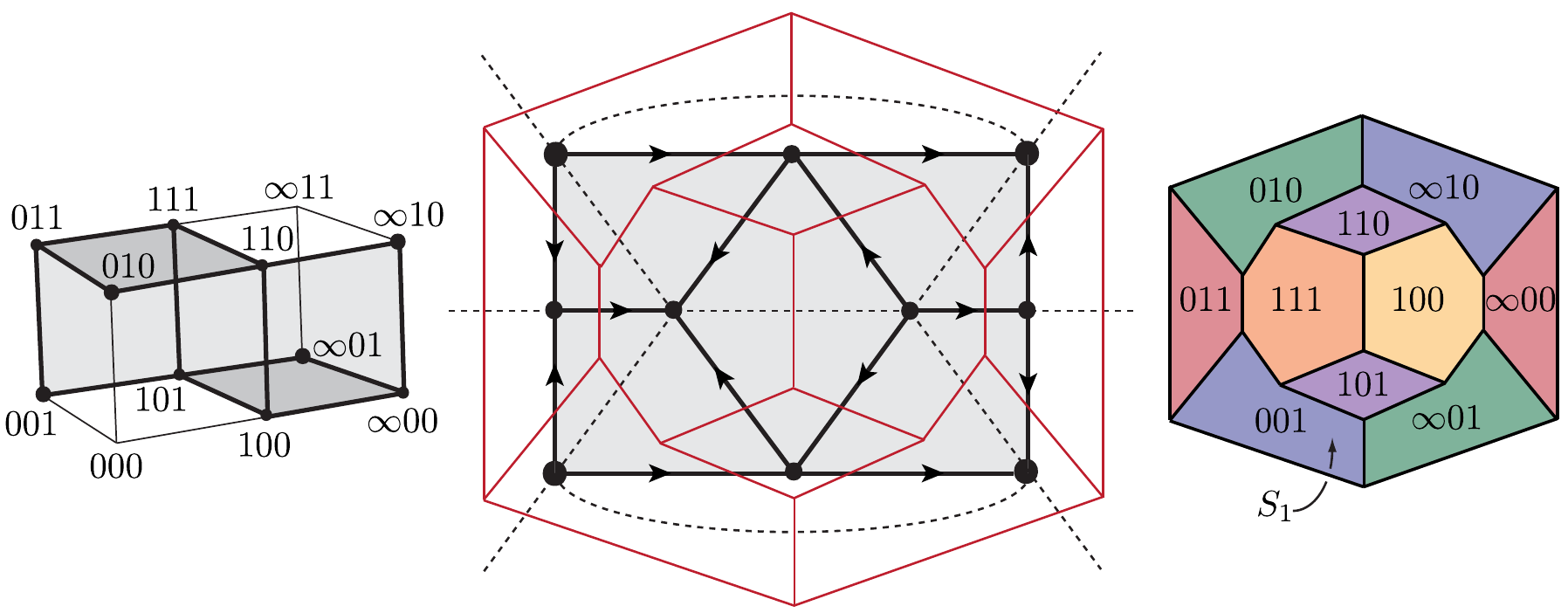}
\caption[The polyhedron of metrics for $\{0,1,\infty\} \times \{0,1\}^2$.]{We construct the boundary of the polyhedron $P_{1,2}$ as the dual of the simplicial complex $\mathcal{C}_{1,2}$.  First, at left, we remove the initial and final nodes from the lattice $\{0,1,\infty \} \times \{0,1\}^2$.  We then flatten the shaded region and take the transitive closure to obtain $\bar{\Gamma}$, represented by the shaded rectangle and compact dotted line segments at center.  Next we add the vertex $I^\prime_1$ at infinity (not shown) and connect it by dotted lines to the six nodes for which $m_1 \neq 1$.  At this stage, we have constructed $\bar{\Gamma}^\prime$, the 1-skeleton of $\mathcal{C}_{1,2}$.  The faces of $\mathcal{C}_{1,2}$ are the 3-cliques (triangles).  Drawing the dual with thin red lines, we obtain the boundary of $P_{1,2}$.  At right, we have redrawn $P_{1,2}$.  The face $S_1$ corresponds to the large hexagonal base under the colorful tortoise shell.  The 12 vertices away from $S_1$ correspond to the 12 paths through the lattice.}
\label{fig:dual}
\end{figure}

At this stage, it may be tempting to conjecture that all the $P_{m,k}$ are permutohedra or associahedra.  We check this against the only remaining 3-dimensional case, namely $P_{1,2}$.  To build this polyhedron, it is useful to return to the viewpoint of Remark \ref{rem:dualgraph}.  Let $\Gamma$ be the oriented graph corresponding to the lattice $\{0,1,\infty\}^m \times \{0,1\}^k$.   Let $\bar{\Gamma}$ be the unoriented graph obtained as the transitive closure of $\Gamma$ with its initial and final nodes removed.  We now add $l$ additional nodes $I^\prime_i$ (representing the $S_i$) to $\bar{\Gamma}$ and connect each $I^\prime_i$ to the others and to  those $I = (m_1,\dots,m_l) \in \bar{\Gamma}$ with $m_i \neq 1$.  By Proposition \ref{prop:hyper2}, the nodes of the resulting graph $\bar{\Gamma}^\prime$ are in bijection with the full set of hypersurfaces, with two nodes connected by an edge if and only if the corresponding hypersurfaces are disjoint.  The graph $\bar{\Gamma}^\prime$ is the 1-skeleton of a simplical complex $\mathcal{C}_{m,k}$ whose face poset is isomorphic to the poset of non-empty cliques in $\bar{\Gamma}^\prime$ under inclusion.  That is, the $d$-dimensional faces of $\mathcal{C}_{m,k}$ are in bijection with the $d$-cliques of $\bar{\Gamma}^\prime$ (the fact that this poset defines a simplicial complex will follow from Theorem \ref{thm:graph}).  The simple polytope dual to $\mathcal{C}_{m,k}$ is then, by definition, the boundary of $P_{m,k}$.  In Figure \ref{fig:dual}, we illustrate this process for  $P_{1,2}$, concluding that it is indeed something new.

The right hand side of Figure \ref{fig:dual} illustrates $P_{1,2}$ as a convex polytope in $\mathbb{R}^3$.  However, our dual-graph perspective does not provide such an explicit realization of $P_{m,k}$ in higher dimensions.  While searching for an alternative construction of $P_{1,2}$, the author discovered beautiful illustrations of similar polyhedra in \cite{carr} and \cite{d}.  Given a connected graph $G$ with $n$ vertices, Carr and Devadoss construct a convex polytope $P_G$ of dimension $n-1$, the {\em graph-associahedron} of $G$, using the following notions.

A {\em tube} of $G$ is a proper, non-empty set of nodes of $G$ whose induced graph is a connected subgraph of $G$.  There are three ways in which tubes $t_1$ and $t_2$ can interact:
\renewcommand{\labelenumi}{(\arabic{enumi})}
\begin{enumerate}
\item Tubes are {\em nested} if $t_1 \subset t_2$ or $t_2 \subset t_1$;
\item Tubes {\em intersect} if $t_1 \bigcap t_2 \neq \emptyset$ and $t_1 \not\subset t_2$ and $t_2 \not\subset t_1$; 
\item Tubes are {\em adjacent} if $t_1 \bigcap t_2 = \emptyset$ and $t_1 \bigcup t_2$ is a tube in $G$.
\end{enumerate}
Tubes are {\em compatible} if they do not intersect and they are not adjacent.  A {\em tubing} $T$ of $G$ is a set of tubes of $G$ such that every pair of tubes in $T$ is compatible.

\begin{figure}[htp]
\centering
\includegraphics[width=155mm]{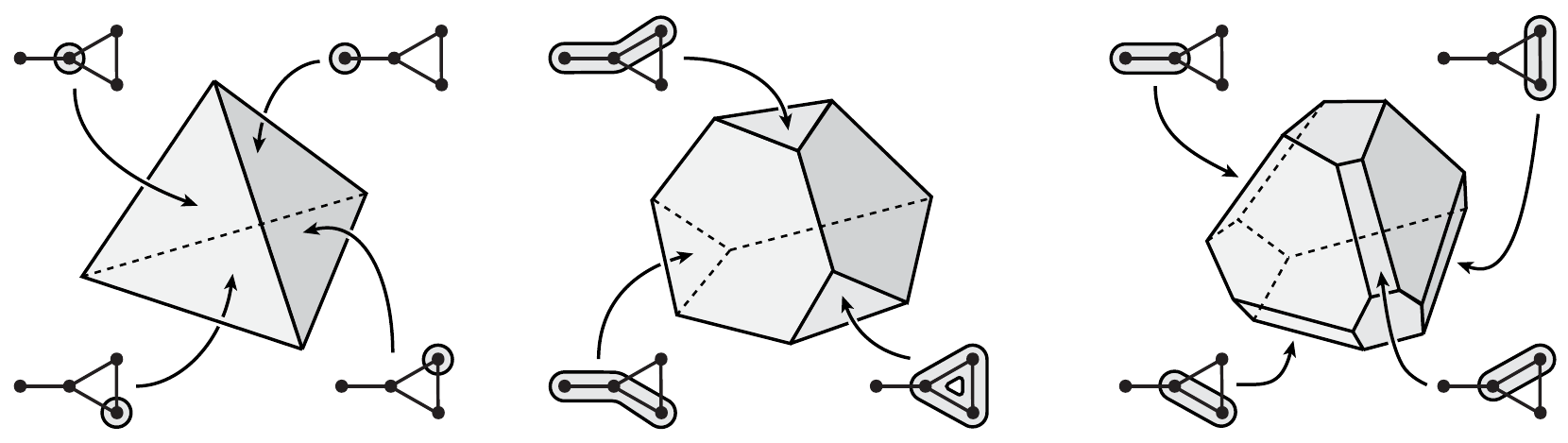}
\caption[The graph associahedron of the 3-clique with one leaf.]{We have modified Figure 6 in \cite{d} to illustrate the sculpting of $P_G$ for the graph $G$ given by the 3-clique with one leaf.  Each node of $G$ slices out a half-space in $\mathbb{R}^3$, leaving the 3-simplex $\triangle_G$ at left.  Next, we shave down those vertices of $\triangle_G$ which correspond to the connected, induced subgraphs of size three.  Finally, at right, we shave down those edges of $\triangle_G$ which correspond to the edges of $G$.  This figure also illustrates the bijection \eqref{eqn:bijection}.}
\label{fig:truncation}
\end{figure}

We now define the graph-associahedron of a connected graph $G$ with $n$ nodes.  Labelling each facet of the $n-1$ simplex $\triangle_G$ by a node of $G$, we have a bijection between the faces of $\triangle_G$ and the proper subsets of nodes of $G$.  By definition, $P_G$ is sculpted from $\triangle_G$ by truncating those faces which correspond to a connected, induced subgraph of $G$ (see Figure \ref{fig:truncation}).  We therefore have a bijection
\begin{align}
\label{eqn:bijection}
\text{\{facets of $P_G$\}} \longleftrightarrow \text{\{tubes of $G$\}}.
\end{align}
More generally, Carr and Devadoss prove that $P_G$ is a simple, convex polytope whose face poset is isomorphic to the set of valid tubings of $G$, ordered such that $ T < T^\prime$ if $T$ is obtained from $T^\prime$ by adding tubes.  Moreover, in  \cite{d}, Devadoss derives a simple, recursive formula for a set of points with integral coordinates in $\mathbb{R}^n$, whose convex hull realizes $P_G$.

\begin{remark}
Carr and Devadoss trace their construction back to the Deligne-Knudsen-Mumford compactification $\overline{\mathcal{M}}_{0,n}(\mathbb{R})$ of the real moduli space of curves.  In this context, the sculpting of $P_G$ is thought of as a sequence of real blow-ups.  When $G$ is a Coxeter graph, $P_G$ tiles the compactification of the hyperplane arrangement associated to the corresponding Coxeter system.  The $n$-node clique, path, and cycle yield the $(n-1)$-dimensional permutohedron, associahedron, and cyclohedron, respectively.
\end{remark}

Comparing Figures \ref{fig:dual} and  \ref{fig:truncation}, we see that $P_{1,2}$ is precisely the graph-associahedron of the 3-clique with one leaf (recall that the $n$-clique is the complete graph on $n$ vertices).  In fact, all of the polytopes $P_{m,k}$ are graph-associahedra:

\begin{theorem}
\label{thm:graph}
The polytope $P_{m,k}$ associated to the lattice $\{0,1,\infty\}^m \times \{0,1\}^k$ is the graph-associahedron of the $(m+k)$-clique with $m$ leaves.  More generally, the polytope naturally associated to the lattice $\{0,...,n_1\} \times \cdots \times \{0,...,n_l\}$, with all $n_i \geq 1$, is the graph-associahedron of the $l$-clique with paths of length $n_1-1,\dots,n_l-1$ attached.
\end{theorem}

\begin{figure}[htp]
\centering
\includegraphics[width=155mm]{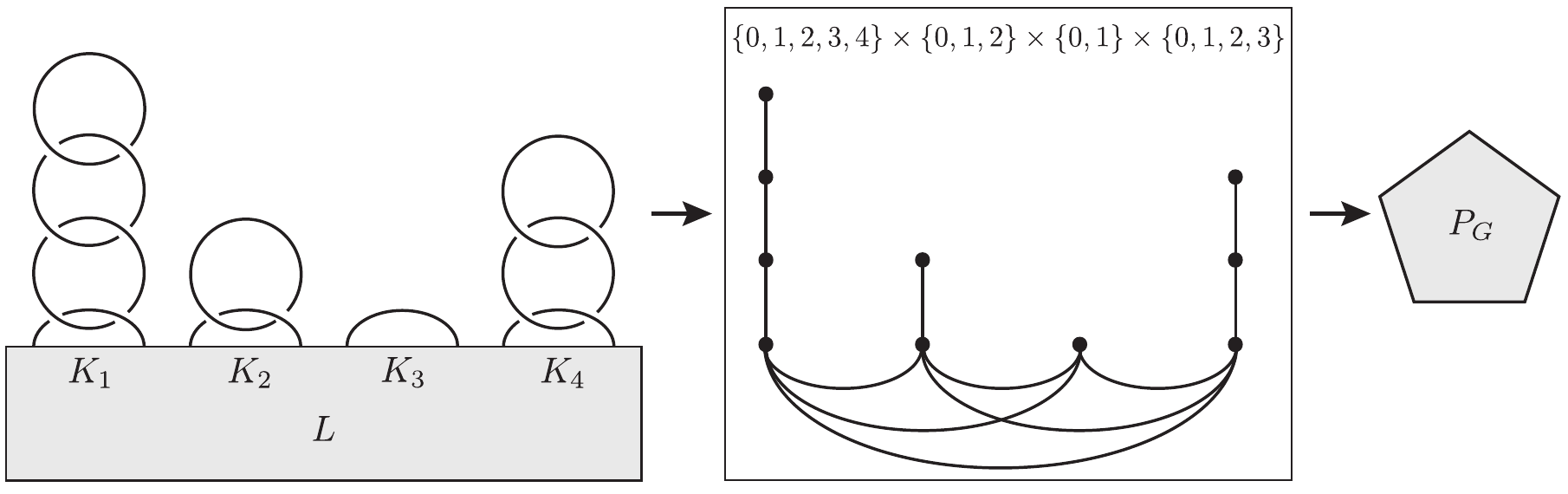}
\caption[From link surgeries to lattices, graphs, and polytopes.]{The figure at left represents a Kirby diagram arising from the 3-periodic surgery sequence applied to each component of a framed link with four components.  The corresponding lattice is the product of four chains, while the graph is obtained by appending paths to the complete graph on four vertices.  The pentagon at right represents the corresponding 9-dimensional graph associahedron.  The above assignment of a polytope $P_G$ to a finite product lattice generalizes the assignment of the permutohedron to the hypercube described in Section \ref{sec:surface}.}
\label{fig:kirbygraph}
\end{figure}

\begin{proof}
An example is given in Figure \ref{fig:kirbygraph}.  We first consider the lattice $\{0,1,\infty\}^m \times \{0,1\}^k$.  In addition to the $3^m 2^k - 2$ internal hypersurfaces $Y_I$, we have $m$ auxialliary hypersurfaces $S_i$.  Let $G$ be the complete graph on nodes $v_1,\dots,v_{m+k}$ with a leaf $v_i^\prime$ attached to $v_i$ for each $i = 1, \dots, m$.  The bijection \eqref{eqn:bijection} is given by
\begin{align*}
Y_I & \mapsto \{v_i \, | \, m_i \geq 1\} \ \bigcup \ \{v^\prime_i \, | \, m_i = \infty\} \\
S_i & \mapsto \{v_i^\prime\}.
\end{align*}
and extends to an isomorphism of posets.

Next, consider the lattice $\{0,1,...,n\}$.  The cobordism $W$ is then built by attaching a single stack of handles $h_1 \bigcup \cdots \bigcup h_n$ to $[0,1] \times Y$ (the $n = 3$ case is shown in Figure \ref{fig:intersections2}, though with different notation).  In addition to the internal hypersurfaces $Y_1,\dots, Y_{n-1}$, we include an auxiliary hypersurface $S^j_k$ between each pair of handles $(h_j, h_k)$ with $1 \leq j < k \leq n$, embedded as the boundary of a tubular neighborhood of the union of the intervening 2-spheres $E_i$.  In fact, if $k - j \equiv 2$ (mod 3), then $S^j_k$ is diffeomorphic to $S^1 \times S^2$.  Otherwise, $S^j_k$ is diffeomorphic to $S^3$.  By a straightforward variation on the theme of Lemma \ref{lem:hyper} and Proposition \ref{prop:hyper2}, these $n - 1 + \binom{n}{2}$ hypersurfaces can all be embedded in $W$ in such a way that
\renewcommand{\labelenumi}{(\roman{enumi})}
\begin{enumerate}
\item the $Y_i$ are all disjoint;
\item $Y_i$ and $S^j_k$ intersect if and only if $j \leq i < k$.  In this case, they intersect in a torus.
\item $S^{j_1}_{k_1}$ and $S^{j_2}_{k_2}$ intersect if and only if the intervals $\{j_1, \dots, k_1\}$ and $\{j_2, \dots, k_2\}$ overlap but are not nested.  In this case, they intersect in a torus.
\end{enumerate}
Now let the graph $G$ be the path with nodes $\{v_0, \dots, v_n\}$.  The bijection \eqref{eqn:bijection} is given by
\begin{align*}
Y_i & \mapsto \{v_0,...,v_i\} \\
S^j_k & \mapsto \{v_j,...,v_k\}.
\end{align*}
and extends to an isomorphism of posets.  As remarked above, $P_G$ is then the $(n-1)$-dimensional associahedron $K_{n+1}$.  The result for a lattice consisting of an arbitrary product of chains follows from a straightforward, subscript-heavy amalgamation of the arguments in the above two cases.
\end{proof}

Now consider the lattice $\Lambda = \{0,...,n_1\} \times \cdots \times \{0,...,n_l\}$ with the corresponding graph $G$ given by Theorem \ref{thm:graph}.  Using a formula in \cite{d}, we can realize $P_G$ concretely as the convex hull of vertices in general position in $\mathbb{R}^d$, where $d = n_1 \cdots n_l - 1$.  Now $P_G$ has one vertex $V_\gamma$ for every maximal collection $\gamma$ of disjoint hypersurfaces in the cobordism $W$ with initial metric $g_0$.  As in Section \ref{sec:surface}, we associate to the vertex $V_\gamma$ a cube of metrics $C_\gamma$ which stretches on the hypersurfaces in $\gamma$.  We can then use $P_\Lambda$ to parameterize a family of metrics on $W$ by identifying each $C_\gamma$ with the cube containing the vertex $V_\gamma$ in the cubical subdivision of $P_\Lambda$.  In particular, $P_{m,k}$ consists of $\sum_{i=0}^m \binom{m}{i} \frac{(2m+k-i)!}{2^{m-i}}$ cubes.

\begin{remark}
\label{rem:genlat}
Using these polytopes of metrics, we can define maps $\dij$ associated to any lattice formed as a product of chains of arbitrary length, where $\{0,\dots,n\}$ has length $n$.  However, we will see that this gives rise to a differential if and only if all the chains have length one or two.  When there is a chain of length three or more, additional terms arise from breaks on auxiliary hypersurfaces.  We will see this phenomenon explicitly for a single chain of length three in the proof of the surgery exact triangle (see Theorem \ref{thm:exact}).
\end{remark}

Having constructed polytopes of metrics for all intervals in the lattice $\{0,1,\infty\}^l$, we proceed to define the complex $(\widetilde{X}, \dd)$.  Fix a metric on the cobordism $W$ which is cylindrical near every hypersurface $\yi$ and auxiliary hypersurface $S_i$, where each $S_i$ has been equipped with the round metric.
We let
$$\widetilde{X} = \bigoplus_{I\in \{0,1, \infty \}^l} \cm(\yi)$$
and define the maps $\dij: \cm(\yi) \to \cm(\yj)$ by {\em exactly} the same construction and matrix \eqref{eqn:dij} as before, with $\dd: \widetilde{X} \to \widetilde{X}$ their sum.

We now prove that $\dd$ is a differential by an argument which parallels that in Section \ref{sec:linksurg}.  We first expand our definition of $M_z^+(\amf,\wij(p)^*,\bmf)$ to intervals of the form $\{0,1,\infty\}^m \times \{0,1\}^k$.  Let $V_i$ denote the copy of $\cptb - \text{int}(D^4)$ cut out by $S_i$.  For each $I < J$, let $S_{j_1}, ..., S_{j_{n(I,J)}}$ be the spheres completely contained in $\wij$.  We denote the corresponding cobordism with $n(I,J)+2$ boundary components by $$U_{IJ} = \wij - \bigcup_{s=1}^{n(I,J)} \text{int}\left(V_{j_s}\right).$$
If $p$ is in the interior of the face $(I_1 < I_2 < \cdots < I_{q-1}, S_1,...,S_r)$, then an element $\bbg$ of $M_z^+(\amf,\wij(p)^*,\bmf)$ is a $(2q+ 2r + 1)$-tuple
\begin{align}
\label{eqn:broken}
(\bbg_0, \gamma_{01}, \bbg_1, \gamma_{12}, \dots, \bbg_{q-1}, \gamma_{q-1\,q}, \bbg_q, \eta_1, \bbd_1, \dots,  \eta_r, \bbd_r)
\end{align}
where
\begin{align*}
\bbg_j & \in \breve{M}^+(\amf_j, \amf_j) \\
\gamma_{j\,j+1} & \in M(\amf_j, \cmf_{j_1}, \dots, \cmf_{j_{n(I_j,I_{j+1})}}, U_{I_j I_{j+1}}^*(p),\amf_{j+1})\\
\bbd_i & \in \breve{M}^+ (\cmf_i, \cmf_i) \\
\eta_i & \in M(V_i^*(p), \cmf_i) \\
\amf_0 &= \amf \\
\amf_q &= \bmf \\
\amf_j & \in \mathfrak{C}(Y_{I_j}) \\
\cmf_i & \in \mathfrak{C}(S_i)
\end{align*}
and $\bbg$ is in the homotopy class $z$ (and similarly when $p$ is in the interior of a face which includes a subset of the $S_i$ other than the first $r$).
The fiber product $$\Modzpijplus = \bigcup_{p \in P} \{p\} \times M_z^+(\amf,W_{IJ}(p)^*,\bmf)$$ is compact.  We then define $\aij: \cm(\yi) \to \cm(\yj)$ by {\em exactly} the same construction and matrix \eqref{eqn:aij} as before.

We will also need the following lemma, consolidated from \cite{kmos} (see Lemma 5.3 there and the preceding discussion).  The essential point is that there is a diffeomorphism of $\cptb - \text{int}(D^4)$ which restricts to the identity on the boundary and induces a fixed-point-free involution on the set of spin$^c$ structures.

\begin{lemma}
\label{lem:cp2}
Fix a sufficiently small perturbation on $S_i$.  Then for each $\cmf \in \mathfrak{C}(S_i)$, $\breve{M}(\cmf, \cmf) = \emptyset$ and the trajectories in the zero-dimensional strata of $M^+(V_i^*, \cmf)$ occur in pairs.
\end{lemma}

When $\Modzpij$ or $\Modzprij$  is 1-dimensional, the number of boundary points in the corresponding compactification is still even (technically, using a generalization of Lemma \ref{lem:ends} to the case of cobordisms with three boundary components, as done in \cite{km} by introducing doubly boundary obstructed trajectories). By Lemma \ref{lem:cp2}, the number of boundary points which break on precisely some non-empty, fixed collection $\{S_{i_1},...S_{i_r}\}$ of the auxiliary hypersurfaces is a multiple of $2^r$, via the pairing of $\eta_{i_j}$ and $\eta_{i_j}^\prime$ in \eqref{eqn:broken}.  Therefore, by inclusion-exclusion, the number of boundary points which do not break on any of the $S_i$ is even as well.   Since these are precisely the boundary points counted by the matrix \eqref{eqn:aij}, $\aij$ still vanishes and the proof of Lemma \ref{lem:aij} goes through without change.  We conclude:

\begin{proposition}
\label{prop:dd2}
$(\widetilde{X}, \dd, F)$ is a filtered chain complex, where $F$ is the filtration induced by weight, namely
\begin{align*}
F^i \widetilde{X} = \bigoplus_{\substack{I\in \{0,1,\infty\}^l \\ w(I) \geq i}} \cm(\yi).
\end{align*}
\end{proposition}

\begin{remark}
In Appendix II, we have written out the $l=1$ case in full.  Note also that, while we were compelled to introduce auxiliary hypersurfaces $S_i$ in order to obtain polytopes, the corresponding facets contribute vanishing terms to $\qqij$ by Lemma \ref{lem:cp2}.  We thereby recover $$\qqij = \sum_{I < K < J} \dkj\dik.$$
\end{remark}

\section{The surgery exact triangle}
\label{sec:surgexact}

We will identify the $E^\infty$ page of the link surgery spectral sequence by applying the surgery exact triangle to the complex of Proposition \ref{prop:dd2}.  Before stating the surgery exact triangle, we first recall the algebraic framework underlying its derivation in both monopole and Heegard Floer homology (see \cite{kmos} and  \cite{osz12}, respectively).

\begin{lemma}
\label{lem:alg}
Let $\{A_i\}_{i=0}^\infty$ be a collection of chain complexes and let $$\{f_i : A_i \to A_{i+1}\}_{i = 0}^\infty$$ be a collection of chain maps satisfying the following two properties:

\renewcommand{\labelenumi}{(\roman{enumi})}
\begin{enumerate}
\item $f_{i+1} \circ f_i$ is chain homotopically trivial, by a chain homotopy $$H_i : A_i \to A_{i+2}$$

\item the map $$\psi_i = f_{i+1} \circ H_i + H_{i+1} \circ f_i : A_i \to A_{i+3}$$ is a quasi-isomorphism.

\end{enumerate}
Then the induced sequence on homology is exact.  Furthermore, the mapping cone of $f_1$ is quasi-isomorphic to $A_3$ via the map with components $H_1$ and $f_2$.
\end{lemma}

Let $Y_0$ be a closed, oriented 3-manifold, equipped with a framed knot $K_0$.  Applying the functor $\hmb$ to the associated 3-periodic sequence of elementary cobordisms $$\{W_n : Y_n \to Y_{n+1} \}_{n \in \mathbb{Z}/3\mathbb{Z}},$$ we obtain the surgery exact triangle:

\begin{theorem}
\label{thm:exact}
With coefficients in $\ftwo$, the sequence
$$\cdots \longrightarrow \hm(Y_{n-1}) \xrightarrow{\hm(W_{n-1})} \hm(Y_n) \xrightarrow{\hm(W_n)} \hm(Y_{n+1}) \longrightarrow \cdots$$
is exact.
\end{theorem}

\begin{proof}
We reorganize the proof in \cite{kmos} to fit it within our general framework of polytopes $\pij$ and identities $\aij$.  We use the notation $\{0,1,\infty,0^\prime\}$ for the lattice $\{1,2,3,4\}$ considered in \cite{kmos}.  The corresponding graph (in the sense of both $\Gamma$ and $G$) is the path of length three, yielding a pentagon of metrics $P_G$ whose sides correspond to $Y_1$, $Y_\infty$, $S_1 = S^1_\infty$, $S_\infty = S^\infty_{0^\prime}$, and $R_1 = S^1_{0^\prime}$ (where the left-hand notation is shorthand for the right-hand notation in the proof of Theorem \ref{thm:graph}).  The auxiliary hypersurface $R_1$ is diffeomorphic to $S^1 \times S^2$ and cuts out $V_1\cong \cptb - \text{int}(D^4)$ from $W$, leaving the cobordism $U_1$ with three boundary components.

\begin{figure}[htp]
\centering
\includegraphics[width=160mm]{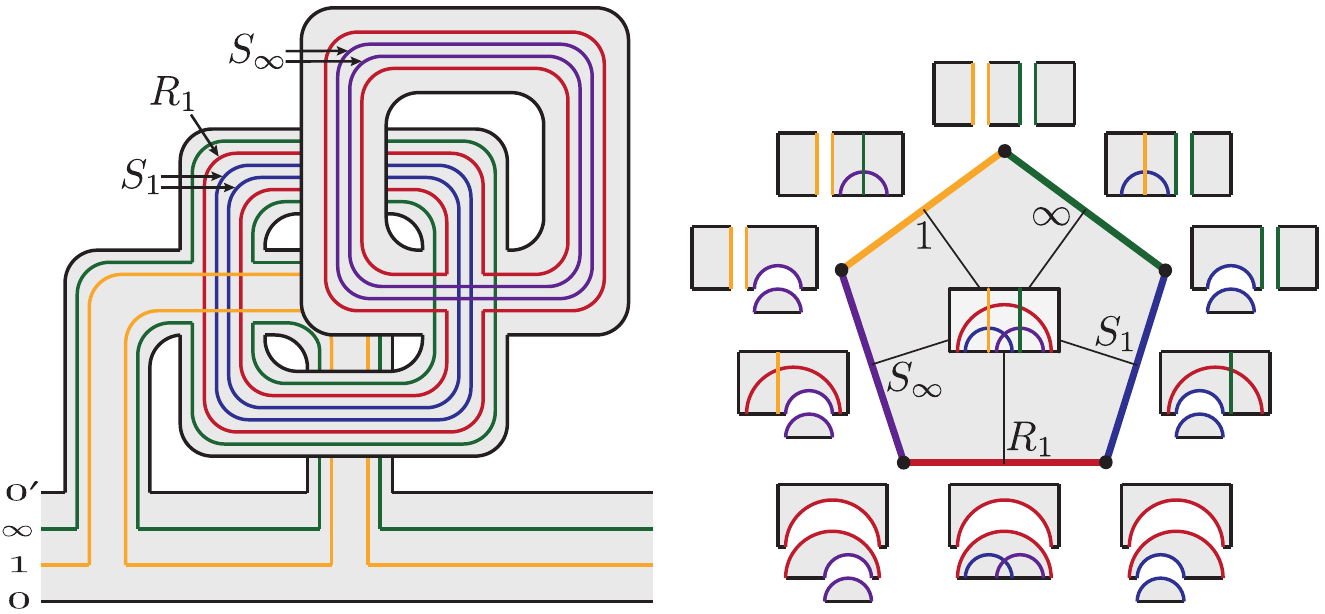}
\caption[The cobordism for $\{0,1,\infty,0^\prime\}$ and the pentagon of metrics.]{At left, the half-dimensional diagram of the cobordism $W$ for the lattice $\{0,1,\infty, 0^\prime\}$.  Note that $S_1$ is represented by two concentric curves, arising as the boundary of the tubular neighborhood of a circle representing the sphere $E_1$ (and similarly for $S_\infty$).  At right, the pentagon $K_4$ of metrics, analogous to the hexagon $P_3$ in Figure \ref{fig:hexagon}.}
\label{fig:intersections2}
\end{figure}

Keeping the 3-periodicity in mind, we prove exactness by applying Lemma \ref{lem:alg} with
\begin{align*}
A_{1+3j} &= \cm(Y_0) & \qquad & f_1 = \dd^0_1 & \qquad & H_1 = \dd^0_\infty & \qquad & \psi_1 = \dd^0_{0^\prime}\\
A_{2+3j} &= \cm(Y_1) & \qquad & f_2 = \dd^1_\infty & \qquad & H_2 = \dd^1_{0^\prime}\\
A_{3+3j} &= \cm(Y_\infty) & \qquad & f_3 = \dd^\infty_{0^\prime}
\end{align*}
where we have yet to define $\dd^0_{0^\prime}$.  The first condition of Lemma \ref{lem:alg} is then satisfied by Proposition \ref{prop:dd2} with $l = 1$.

Let $R$ denote the edge of the pentagon corresponding to $R_1$, considered as a one-parameter family of metrics on $V_1$ stretching from $S_1$ to $S_2$.  Viewing $V_1$ as a cobordism from the empty set to $R_1$, with the family of metrics $R$, we have components
\begin{align*}
n_o &\in C^o_\bullet(R_1) \qquad n_s  \in C^s_\bullet(R_1) \qquad \bar{n}_s \in C^s_\bullet(R_1) \qquad \bar{n}_u \in C^u_\bullet(R_1)
\end{align*}
In other words, these elements count isolated trajectories in moduli spaces of the form $M_z(V_1^*,\cmf)_R$ and $M_z^\text{red}(V_1^*,\cmf)_R$.  In fact, by Lemma 5.4 of \cite{kmos}, when the perturbation on $R_1$ is sufficiently small, there are no irreducible critical points and all components of the differential on $\cm(R_1)$ vanish, as do $n_o$ and $\bar{n}_s$.

We define the maps $D^*_*({}^0_{0^\prime})$ exactly as before.  We similarly define maps $\bar{D}^{ss}_s({}^0_{0^\prime})$ and $\bar{D}^{ss}_u({}^0_{0^\prime})$ which count isolated trajectories in $M^\text{red}_z(\amf,\cmf,U_1^*,\bmf)$:
\begin{align*}
\bar{D}^{ss}_s({}^0_{0^\prime}) : C^s_\bullet (R_1) \otimes C^s_\bullet(Y_0) \to C^s_\bullet(Y_{0^\prime}) \qquad \quad \ \ \bar{D}^{ss}_s({}^0_{0^\prime}) (e_\cmf \otimes e_\amf) = \sum_{\bmf \in \Csmf(Y_{0^\prime})} \sum_z \bar{m}_z(\amf, \cmf, U_1^*, \bmf) e_\bmf; \\
\bar{D}^{ss}_u({}^0_{0^\prime}) : C^s_\bullet (R_1) \otimes C^s_\bullet(Y_0) \to C^u_\bullet(Y_{0^\prime}) \qquad \quad \ \ \bar{D}^{ss}_u({}^0_{0^\prime}) (e_\cmf \otimes e_\amf) = \sum_{\bmf \in \Cumf(Y_{0^\prime})} \sum_z \bar{m}_z(\amf, \cmf, U_1^*, \bmf) e_\bmf.
\end{align*}
We combine these components to define the map $\dd^0_{0^\prime}: \cm(Y_0) \to \cm(Y_{0^\prime})$ by
\begin{align}
\nonumber
\dd^0_{0^\prime} & = \left[\begin{array}{rr}
D^o_o({}^0_{0^\prime}) & \sum_{0 \leq K \leq 0^\prime} D^u_o({}^K_{0^\prime}) \bar{D}^s_u({}^0_K) \\
D^o_s({}^0_{0^\prime}) & D^s_s({}^0_{0^\prime}) + \sum_{0 \leq K \leq 0^\prime} D^u_s({}^K_{0^\prime}) \bar{D}^s_u({}^0_K)
\end{array}\right] \\
\label{eqn:exactd2}
&+  \left[\begin{array}{rr}
0 & D^u_o ({}^{0^\prime}_{0^\prime})\bar{D}^{ss}_u({}^0_{0^\prime})(n_s \otimes \cdot) \\
0 & \bar{D}^{ss}_s({}^0_{0^\prime})(n_s \otimes \cdot) + D^u_s ({}^{0^\prime}_{0^\prime})\bar{D}^{ss}_u({}^0_{0^\prime})(n_s \otimes \cdot)
\end{array} \right],
\end{align}
\noindent which is written out in full in Appendix II.  The terms in \eqref{eqn:exactd2} break on a boundary-stable critical point in $\cm(R_1)$.  Of these, the term $\bar{D}^{ss}_s({}^0_{0^\prime})(n_s\otimes\cdot)$ is singly boundary-obstructed, while the other two are compositions of a non-boundary obstructed operator and a doubly boundary-obstructed operator (see Definition 24.4.4 in \cite{km}).  Finally, we introduce the chain map $\check{L}: \cm(Y_0) \to \cm(Y_{0^\prime})$ defined by
\begin{align}
\label{eqn:l}
\check{L} = \left[\begin{array}{rr}
L^o_o & L^u_o \bar{D}^s_u({}^0_0) + D^u_o({}^{0^\prime}_{0^\prime}) \bar{L}^s_u \\
L^o_s & \bar{L}^s_s +  L^u_s \bar{D}^s_u({}^0_0) + D^u_s({}^{0^\prime}_{0^\prime}) \bar{L}^s_u
\end{array}\right]
\end{align}
where $L^*_* = D^{u*}_*(\bar{n}_u \otimes \cdot)$ and  $\bar{L}^*_* = \bar{D}^{u*}_*(\bar{n}_u \otimes \cdot)$.  So the coefficient of $\bmf$ in $\check{L}(e_\amf)$ is a count of the zero-dimensional stratum of $M^+_z(\amf, \cmf, U_1^*, \bmf)$, over all $\cmf$ such that $e_\cmf$ is a summand of $\bar{n}_u$.

By Lemma \ref{lem:aij2} below, these maps are related by
\begin{align}
\label{eqn:aijlij1}
\dd{}^{0^\prime}_{0^\prime} \dd{}^0_{0^\prime} + \dd{}^0_{0^\prime} \dd{}^0_0 = \dd{}^1_{0^\prime} \dd{}^0_1 + \dd{}^\infty_{0^\prime} \dd{}^0_\infty + \check{L}.
\end{align}
Furthermore, by Proposition 5.6 of \cite{kmos}, $\check{L}$ is a quasi-isomorphism.  We conclude that $\dd{}^1_{0^\prime} \dd{}^0_1 + \dd{}^\infty_{0^\prime} \dd{}^0_\infty$ is a quasi-isomorphism as well.  This is precisely the second condition of Lemma \ref{lem:alg}, which then implies the theorem.
\end{proof}

\begin{remark}
In fact, the authors of \cite{kmos} show that the map induced by $\check{L}$ on $\hmb(Y_0)$ is given by multiplication by the power series $$\sum_{k \geq 0} U_\dagger^{k(k+1)/2}.$$  The proof is related to that of the blow-up formula, Theorem 39.3.1 of \cite{km}.
\end{remark}

Equation \eqref{eqn:aijlij1} is proved by counting ends.  The maps $A^*_*({}^0_{0^\prime})$ and $\bar{A}^*_*({}^0_{0^\prime})$ are defined using the vanishing elements $n_z(\amf, W^*,\bmf)_{P_{IJ}}$ and $\bar{n}_z(\amf, W^*,\bmf)_{P_{IJ}}$ exactly as before.  By analogy with the maps $D^{ss}_*$ above, we also define vanishing maps $\bar{A}^{ss}_s({}^0_{0^\prime})$ and $\bar{A}^{ss}_u({}^0_{0^\prime})$ which count boundary points of $M_z^{\text{red}+}(\amf,\cmf,U_1^*,\bmf)$.  Finally, we define $\check{A}^0_{0^\prime}: \cm(Y_0) \to \cm(Y_{0^\prime})$ by
\begin{align}
\nonumber
\check{A}^0_{0^\prime} & = \left[\begin{array}{rr}
A^o_o({}^0_{0^\prime}) & \sum_{0 \leq K \leq 0^\prime} \left(A^u_o({}^K_{0^\prime}) \bar{D}^s_u({}^0_K)   + D^u_o({}^K_{0^\prime}) \bar{A}^s_u({}^0_K) \right) \\
A^o_s({}^0_{0^\prime}) & A^s_s({}^0_{0^\prime}) + \sum_{0 \leq K \leq 0^\prime} \left(A^u_s({}^K_{0^\prime}) \bar{D}^s_u({}^0_K)   + D^u_s({}^K_{0^\prime}) \bar{A}^s_u({}^0_K) \right)
\end{array}\right] \\
\label{eqn:exacta2}
&+  \left[\begin{array}{rr}
0 & A^u_o ({}^{0^\prime}_{0^\prime})\bar{D}^{ss}_u({}^0_{0^\prime})(n_s \otimes \cdot) + D^u_o ({}^{0^\prime}_{0^\prime})\bar{A}^{ss}_u({}^0_{0^\prime})(n_s \otimes \cdot)\\
0 & \bar{A}^{ss}_s({}^0_{0^\prime})(n_s \otimes \cdot) + A^u_s ({}^{0^\prime}_{0^\prime})\bar{D}^{ss}_u({}^0_{0^\prime})(n_s \otimes \cdot) + D^u_s ({}^{0^\prime}_{0^\prime})\bar{A}^{ss}_u({}^0_{0^\prime})(n_s \otimes \cdot)
\end{array} \right],
\end{align}
which therefore vanishes as well.  The form of $\check{A}^0_{0^\prime}$ follows from the model case of Morse theory for manifolds with boundary, as described in Appendix I.  Note that all the terms in \eqref{eqn:exacta2} break on a boundary-stable critical point in $\cm(R_1)$.  The term $\bar{A}^{ss}_s({}^0_{0^\prime})(n_s \otimes \cdot)$ is singly boundary-obstructed, while the other four are compositions of a non-boundary-obstructed operator and a doubly-boundary-obstructed operator.  In Appendix II, we have written out these terms in expanded form in order to verify the following lemma.

\begin{lemma}
\label{lem:aij2}
The map $\check{A}^0_{0^\prime} +  \check{L}$ is equal to the component of $\dd^2$ from $\cm(Y_0)$ to $\cm(Y_{0^\prime})$:
$$\check{A}^0_{0^\prime} +  \check{L} = \sum_{0 \leq K \leq 0^\prime} \dd{}^K_{0^\prime} \dd{}^0_K.$$
\end{lemma}
\begin{proof}
As in the proof of Lemma \ref{lem:aij}, all terms on the right appear exactly once on the left, with the additional terms on the left being those which do not have a good break on any $Y_I$.  We divide these extra terms into those with
\renewcommand{\labelenumi}{(\roman{enumi})}
\begin{enumerate}
\item no break on $R_1$;
\item a boundary-stable break on $R_1$;
\item a boundary-unstable break on $R_1$. 
\end{enumerate}
Terms of type (i) can be enumerated just as in the proof of Lemma \ref{lem:aij}, so each occurs twice in $\check{A}^0_{0^\prime}$.  Dropping indices where it causes no ambiguity, the terms of type (ii) occur in six pairs:
$$\begin{array}{lclcr}
D^u_o \bar{D}^{su}_u(n_s \otimes \bar{D}^s_u(\cdot)) & \text{in} & A^u_o\bar{D}^s_u & \text{and} & D^u_o \bar{A}^{ss}_u(n_s \otimes \cdot); \\
D^u_o \bar{D}^u_u \bar{D}^{ss}_u(n_s \otimes \cdot) &\text{in} & A^u_o \bar{D}^{ss}_u(n_s \otimes \cdot) & \text{and} & D^u_o\bar{A}^{ss}_u(n_s \otimes \cdot); \\
D^u_s \bar{D}^{su}_u(n_s \otimes \bar{D}^s_u(\cdot)) &\text{in} & A^u_s\bar{D}^s_u & \text{and} & D^u_s \bar{A}^{ss}_u(n_s \otimes \cdot); \\
D^u_s \bar{D}^u_u \bar{D}^{ss}_u(n_s \otimes \cdot) &\text{in} & A^u_s \bar{D}^{ss}_u(n_s \otimes \cdot) & \text{and} & D^u_s\bar{A}^{ss}_u(n_s \otimes \cdot); \\
\bar{D}^{su}_u(n_s \otimes \bar{D}^s_u(\cdot)) &\text{in} & A^u_s\bar{D}^s_u & \text{and} & \bar{A}^{ss}_s(n_s \otimes \cdot); \\
\bar{D}^u_s \bar{D}^{ss}_u(n_s \otimes \cdot) &\text{in} & A^u_s \bar{D}^{ss}_u(n_s \otimes \cdot) & \text{and} & \bar{A}^{ss}_s(n_s \otimes \cdot).
\end{array}$$
Finally, the terms of type (iii) occur in five pairs:
$$\begin{array}{lcrcr}
D^u_o D^{us}_u(\bar{n}_u \otimes \cdot) & \text{in} & D^u_o\bar{A}^s_u & \text{and} & D^u_o \bar{L}^s_u; \\
D^{uu}_o(\bar{n}_u \otimes \bar{D}^s_u(\cdot)) &\text{in} & A^u_o \bar{D}^s_u & \text{and} & L^u_o \bar{D}^s_u; \\
D^u_s D^{us}_u(\bar{n}_u \otimes \cdot) & \text{in} & D^u_s\bar{A}^s_u & \text{and} & D^u_s \bar{L}^s_u; \\
D^{uu}_s(\bar{n}_u \otimes \bar{D}^s_u(\cdot)) &\text{in} & A^u_s \bar{D}^s_u & \text{and} & L^u_s \bar{D}^s_u; \\
\bar{D}^{us}_s(\bar{n}_u \otimes \cdot) &\text{in} & \bar{A}^s_s & \text{and} & \bar{L}^s_s.
\end{array}$$
We conclude that terms of types (i) and (ii) are double counted by $\check{A}^0_{0^\prime}$ while those of type (iii) are counted once each by $\check{A}^0_{0^\prime}$ and $\check{L}$.  We therefore have equality over $\ftwo$.
\end{proof}

\begin{remark}
If we consider a boundary-unstable break on $R_1$ to be a good break as well, then Remark \ref{rem:badbreak} goes through exactly as before.  Furthermore, $\check{L}$ counts those trajectories which break well on $R_1$ (see also the discussion following Proposition 5.5 in \cite{kmos}).
\end{remark}

\begin{remark}
For the lattice $\{0,1,\infty, 0^\prime\}$, we introduced the auxiliary hypersurfaces $S_1$, $S_2$, and $R_1$ in order to build the pentagon of metrics.  The $S_i$ edges contribute vanishing terms to $\qqij$ by  Lemma \ref{lem:cp2}, whereas the $R_1$ edge contributes the term $\check{L}$.  Thus,  $$\qqij = \dd{}^1_{0^\prime} \dd{}^0_1 + \dd{}^\infty_{0^\prime} \dd{}^0_\infty + \check{L}$$ and once more we can view \eqref{eqn:aijlij1} as a ``generalization'' of \eqref{eqn:qnull}.
\end{remark}

\section{The link surgery spectral sequence: convergence}
\label{sec:linksurg2}

We are now positioned to identify the limit of the link surgery spectral sequence.

\begin{proof}[Proof of Theorem \ref{thm:linksurg}]
For $1 \leq k \leq l$, define the map $$F_k : \bigoplus_{I \in \{\infty\}^{l-k} \times \{0,1\} \times \{0,1\}^{k-1}} \cm(\yi) \longrightarrow \bigoplus_{I \in \{\infty\}^{l-k} \times \{\infty\} \times \{0,1\}^{k-1}} \cm(\yi)$$ as the sum of all relevant components of the differential $\dd$ on the subcomplex
$$\bigoplus_{I \in \{\infty\}^{l-k} \times \{0,1,\infty\} \times \{0,1\}^{k-1}} \cm(\yi)$$ of $\tilde X$.  Then $\dd^2 = 0$ implies that $F_k$ is a chain map.   Consider the filtration given by the weight of the last $k-1$ digits of $I$.  By applying the final assertion of Lemma \ref{lem:alg} to the surgery exact triangles arising from the component $K_{l-k+1}$, we conclude that $F_k$ induces an isomorphism between the $E^1$ pages of the associated spectral sequences.
Therefore, $F_k$ is a quasi-isomorphism, as is the composition
\begin{align}
\label{eqn:quasi}
F = F_1 \circ F_2 \circ \cdots \circ F_l : X \to \cm(Y_\infty).
\end{align}
\end{proof}

\begin{remark}
The proof of Theorem \ref{thm:linksurg} hinges on two facts:
\renewcommand{\labelenumi}{(\roman{enumi})}
\begin{enumerate}
\item lattices of the form $\{0,1\}^k$ and $\{0,1,\infty\} \times \{0,1\}^k$ give rise to filtered complexes;
\item the lattice $\{0,1,\infty, 0^\prime \}$ gives rise to an exact sequence.
\end{enumerate}
We considered more general lattices in Theorem \ref{thm:graph} and Proposition \ref{prop:dd2} in part to make clear how both these facts arise as special cases of the same polytope constructions.  The lattice $\{0,1,\infty,0'\} \times \{0,1\}$ will arise naturally in Section \ref{sec:umap}.
\end{remark}

\subsection{Grading}
\label{sec:grading}

The group $\hmb(Y)$ is endowed with an absolute mod 2 grading $\grt$, as explained in Sections 22.4 and 25.4 of \cite{km}.  This gradings is uniquely characterized by two properties.  First, the group $\hmb(S^3)$ is supported in even grading.  Second, if $W$ is a cobordism from $Y_-$ to $Y_+$, then the map $\check m(W): \cm(Y_-) \to \cm(Y_+)$ shifts $\grt$ according to the parity of the integer
\begin{align}
\label{eqn:iota}
\iota(W) = \frac{1}{2}\left(\chi(W) + \sigma(W) +  b_1(Y_+) - b_1(Y_-)\right),
\end{align}
where $\chi$ is the Euler number and $\sigma$ is the signature of the intersection form on $I^2(W) = \text{Im}\left(H^2(W, \partial W) \to H^2(W)\right)$.
Note that $\iota$ is additive under composition, since both the signature and Euler characteristic are additive in this context.  Furthermore, if $P$ parameterizes an $n$-dimensional family of metrics on $W$, then the map $\check m(W)_P$ shifts $\grt$ by $\iota(W) + n$.

We now introduce an absolute mod 2 grading $\check{\delta}$ on the hypercube complex $(X,\dd)$ which reduces to $\grt$ in the case $l = 0$.  In fact, it will be useful to define $\hat{\delta}$ on the larger complex $(\widetilde{X}, \dd)$ associated to the lattice $\{0,1,\infty\}^l$.  Let $x \in \cm(\yi)$ be homogeneous with respect to the $\grt$ grading.  Then for $l > 0$, we define
\begin{align}
\check{\delta}(x)
\nonumber
&= \grt(x) + \left(\iota(W_{\Ii I}) - w(I)\right) - \left(\iota(W_{\Ii \infty}) - 2l\right) - l \mod 2 \\
\label{eqn:form1} &=\grt(x) - \left(\iota(W_{I \infty}) + w(I)\right) + l \mod 2.
\end{align}
Here the subscripts $\Ii$ and $\infty$ are shorthand for the initial and final vertices of $\{0,1,\infty\}^l$.

\begin{lemma}
\label{lem:byone}
The differential $\dd$ on $\widetilde{X}$ and $X$ lowers $\hat{\delta}$ by 1.
\end{lemma}
\begin{proof}
Since $\dij$ is defined using a family of metrics of dimension $w(J) - w(I) - 1$ on $\wij$, it shifts $\grt$ by $$\iota(\wij) + \left(w(J) - w(I) - 1\right) = \left(\iota(W_{I\infty}) + w(I)\right) - \left(\iota(W_{J\infty}) + w(J)\right)  - 1.$$  The claim now follows from \eqref{eqn:form1}.
\end{proof}

We now complete the proof of Theorem \ref{thm:linksurg1}.

\begin{proposition}
\label{prop:quasi}
The gradings $\check \delta$ and $\grt$ coincide under the quasi-isomorphism $$F : X \to \cm(Y_\infty)$$ defined in $\eqref{eqn:quasi}$.
\end{proposition}
\begin{proof}
The weight of the vertex $\{\infty\}^l$ is $2l$.  Therefore, given $x \in \cm(Y_\infty)$, by \eqref{eqn:form1} we have
$$\check{\delta}(x) = \grt(x) - l \mod 2.$$
So it suffices to show that the quasi-isomorphism $F: X \to \cm(Y_\infty)$ lowers $\check \delta$ by $l$.  But $F$ is a composition of $l$ maps $F_k$, each of which is a sum of components of $\dd$.  So we are done by the Lemma \ref{lem:byone}.
\end{proof}

\subsection{Invariance}

The construction of the hypercube complex
$$X(g,q) = \bigoplus_{I \in \{0,1\}^l} \cm(\yi(g|_I, q|_I))$$
depends immaterially on numbering the components of $L$, and materially on a choice of regular metric $g$ and perturbation $q$ on the full cobordism $W$, where the metric is cylindrical near each of the hypersurfaces $\yi$.  Let $(g_0,q_0)$ and $(g_1, q_1)$ be two such choices.

\begin{theorem}
\label{thm:anal}
There exists a $\check t$-filtered, $\check \delta$-graded chain homotopy equivalence $$\phi: X(g_0,q_0) \to X(g_1,q_1),$$ which induces a $(\check t, \check \delta)$-graded isomorphism between the associated $E^i$ pages for all $i \geq 1$.
\end{theorem}

\begin{proof}
We start by embedding a second copy of each $Y_I$ in $W$ as follows (see Figure \ref{fig:metric} for the case $l = 2$).  First, relabel the incoming end $Y_\Ii$ as $Y_{\Ii \times \{0\}}$ and every other $Y_I$ as $Y_{I \times \{1\}}$.  Then embed a second copy of $Y_{\Ii \times \{0\}}$, labeled $Y_{\Ii \times \{1\}}$, just above the original.  Finally, embed a second copy of each $Y_{I \times \{1\}}$, labeled $Y_{I \times \{0\}}$, just below the original.  We now have an embedded hypersurface $Y_{I \times \{i\}}$ for each $I \times \{i\}$ in the hupercube  $\{0,1\}^l \times \{0,1\}$, with diffeomorphisms
\begin{align}
\label{eqn:wiidiff}
W_{I \times \{0\}, I \times \{1\}} & \cong \yi \times [0,1] \\
\label{eqn:wijdiff}
W_{I \times \{i\}, J \times \{j\}} & \cong \wij
\end{align}
where in \eqref{eqn:wijdiff} we assume $I<J$.  Furthermore, $Y_{I \times \{i\}}$ and $Y_{J \times \{j\}}$ are disjoint if $I \times \{i\}$ and $J \times \{j\}$ are ordered.

\begin{figure}[htp]
\centering
\includegraphics[width=160mm]{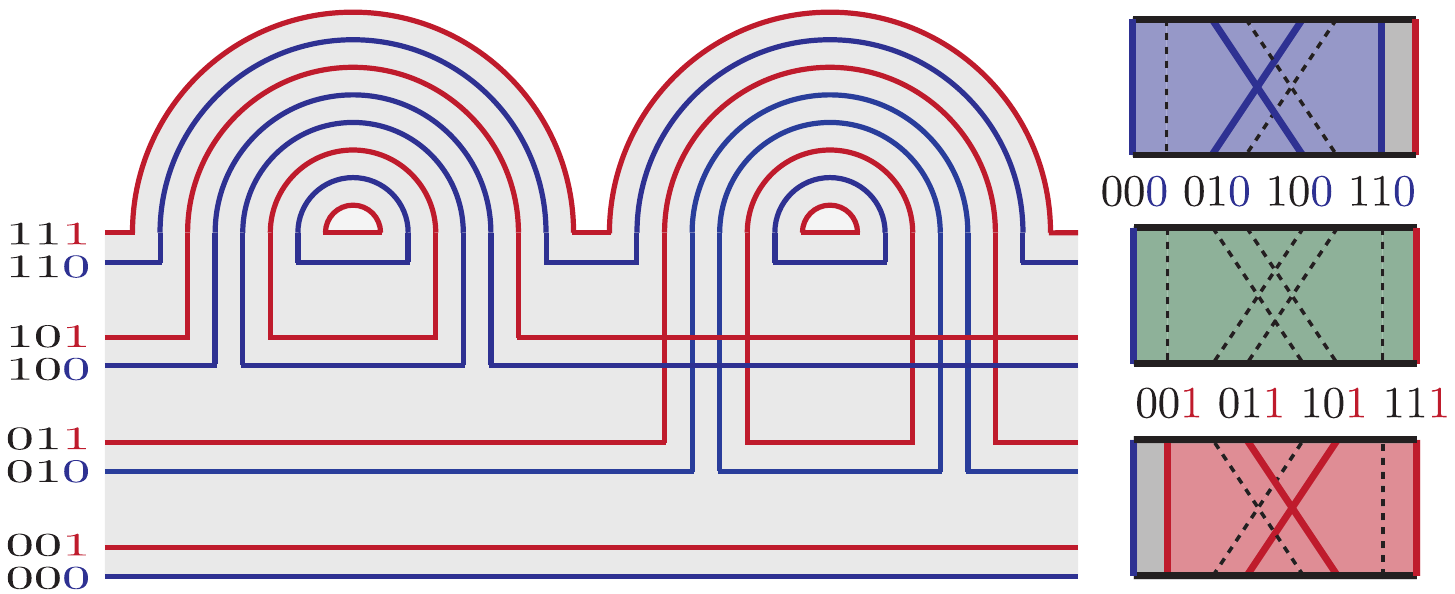}
\caption[The cobordism used to the construct the homotopy equivalence for $\{0,1\}^2$.]{At left, we have the half-dimensional diagram of the cobordism $W$ used to prove analytic invariance in the case $l = 2$.  For each $I \in \{0,1\}^2$, the hypersurfaces $Y_{I\times\{0\}}$ (in blue) and $Y_{I\times\{1\}}$ (in red) bound a cylindrical cobordism.  At right, we can fix the blue metric $g_0$ on $W_{000,110}$ (top), or the red metric $g_1$ on $W_{001,111}$ (bottom).  The green metric on the middle rectangle represents an intermediate state.  To construct the homotopy, we slide the metric from that on the top rectangle to that on the bottom rectangle in a controlled manner, as explained in Figure \ref{fig:invariance}.}
\label{fig:metric}
\end{figure}

Our strategy is as follows.  We define a complex
\begin{align*}
\underline {\check X} = \displaystyle \bigoplus_{I \in \{0,1\}^l, i \in \{0,1\}} \cm(Y_{I \times \{i\}}),
\end{align*}
where the differential $\underline\dd$ is defined as a sum of components $$\underline{\dd}^{I \times \{i\}}_{J \times \{j\}} : \cm(Y_{I \times \{i\}}) \to \cm(Y_{J \times \{j\}}).$$  Those components of the form $\underline{\dd}^{I \times \{i\}}_{J \times \{i\}}$ are inherited from $X(g_i,q_i)$.  So we may view $X(g_0, q_0)$ as the complex over $\{0,1\}^l \times \{0\}$ obtained from quotienting $\underline {\check X}$ by the subcomplex $X(g_1, q_1)$ over $\{0,1\}^l \times \{1\}$.   The component $\underline{\dd}^{I \times \{0\}}_{J \times \{1\}}$ is induced by the cobordism $W_{I \times \{i\}, J \times \{j\}}$ over a family of metrics and perturbations parameterized by a permutohedron $\underline{ \check P}_{I \times \{0\}, J \times \{1\}}$, to be defined momentarily.  Then $\underline \dd^2 = 0$ implies that $$\phi = \sum_{I \leq J} \underline{\dd}^{I \times \{0\}}_{J \times \{1\}} : X(g_0, q_0) \to X(g_1, q_1)$$ is a chain map.  If we extend the $\check \delta$ grading verbatim to $\underline {\check X}$, then $\phi$ is odd as a map on  $\underline {\check X}$ by Proposition \ref{prop:quasi}, and thus even as a map from $X(g_0,q_0)$ and $X(g_1,q_1)$.  Thus, $\phi$ is $\check \delta$-graded, and it is clearly $\check t$-filtered.  By \eqref{eqn:wiidiff}, the map
$$\underline{ \dd}^{I \times \{0\}}_{I \times \{1\}} : \cm(Y_{I \times \{0\}}) \to \cm(Y_{I \times \{1\}})$$
induces an isomorphism on homology.  Thus, filtering by the horizontal weight $\underline w$ defined by $\underline w(I \times \{i\}) = w(I)$, $\phi$ induces a $(\check t, \check \delta)$-graded isomorphism between the $E^1$ pages of the corresponding spectral sequences.  By Theorem 3.5 of  \cite{mcc}, we conclude that $\phi$ induces a $(\check t, \check \delta)$-graded isomorphism between the $E^i$ pages for each $i \geq 1$.  Thus, $\phi$ is a quasi-isomorphism, and therefore (since we are working over a field) a homotopy equivalence.

It remains to construct the family parameterized by each $\underline{\check P}_{I \times \{0\}, J \times \{1\}}$ and to prove that $\underline\dd^2 = 0$.  We start by fixing a metric $g^I_I$ on each cylindrical cobordism $W_{I \times \{0\}, I \times \{1\}}$ for which $g^I_I(Y_{I \times \{0\}}) = g_0(\yi)$ and $g^I_I(Y_{I \times \{1\}}) = g_1(\yi)$ (we proceed similarly with regard to the perturbations, though we will suppress this).  Here the notation $g(Y)$ denotes the restriction of $g$ to $Y$.   The point $\underline{\check P}_{I \times \{0\}, I \times \{1\}}$ is defined to correspond to the metric $g^I_I$.  Now for each $I \in \{0,1\}^l$, we specify a metric $g_I$ on $W$ by its restriction to each of three pieces:
\begin{align*}
g_I(W_{\Ii \times \{0\}, I \times \{0\}}) &= g_0(W_{\Ii I}) \\
g_I(W_{I \times \{0\}, I \times \{1\}}) &= g^I_I \\
g_I(W_{I \times \{1\}, \If \times \{1\}}) &= g_1(W_{I\If}).
\end{align*}
We will use these metrics to construct the family parameterized by $\underline{\check P}_{\Ii \times \{0\}, \If \times \{1\}}$ in several stages.  The case $l=2$ is illustrated in Figure \ref{fig:invariance}.

\begin{figure}[htp]
\centering
\includegraphics[width=160mm]{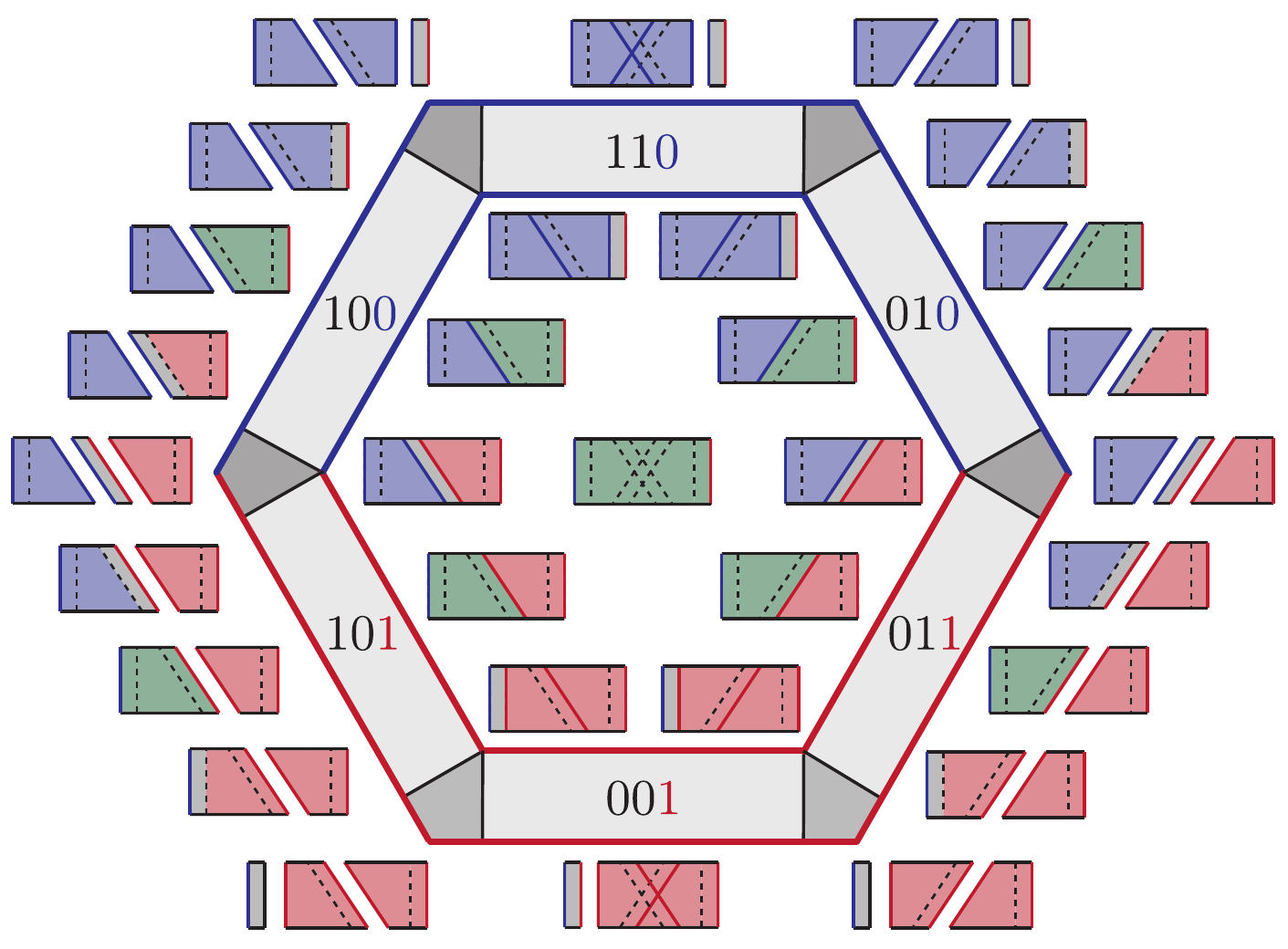}
\caption[The hexagon of metrics used to the construct the homotopy equivalence for $\{0,1\}^2$.]{The hexagon $\underline{\check P}^{000}_{111}$ is drawn so that increasing the vertical coordinate is suggestive of moving from the red metrics to the blue metrics.  Gray represents the cylindrical metrics $g^I_I$, while green represents an intermediate mixture of red, blue, and gray.}
\label{fig:invariance}
\end{figure}

\begin{figure}[htp]
\centering
\includegraphics[width=145mm]{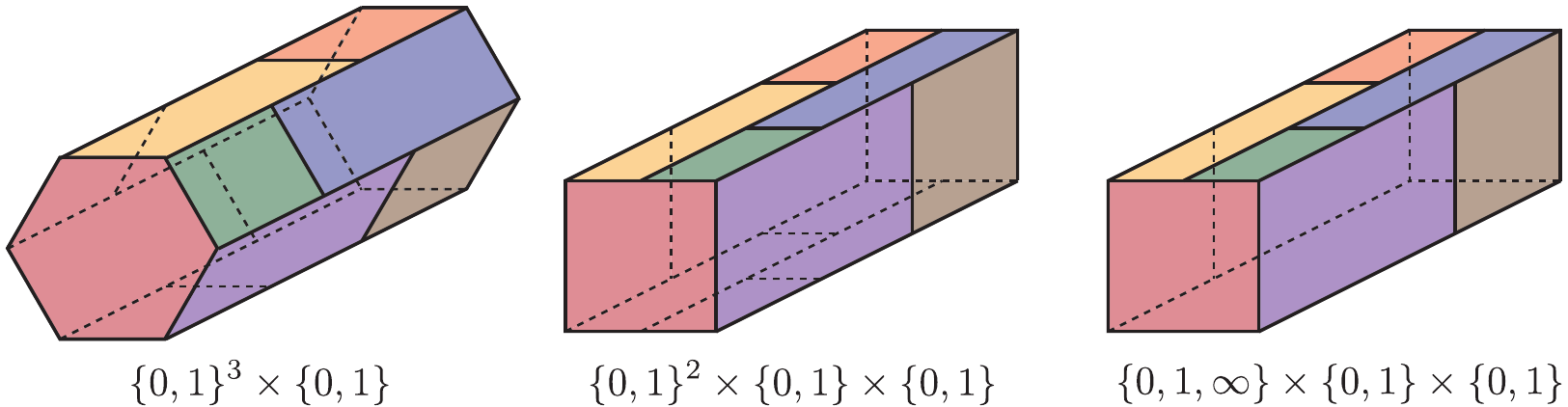}
\caption[Realizations of the permutohedron and a graph associahedron.]{Two alternative realizations of the permutohedron from Figure \ref{fig:permutohedron}, and one of the graph associahedron from Figures \ref{fig:dual} and \ref{fig:truncation}.}
\label{fig:induct}
\end{figure}

We first describe a family $\mathcal{F}$ of non-degerate metrics on $W$, parameterized by the permutohedron $P_{l+1}$.  Let $Q_I$ denote the facet of $P_l$ corresponding to the internal vertex $I$.  $P_{l+1}$ may be obtained from $P_l \times [0,l]$ by subdividing each facet $Q_I \times [0,1]$ by the ridge $Q_I \times \{w(I)\}$.  In the $l=2$ case, this amounts to adding a vertex at the midpoint of each vertical edge in a square.  In the $l = 3$ case,  shown at right in Figure \ref{fig:induct}, we have cross the hexagon with an interval and add an edge to each lateral face.  We next label the facets $Q_I \times [0,w(I)]$ and $Q_I \times [w(I),l]$ by  $I \times \{1\}$ and $I \times \{0\}$, respectively.  Furthermore, we label $P_l \times \{0\}$ and $P_l \times \{l\}$ by $\Ii \times \{1\}$ and $\If \times \{0\}$, respectively.

We then associate the metric $g_I$ to each vertex of $P_{l+1}$ lying on $Q_I \times \{w(I)\}$.  The remaining vertices of $P_{l+1}$ lie on $P_l \times \{0\}$ or $P_l \times \{l\}$.  We associate to these vertices the metrics $g_\Ii$ and $g_\If$, respectively (note that $w(\Ii) = 0$ and $w(\If) = l$).  At this stage, we have defined $\mathcal{F}$ on the 0-skeleton of $P_{l+1}$.  We proceed inductively: having extended $\mathcal{F}$ to the boundary $\partial F$ of a k-dimensional face $F$ of $P_{l+1}$, we extend $\mathcal{F}$ to the interior of $F$, subject to the following constraint:
\begin{align}
\label{eqn:constraint}
\text{If $\mathcal{F}|_{\partial F}$ is constant over some hypersurface or component of $W$, then so is $\mathcal{F}|_F$.}
\end{align}
In particular, the family $\mathcal{F}$ is constant when restricted to each of the facets $P_l \times \{0\}$ and $P_l \times \{l\}$ and each of the ridges $Q_I \times \{w(I)\}$.

The family $\mathcal{F}$ over $P_{l+1}$ slides the metric (and perturbation) on $W$ in stages (in Figure  \ref{fig:invariance}, $P_{2+1}$ is the inner hexagon).  We now extend $\mathcal{F}$ to a family $\mathcal{G}$ which incorporate stretching.  To each facet $Q_{I \times \{i\}}$ of $P_{l+1}$, we glue the polytope $Q_{I \times \{i\}} \times [0,\infty]$ along the facet $Q_{I \times \{i\}} \times \{0\}$ (in Figure  \ref{fig:invariance}, these are the six lightly shaded rectangles).  We extend $\mathcal{G}$ over $Q_{I \times \{i\}} \times [0,\infty]$ by stretching on $Y_{I \times \{i\}}$ in accordance with the latter coordinate (recall that the metric on $Y_{I \times \{i\}}$ is constant over $Q_{I \times \{i\}}$).  Next, along each ridge $Q_{I \times \{i\} < J \times \{j\}}$ in $P_{l+ 1}$, we glue on the polytope $Q_{I \times \{i\} < J \times \{j\}} \times [0, \infty] \times [0,\infty]$ in the obvious manner (in Figure \ref{fig:invariance}, these are the six heavily shaded squares).  The first interval parameterizes stretching on $Y_{I \times \{i\}}$ while the second interval parameterizes stretching on $Y_{J \times \{j\}}$.  We continue this process until the last stage, when we glue one cube $[0,\infty]^l$ at each vertex of $P_{l+1}$, over which $\mathcal{G}$ stretches on the corresponding maximal chain of internal hypersurfaces.

In the end, we have simply thickened the boundary of $P_{l+1}$ to describe a family $\mathcal{G}$ of metrics on $W$ parameterized by the permutohedron $\underline{\check P}_{\Ii \times \{0\}, \If \times \{1\}}$ (the full hexagon in Figure \ref{fig:invariance}).  This family is degenerate over the boundary of $\underline{\check P}_{\Ii \times \{0\}, \If \times \{1\}}$ precisely as described by Proposition \ref{prop:perm}.  Now, for each $I \leq J$, we construct a family of metrics $\mathcal{G}_{IJ}$ over $\underline{\check P}_{I \times \{0\}, J \times \{1\}}$ by restricting the family $\mathcal{G}$ to $W_{I \times \{0\}, J \times \{1\}}$ over an appropriate face of $\underline{\check P}_{\Ii \times \{0\}, \If \times \{1\}}$ (here the constraint \eqref{eqn:constraint} is essential).

The proof that $\underline\dd^2 = 0$ now lifts directly from the original proof that $\dd^2 = 0$, with one new point that we now explain.  The component of $\underline\dd^2$ from $\cm(Y_{\Ii \times \{0\}})$ to $\cm(Y_{\If \times \{1\}})$ vanishes if and only if
\begin{align}
\nonumber
\underline{\dd}^{\Ii \times \{0\}}_{\If \times \{1\}} \dd^{\Ii \times \{0\}}_{\Ii \times \{0\}} + \dd^{\If \times \{1\}}_{\If \times \{1\}}\underline{\dd}^{\Ii \times \{0\}}_{\If \times \{1\}} &= \underline{\dd}^{\If \times \{0\}}_{\If \times \{1\}} \dd^{\Ii \times \{0\}}_{\If \times \{0\}} + \dd^{\Ii \times \{1\}}_{\If \times \{1\}} \underline{\dd}^{\Ii \times \{0\}}_{\Ii \times \{1\}} \\
\label{eqn:lateral}
&+ \sum_{\Ii < I < \If} \underline{\dd}^{I \times \{0\}}_{\If \times \{1\}} \dd^{\Ii \times \{0\}}_{I \times \{0\}} + \dd^{I \times \{1\}}_{\If \times \{1\}} \underline{\dd}^{\Ii \times \{0\}}_{I \times \{1\}}.
\end{align}
Consider the composite map corresponding to the family $\mathcal{G}$ over the facet $\If \times \{0\}$ of $\underline{\check{P}}_{\Ii \times \{0\}, \If \times \{1\}}$.  Since the family $\mathcal{F}$ over the corresponding facet of $P_{l+1}$ is constant, the only sections of the facet $\If \times \{0\}$ which contributes non-trivially to this map are those of the form $\{\infty\} \times [0,\infty]^{l-l}$ in the boundary of the cubes $[0,\infty]^{l}$ (in Figure \ref{fig:invariance}, these are the two segments of the top edge of the hexagon which lie in the boundary of the heavily shaded squares).  The other sections cannot give rise to 0-dimensional moduli spaces, since they involve at least one parameter which does not change the metric.  We can therefore identify the map associated to the facet $\If \times \{0\}$ with $\underline{\dd}^{\If \times \{0\}}_{\If \times \{1\}} \dd^{\Ii \times \{0\}}_{\If \times \{0\}}$ (in Figure \ref{fig:invariance}, we are contracting out the middle segment of the top edge).  Similarly, the map associated to the facet $\Ii \times \{1\}$ coincides with $\underline{\dd}^{\If \times \{0\}}_{\If \times \{1\}} \dd^{\Ii \times \{0\}}_{\If \times \{0\}}$, and the sum on line \eqref{eqn:lateral} coincides with the map associated to the remaining lateral facets of $\underline{\check{P}}_{\Ii \times \{0\}, \If \times \{1\}}$.  Thus, the full equation expresses the fact that the map $\underline{\dd}^{\Ii \times \{0\}}_{\If \times \{1\}}$ associated to the full permutohedron is a null-homotopy for the map associated to its boundary.  The other components of $\dd^2$ vanish by a completely analogous argument.
\end{proof}

\begin{remark}
Recall the top and bottom rectangles at right in Figure \ref{fig:metric}.  Suppose that the red and blue metrics agree where they overlap, so that the family $\mathcal{F}$ on $P_{l+1}$ can be made completely constant.  Then only the cubes $[0,\infty]^l$ contribute non-trivially to the map $\underline{\dd}^{\Ii \times \{0\}}_{\If \times \{1\}}$.  Discarding the rest of $\underline{\check P}_{\Ii \times \{0\}, \If \times \{1\}}$ and gluing these cubes together, we build a permutohedron giving rise to the same map.  This viewpoint highlights the connection between the permutohedra $\underline{\check P}_{I \times \{0\}, J \times \{1\}}$ and the permutohedra $\pij$ that we first constructed in Section \ref{sec:surface}, using only cubes which stretch the metric along maximal chains of internal hypersurfaces.
\end{remark}


\section{The $U_\dagger$ map and $\hmhat(Y)$}
\label{sec:umap}

Given a cobordism $W : Y_0 \to Y_1$, Kronheimer and Mrowka construct a map $\hmb(U_\dagger | W) : \hmb(Y_0) \to \hmb(Y_1)$.   In \cite{km}, this map is defined by pairing each moduli space $M_z(\mathfrak{a}, W^*, \mathfrak{b})$ with the first Chern class of the natural complex line bundle on $\mathcal{B}^\sigma(W)$.  A dual description of the map is given in \cite{kmos}.  We will use the notation $\check{m}(U_\dagger|W)$ for this map on the chain level.

We introduce a third description which fits in neatly with our previous constructions.  We first recall some facts about the monopole Floer homology of the 3-sphere (see Sections 22.7 and 25.6 of \cite{km}).  With round metric and small perturbation, the monopoles on the 3-sphere consist of a single bi-infinite tower $\{\mathfrak{c}_i\}_{i \in \zz}$ of reducibles, with $\mathfrak{c}_i$ boundary-stable if and only if $i \geq 0$, and $\grq(\mathfrak{c}_i) = 2i$.  Furthermore, $U_\dagger$ sends $\mathfrak{c}_i$ to $\mathfrak{c}_{i-1}$, and in particular, $\langle U_\dagger, M(D^{4*}, \mathfrak{c}_i)\rangle$ is non-zero if and only if $i = -1$.  It is this last property which motivates the following reformulation.

Given a cobordism $W : Y_0 \to Y_1$, let $W^{**}$ denote the manifold obtained by removing a ball from the interior of $W$ and attaching cylindrical ends to all three boundary components, with the new $S^3 \times [0,\infty)$ end regarded as incoming.  Choose the metric and perturbation on $W$ so that we return to the situation described in the last paragraph over $S^3$.  We define the map $\check{m}(U | W): \cm(Y_0) \to \cm(Y_1)$ by replacing each moduli space $M_z(\mathfrak{a}, W^*, \mathfrak{b})$ in the definition of $\check{m}(W)$ with the moduli space $M_z(\mathfrak{a}, \mathfrak{c}_{-1}, W^{**}, \mathfrak{b})$.  In other words,
\begin{align*}
\check{m}(U | W) = \left[\begin{array}{rr}
m^{uo}_o(\mathfrak{c}_{-1} \otimes \cdot) & m^{uu}_o(\mathfrak{c}_{-1} \otimes \esu(\cdot)) + \duo\bar{m}^{us}_{u}(\mathfrak{c}_{-1} \otimes \cdot)\\
m^{uo}_s(\mathfrak{c}_{-1} \otimes \cdot) & \bar{m}^{us}_s(\mathfrak{c}_{-1} \otimes \cdot) + m^{uu}_s(\mathfrak{c}_{-1} \otimes \esu(\cdot)) + \dus\bar{m}^{us}_{u}(\mathfrak{c}_{-1} \otimes \cdot)
\end{array}\right].
\end{align*}
One sees that $\check{m}(U | W)$ is a chain map and well-defined up to homotopy equivalence by the same argument used for $\check{m}(W)$, together with the fact that there are no isolated trajectories from $\mathfrak{c}_{-1}$ to any other $\mathfrak{c}_i \in \mathfrak{C}(S^3)$.  Note that the choice of ball in $W$ may be interpreted as a metrical choice, since a diffeomorphism $\phi$ of $W$ sending one ball to another is an isometry from $(W,g)$ to $(\phi(W),(\phi^{-1})^*g)$.

\begin{proposition}
The map $\check{m}(U | W)$ is homotopy equivalent to the map $\check{m}(U_\dagger | W)$.
\end{proposition}
\begin{proof}
Given a cobordism $W$, we may assume the cochain $u$ representing $U_\dagger$ is supported over the configuration space of a small ball.  The homotopy relating the two maps is now given by stretching the metric normal to the 3-sphere bounding this ball.
\end{proof}
This justifies a return to the notation $\check{m}(U_\dagger | W)$.

We now define a fourth version of monopole Floer homology, denoted $\hmhat(Y)$ and motivated by the properties of $\hf(Y)$ in Heegaard Floer homology.  We use the shorthand $U_\dagger$ for the map $\check{m}(U_\dagger | Y \times [0,1]) : \cm(Y) \to \cm(Y)$ induced by the cylindrical cobordism, with some regular choice of metric and perturbation (which need not have the same restrictions on both ends).
The complex $\cmhat(Y)$ is defined to be the mapping cone of $U_\dagger$:
$$\cmhat(Y) = \cm(Y) \textstyle\bigoplus \cm(Y)[1] \hspace{1in} \tilde\partial = \left[\begin{array}{cc}\check\partial & 0 \\ U_\dagger & \check\partial\end{array}\right]$$
Since $U_\dagger$ is an even map, the differential $\tilde\partial$ is odd, and therefore $\grt$ naturally extends to $\cmhat(Y)$ (as does $\grq$ for torsion spin$^c$ structures).  We then define $\hmhat(Y)$ as the homology $H_*(\cmhat(Y), \tilde\partial)$.  While the completion involved in the definition of the other versions has no effect here, we keep the bullet for notational consistency.  We will see that $\hmhat$ describes a functor in the same spirit as $\hmb$.  By construction, there is an exact sequence
\begin{align}
\label{eqn:hatlong}
\cdots  \xrightarrow{j} \hmhat(Y) \xrightarrow{i} \hmb(Y) \xrightarrow{U_\dagger} \hmb(Y) \xrightarrow{j} \cdots
\end{align}
of $\ftwo[[U_\dagger]]$ modules where $U_\dagger$ acts by zero on $\hmhat(Y)$.  Here the maps $i$ and $U_\dagger$ are even, while $j$ is odd.

\begin{figure}[htp]
\centering
\includegraphics[width=155mm]{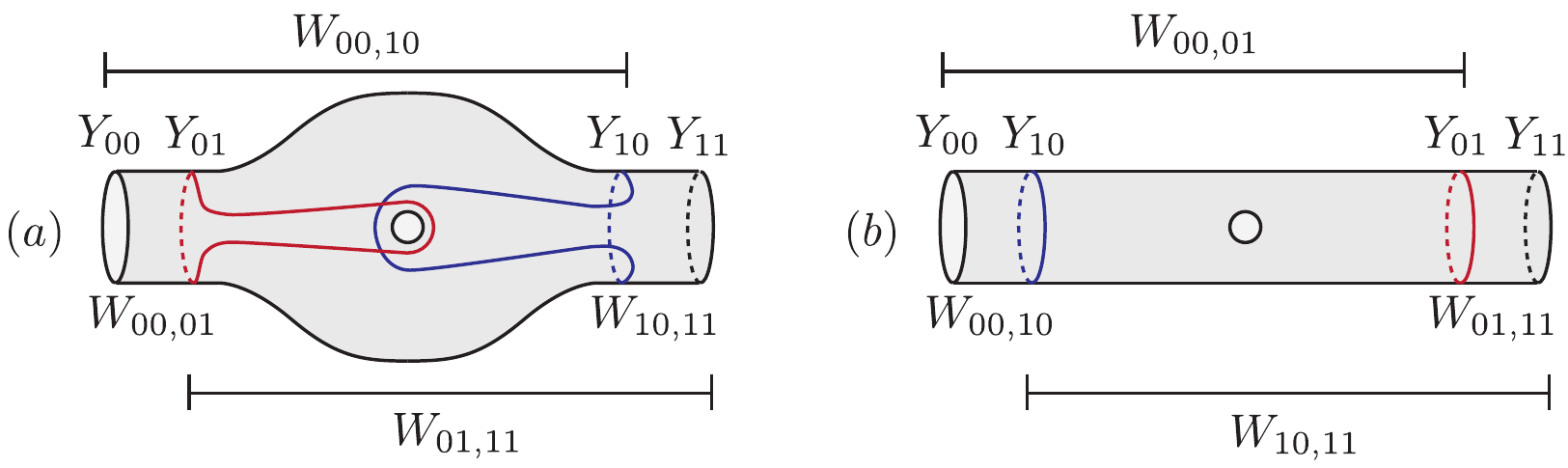}
\caption[The cobordisms used to construct the $U_\dagger$ map and prove $\hmhat$ is well-defined.]{Surface $(a)$ represents the cobordism used in the construction of the chain map $\tilde m(W) : \tilde C(Y_0) \to \tilde C(Y_1)$.  Surface $(b)$ represents the cylindrical cobordism used to prove that $\hmhat(Y)$ is well-defined, as discussed at the end of this section.}
\label{fig:uchain}
\end{figure}

The construction of the map $\hmhat(W) : \hmhat(Y_0) \to \hmhat(Y_1)$ induced by a cobordism $W$ is essentially the same as the $l = 1$ case of the $\hmhat$ spectral sequence to follow, but we describe it separately for concreteness and to motivate what follows.  First fix a small ball in the interior of $W$ (which we subsequently excise).  We relabel the ends of $W$ as $Y_{00}$ and $Y_{11}$ and embed a second copy of each in the interior of $W$ as follows (see Figure \ref{fig:uchain}a).  Consider a path $\gamma$ in $W$ from $Y_{00}$ to $Y_{11}$ such that a small tubular neighborhood $\nu(\gamma)$ of the path contains the ball.  $Y_{01}$ is obtained by taking a parallel copy of $Y_{00}$ just inside the boundary and smoothly pushing the region in $\nu(\gamma)$ past the ball, so that cutting along $Y_{01}$ leaves the ball in the first component $W_{00,01} \cong Y_0 \times [0,1]$.  Similarly, $Y_{10}$ is obtained by taking a parallel copy of $Y_{11}$ near the boundary and smoothly pushing the region inside $\nu(\gamma)$ inward past the ball, so that cutting along $Y_{10}$ leaves the ball in the second component $W_{10,11} \cong Y_1 \times [0,1]$.  Note also that both $W_{00,10}$ and $W_{01,11}$ are diffeomorphic to $W$.

The intersection of $Y_{01}$ and $Y_{10}$ is modeled on $S^2 \times \{0\} \times \{0\} \subset S^2 \times (-\epsilon, \epsilon) \times (-\epsilon, \epsilon)$, so we can choose the metric on $W$ to be cylindrical near both interior hypersurfaces.  Consider the interval of metrics $\tilde P_{00,11} = [-\infty, \infty]$ which stretches from $Y_{01}$ to $Y_{10}$.  We use this interval to define eight operators $H^{u*}_*(\cmf_{-1} \otimes \cdot)$ which count isolated trajectories in the moduli space

\begin{align*}
\Mod_{\tilde P_{00,11}} = \bigcup_{p \in \tilde P_{00,11}} \bigcup_z M_z(\mathfrak{a}, \mathfrak{c}_{-1}, W^{**}(p), \mathfrak{b})
\end{align*}
where 
\begin{align*}
M_z(\amf, W(-\infty)^*, \bmf) = \bigcup_{\cmf \in \cm(Y_{01})}  \bigcup_{z_1, z_2} M_{z_1}(\amf, \mathfrak{c}_{-1}, W_{00,01}^{**}, \cmf) \times M_{z_2}(\cmf, W_{01,11}^*, \bmf),
\end{align*}
and
\begin{align*}
M_z(\amf, W(\infty)^*, \bmf) = \bigcup_{\cmf \in \cm(Y_{10})}  \bigcup_{z_1, z_2} M_{z_1}(\amf, W_{00,10}^*, \cmf) \times M_{z_2} (\cmf , \cmf_{-1}, W_{10,11}^{**}, \bmf).
\end{align*}
We then define $\check{H}(U|W_{00,11})$ by the same expression as $\dd^{00}_{11}$ in \eqref{eqn:dij}, except that if $I$ ends in $0$ and $J$ ends in $1$, then $D^*_*({}^I_J)$ is replaced by $D^{u*}_*({}^I_J)(\mathfrak{c}_{-1} \otimes \cdot).$  So in full, we have
\begin{align*}
\check{H}(U|W_{00,11}) &= \left[\begin{array}{rr}
H^{uo}_o(\mathfrak{c}_{-1} \otimes \cdot) & H^{uu}_o(\mathfrak{c}_{-1} \otimes \esu(\cdot)) + \duo\bar{H}^{us}_{u}(\mathfrak{c}_{-1} \otimes \cdot)\\
H^{uo}_{s}(\mathfrak{c}_{-1} \otimes \cdot) & \bar{H}^{us}_s(\mathfrak{c}_{-1} \otimes \cdot) + H^{uu}_s(\mathfrak{c}_{-1} \otimes \esu(\cdot)) + \dus\bar{H}^{us}_{u}(\mathfrak{c}_{-1} \otimes \cdot) \end{array}\right] \\
&+ \left[\begin{array}{rr}
0 &  \muob\bar{m}^{us}_u({}^{00}_{01})(\mathfrak{c}_{-1} \otimes \cdot) + m^{uu}_o({}^{10}_{11})(\mathfrak{c}_{-1} \otimes \nsuc(\cdot)) \\
0 &  \musb\bar{m}^{us}_u({}^{00}_{01})(\mathfrak{c}_{-1} \otimes \cdot) + m^{uu}_s({}^{10}_{11})(\mathfrak{c}_{-1} \otimes \nsuc(\cdot)) \end{array}\right]
\end{align*}
From this perspective, the differentials on $\tilde{C}(Y_0)$ and  $\tilde{C}(Y_1)$ are
$$\tilde \partial (Y_0)
= \left[\begin{array}{cc} \check\partial(Y_{00}) & 0 \\ \check{m}(U | W_{00,01}) & \check\partial(Y_{01})\end{array}\right] \quad \text{and} \quad \tilde \partial (Y_1) 
= \left[\begin{array}{cc} \check\partial(Y_{10}) & 0 \\ \check{m}(U | W_{10,11}) & \check\partial(Y_{11})\end{array}\right],$$
respectively.
Finally, the map $\tilde{m}(W): \tilde{C}(Y_0) \to \tilde{C}(Y_1)$ is defined by
$$\tilde{m}(W) = \left[\begin{array}{cc} \check{m}(W_{00,10}) & 0 \\ \check{H}(U|W_{00,11}) & \check{m}(W_{01,11}) \end{array}\right].$$

We now turn to the general construction of the total complex underlying the $\hmhat$ version of the link surgery spectral sequence.  In this section, we will denote this complex by $(\underline{X}, \underline D)$, though in other sections we will return to the notation $(X,\dd)$ when it is clear from context which version is intended.  The same goes for the pages $\underline E^i,$ etc.

\begin{figure}[htp]
\centering
\includegraphics[width=160mm]{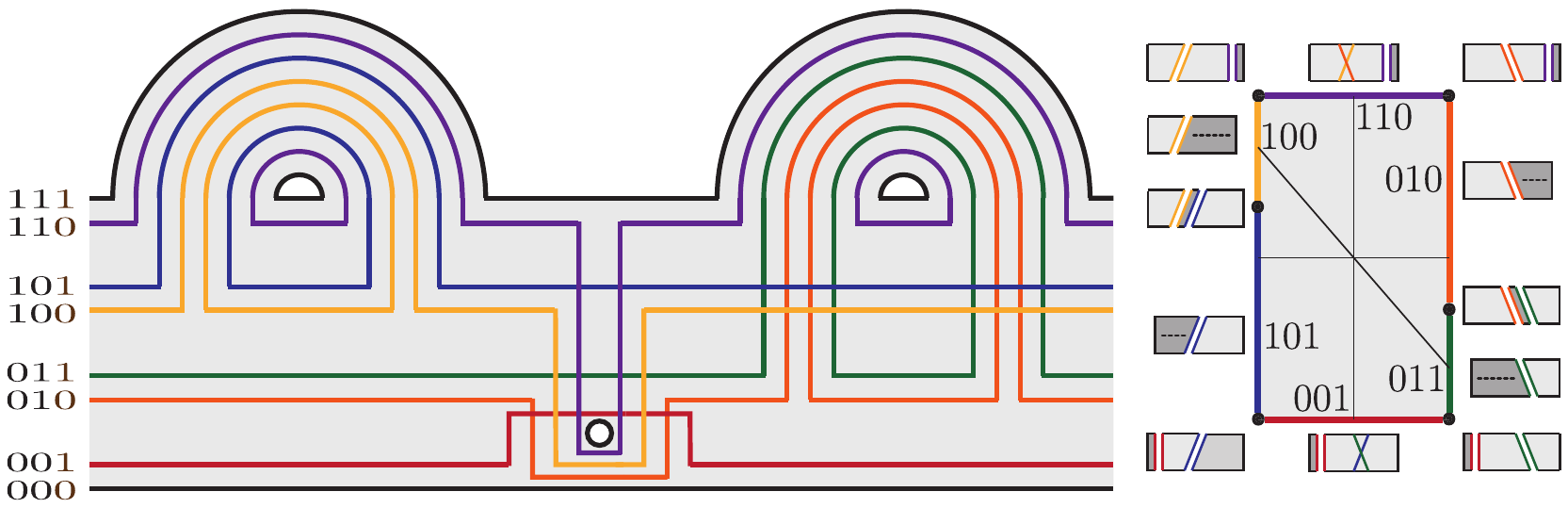}
\caption[The cobordism and hexagon of metrics used to construct the $\hmhat$ version of $\{0,1\}^2$.]{At left, we have the half-dimensional diagram of the cobordism $W$ with a small ball removed in the case $l = 2$.  For each $I \in \{0,1\}^2$, the pair $Y_{I\times\{0\}}$ and $Y_{I\times\{1\}}$ bound a cylindrical cobordism containing the ball.  At right, we have drawn the corresponding hexagon $\underline P^{000}_{111}$ so that increasing the vertical coordinate is suggestive of translating the sphere through $W$.  The small figures at the vertices and edges illustrate the metric degenerations, read as composite cobordisms from left to right.  In each, the component containing the sphere is more heavily shaded.}
\label{fig:intersectionsu}
\end{figure}

Given an $l$-component framed link $L \subset Y$, we embed a small ball $D^4$ in the interior of $(Y - \nu(L)) \times [0,1] \subset W$, centered at a point $\{x\} \times \{t\}$.  Next we relabel the incoming end $Y_{\{0\}^l}$ as $Y_{\{0\}^l \times \{0\}}$ and every other $Y_I$ as $Y_{I \times \{1\}}$.  We then embed a second copy of $Y_{\{0\}^l \times \{0\}}$, labeled $Y_{\{0\}^l \times \{1\}}$, just above the first, modified so that it now passes above the ball.  Finally, we embed a second copy of each $Y_{I \times \{1\}}$, labeled $Y_{I \times \{0\}}$, just below the first, modified so that it now passes below the ball, using the path $\{x\} \times [0,1]$ as a guide.  See Figure \ref{fig:intersectionsu} for the case $l = 2$.  We now have an embedded hypersurface $Y_I$ for each $I \in \{0,1\}^{l+1}$.  Furthermore, the intersection data is precisely what we expect for this hypercube, namely that $Y_I$ intersects $Y_J$ if and only if $I$ and $J$ are not ordered.  Therefore, given any $I < J$ we may construct a family of metrics on $\wij$ parameterized by a permutohedron $\underline P^I_J$ of dimension $w(J) - w(I) - 1$.

Now fix a regular metric and perturbation on the cobordism $W$ which is cylindrical near every hypersurface $\yi$ and round near $S^3$.  We will define a complex
\begin{align}
\label{eqn:xu}
\underline X = \bigoplus_{I\in \{0,1\}^l \times \{0,1\}} \cm(\yi),
\end{align}
where the differential $\underline{D}: \underline X \to \underline X$ is the sum of components $\underline D^I_J: \cm(\yi) \to \cm(\yj)$ over all $I \leq J$.  We have set things up so that the ball is contained in $\wij$ if and only if $I$ ends in $0$ and $J$ ends in $1$.  So when $I$ and $J$ end in the same digit, the operators $D^*_*({}^I_J)$ may be defined exactly as before (see \eqref{eqn:thed}).  In the other case, we construct operators $D^{u*}_*({}^I_J)(\mathfrak{c}_{-1} \otimes \cdot)$ using moduli spaces $M_z(\amf, \cmf_{-1}, W_{IJ}^{**}, \bmf)_{\underline P_{IJ}}$ which are defined by slightly modifying the definition of the moduli spaces  $\Modzpij$.  Namely, if $p \in \underline P^I_J $ is in the interior of the face $I_1 < I_2 < \cdots < I_{q-1}$, with the last digit changing between $I_k$ and $I_{k+1}$, then an element of $M_z(\amf,\cmf_{-1},\wij(p)^{**},\bmf)$ is a $q$-tuple
$$(\gamma_{01}, \gamma_{12}, \dots, \gamma_{q-1\,q})$$
as before except that
\begin{align*}
\gamma_{k\,k+1} & \in M(\amf_j, \cmf_{-1}, W_{I_k I_{k+1}}^{**}(p),\amf_{k+1}).
\end{align*}
Let $\underline D^*_*({}^I_J)$ be synonymous with $D^*_*({}^I_J)$ if $I$ and $J$ end in the same digit, and with $D^{u*}_*({}^I_J)(\mathfrak{c}_{-1} \otimes \cdot)$ otherwise.  Similar remarks apply to $\bar D^*_*({}^I_J)$ and $\bar D^{u*}_*({}^I_J)(\mathfrak{c}_{-1} \otimes \cdot)$.
We then define $\underline D^I_J: \cm(\yi) \to \cm(\yj)$ by precisely the same expression as \eqref{eqn:dij}, with each $D$ underlined.

The proof that $\underline D^2 = 0$ goes along familiar lines.   The operators $A^*_*({}^I_J)$ may be defined exactly as before when $I$ and $J$ end in the same digit.  When $I$ ends in $0$ and $J$ ends in $1$, we define operators $A^{u*}_*({}^I_J)(\mathfrak{c}_{-1} \otimes \cdot)$ which count ends of 1-dimensional moduli spaces $M_z^+(\amf, \cmf_{-1}, W_{IJ}^{**}, \bmf)_{\underline P_{IJ}}$, which in turn are defined by slightly modifying the definition of the moduli spaces  $\Modzpijplus$ in the same manner as above.  As before, these operators all vanish.  Now let $\underline A^*_*({}^I_J)$ be synonymous with $A^*_*({}^I_J)$ if $I$ and $J$ end in the same digit, and with $A^{u*}_*({}^I_J)(\mathfrak{c}_{-1} \otimes \cdot)$ otherwise.  Similar remarks apply to $\bar A^*_*({}^I_J)$ and $\bar A^{u*}_*({}^I_J)(\mathfrak{c}_{-1} \otimes \cdot)$.  We then define $\underline A^I_J: \cm(\yi) \to \cm(\yj)$ by precisely the same expression as \eqref{eqn:aij}, with each $D$ and $A$ underlined.

\begin{lemma} 
\label{lem:aij3}
$\underline A^I_J$ is equal to the component of $\underline D^2$ from $\cm(\yi)$ to $\cm(\yj)$:
$$\underline A^I_J = \sum_{I \leq K \leq J} \underline D^K_J \underline D^I_K.$$
Thus, $\underline D$ is a differential.  
\end{lemma}
\begin{proof}
Recall that $\cmf_{-1}$ is an unstable reducible and that there are no isolated trajectories from $\mathfrak{c}_{-1}$ to any other $\mathfrak{c}_i \in \mathfrak{C}(S^3)$.  It follows that 1-dimensional moduli spaces $M_z^+(\amf, \cmf_{-1}, W_{IJ}^{**}, \bmf)_{\tilde P_{IJ}}$ have the same types of ends as $M_z^+(\amf, W_{IJ}^*, \bmf)_{P_{IJ}}$, as described in Lemma \ref{lem:ends}.  Similarly, $M_z^{\text{red}+}(\amf, \cmf_{-1}, W_{IJ}^{**}, \bmf)_{\underline P_{IJ}}$ has the same types of ends as $M_z^{\text{red}+}(\amf, W_{IJ}^*, \bmf)_{\underline P_{IJ}}$.  Now we simply repeat the proof of Lemma \ref{lem:aij} with everything underlined.
\end{proof}

\begin{remark}
\label{rem:chain}
The case $l = 0$ shows that $U_\dagger : \cm(Y_0) \to \cm(Y_0)$ is a chain map.  The case $l = 1$ shows that $\tilde{m}(W): \cmhat(Y_0) \to \cmhat(Y_1)$ is a chain map when $W$ an elementary 2-handle cobordism, and goes through without change for arbitrary cobordisms.
\end{remark}

In order to interpret Lemma \ref{lem:aij3} as a result in hat theory, we collapse $$\underline X = \bigoplus_{I \in \{0,1\}^l} \cmhat(\yi)$$ along the final digit, with $\underline D$ given by the sum of maps $\tilde D^I_J = \tilde C(\yi) \to \tilde C(\yj)$  where
\begin{align}
\label{eqn:tilde}
\tilde D^I_J = \left[\begin{array}{cc} \dd^{I\times\{0\}}_{J \times \{0\}} & 0 \\ \underline D^{I\times\{0\}}_{J \times \{1\}} & \dd^{I\times\{1\}}_{J \times \{1\}}\end{array}\right].
\end{align}
Define the horizontal weight $\underline w(I)$ of a vertex $I$ to be the sum of all but the final digit.  Filtering $(\underline X, \underline D)$ by $\underline w$, we obtain the $\hmhat$ version of the link surgery spectral sequence.  In particular,$$\underline E^1 = \bigoplus_{I\in \{0,1\}^l} \hmhat(\yi)$$ and the $\underline d^1$ differential is given by $$\underline d^1 = \bigoplus_{w(J) - w(I) = 1} \hmhat(\wij).$$

In order to identify $\underline E^\infty$ with $\hmhat(Y_\infty)$, we expand to the larger complex
\begin{align}
\label{eqn:xut}
\underline {\widetilde X} = \bigoplus_{I\in \{0,1, \infty\}^l \times \{0,1\}} \cm(\yi),
\end{align}
where we have again relabeled the $Y_I$ and embedded a second copy of each which passes on the opposite side of the ball.  These hypersurfaces all avoid the auxiliary $S^3$ hypersurfaces which are confined in the handles.  The intersection data is as predicted for the shape of the lattice by Theorem \ref{thm:graph}, so we may build families of metrics parameterized by graph associahedra which define maps $\underline D^I_J$ and $\underline A^I_J$.  The auxiliary hypersurfaces still cut out $\cptb - \text{int}(D^4)$, so the corresponding facets contribute vanishing operators as before.  We conclude from \ref{lem:aij3} that $(\underline {\widetilde X}, \underline D)$ forms a complex.

\begin{figure}[htp]
\centering
\includegraphics[width=155mm]{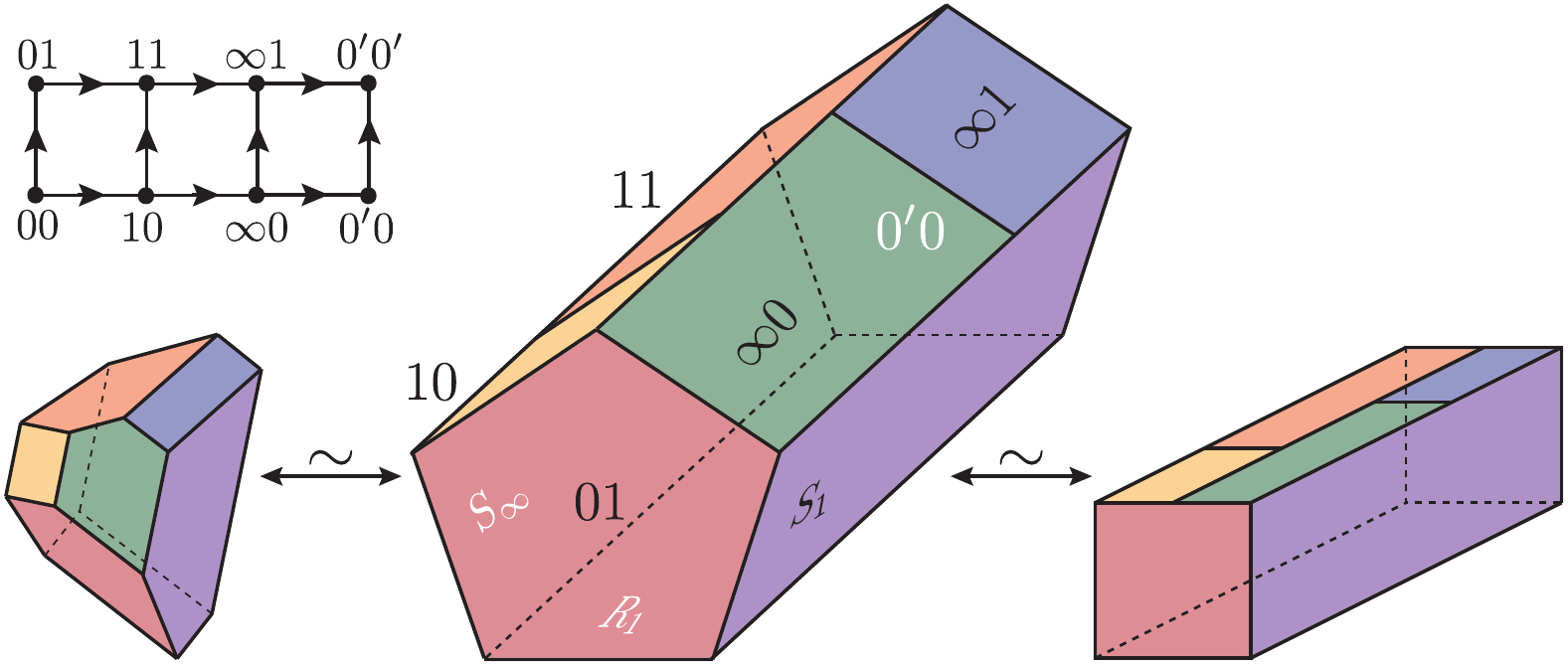}
\caption[The associahedron of metrics for $\{0,1,\infty,0'\} \times \{0,1\}$.]{The lattice $\{0,1,\infty,0'\} \times \{0,1\}$ corresponds to Stasheff's polytope, drawn at center and right so that the depth coordinate is suggestive of translating the sphere through $W$.  Recall that the same polytope is associated to the lattice $\{0,1,\infty\} \times \{0,1,\infty \}$ in Figure \ref{fig:assoc}, redrawn at left.  In fact, by Theorem \ref{thm:graph}, an associahedron arises whenever the lattice is a product of at most two chains.  The map corresponding to $K_5$ above is a null-homotopy for the sum of the maps associated to the faces.  Those associated to $S_1$ and $S_2$ vanish as before, leaving a nine-term identity on the chain level.}
\label{fig:assocu}
\end{figure}

Finally, we turn to the surgery exact triangle.  Recall the hypersurfaces $Y_0$, $Y_1$, $Y_\infty$, $Y_{0'}$ and auxiliary $S_1$, $S_2$, and $R_1$.  After relabeling, these become $Y_{00}$, $Y_{10}$, $Y_{\infty 0}$, $Y_{0'1}$, $S_1$, $S_2$, and $R_1$, to which we add $Y_{01}$, $Y_{11}$, $Y_{\infty 1}$, and $Y_{0'0}$.  The nine hypersurfaces in the interior of $W$ intersect as predicted by the shape of the lattice $\{0,1,\infty, 0'\} \times \{0,1\}$, yielding the graph associahedron on a chain of length $4$, namely, the 3-dimensional associahedron $K_5$, shown in Figure \ref{fig:assocu}.

The map $\underline L^{00}_{0'1}: \cm(Y_{00}) \to \cm (Y_{0'1})$ associated to $R_1$ is given by the same expression as $\check{L}$, but with $L^*_* = D^{uu*}_*(\cmf_{-1} \otimes \bar{n}_u \otimes \cdot)$ and  $\bar{L}^*_* = \bar{D}^{uu*}_*(\cmf_{-1} \otimes \bar{n}_u \otimes \cdot)$, over the one-parameter family of metrics stretching from $Y_{01}$ to $Y_{0'0}$. This gives the identity 
\begin{align}
\label{eqn:lid}
\dd^{0'1}_{0'1} \underline L^{00}_{0'1} + \underline L^{00}_{0'1} \dd^{00}_{00} = \check L^{01}_{0'1} \underline D^{00}_{01} + \underline D^{0'0}_{0'1} \check L^{00}_{0'0},
\end{align}
where $\check L^{00}_{0'0}$ and $\check L^{01}_{0'1}$ are the analogues of $\check L$ corresponding to the four hypersurfaces $Y_I$ ending in $0$ and $1$, respectively, together with $S_1$, $S_2$, and $R_1$.  Similarly, $\underline{A}^0_{0^\prime}$ is modeled on the expression \eqref{eqn:exacta2} for $\check{A}^0_{0^\prime}$ and includes terms counting monopoles on cobordisms with four boundary components.  By the same reasoning as in the proof of Lemma \ref{lem:aij3}, the analogue of Lemma \ref{lem:aij2} goes through essentially unchanged (although with nearly twice as many terms), leading to the nine-term identity given by the lower left entry of
\begin{align*}
\left[\begin{array}{cc} \dd^{0'0}_{0'0} & 0 \\ \underline D^{0'0}_{0'1} & \dd^{0'1}_{0'1}\end{array}\right]
\left[\begin{array}{cc} \dd^{00}_{0'0} & 0 \\ \underline D^{00}_{0'1} & \dd^{01}_{0'1}\end{array}\right] +
\left[\begin{array}{cc} \dd^{00}_{0'0} & 0 \\ \underline D^{00}_{0'1} & \dd^{01}_{0'1}\end{array}\right]
\left[\begin{array}{cc} \dd^{00}_{00} & 0 \\ \underline D^{00}_{01} & \dd^{01}_{01}\end{array}\right] \hspace{31.45mm} \\
= \left[\begin{array}{cc} \dd^{10}_{0'0} & 0 \\ \underline D^{10}_{0'1} & \dd^{11}_{0'1}\end{array}\right]
\left[\begin{array}{cc} \dd^{00}_{10} & 0 \\ \underline D^{00}_{11} & \dd^{01}_{11}\end{array}\right] +
\left[\begin{array}{cc} \dd^{\infty 0}_{0'0} & 0 \\ \underline D^{\infty 0}_{0'1} & \dd^{\infty 1}_{0'1}\end{array}\right]
\left[\begin{array}{cc} \dd^{00}_{\infty0} & 0 \\ \underline D^{00}_{\infty 1} & \dd^{01}_{\infty 1}\end{array}\right] +
\left[\begin{array}{cc} \check L^{00}_{0'0} & 0 \\ \underline L^{00}_{0'1} & \check L^{01}_{0'1}\end{array}\right]
\end{align*}
The upper left and lower right identities are precisely those given by Lemma \ref{lem:aij2}.  Rewriting this identity via \eqref{eqn:tilde}, we have the $\hmhat$ analog of \eqref{eqn:aijlij1}:
\begin{align*}
\tilde D^{0^\prime}_{0^\prime} \tilde D^0_{0^\prime} + \tilde D^0_{0^\prime} \tilde D^0_0 = \tilde D^1_{0^\prime} \tilde D^0_1 + \tilde D^\infty_{0^\prime} \tilde D^0_\infty + \tilde L.
\end{align*}
The final map $\tilde L : \cmhat(Y_0) \to \cmhat(Y_{0'})$ is a chain map by \eqref{eqn:lid}.  Filtering the corresponding square complex $$Z = \bigoplus_{I \in \{00, 0'0, 01, 0'1\}} \cm(\yi)$$ by the second digit, and recalling that $\check L^{00}_{0'0}$ and $\check L^{01}_{0'1}$ are quasi-isomorphisms, we conclude that $H_*(Z) = 0$.  Therefore, $\tilde L$ is a quasi-isomorphism as well.  Now exactly the same algebraic arguments yield the surgery exact triangle, and more generally the full statement of the link surgery spectral sequence, for $\hmhat$.

The grading and invariance results of Section \ref{sec:linksurg2} readily extend to $\hmhat$ by viewing the underlying complex in terms of $\hmb$ as in \eqref{eqn:xu} and \eqref{eqn:xut}.  In this way, we may extend $\check \delta$ to a mod 2 grading on $\underline {\tilde X}$ using the same definition.  Since $U_\dagger$ cuts down the dimension of moduli spaces by two, the maps $\underline D^I_J$ on $\underline{\tilde X}$ obey the same mod 2 grading shift formula as the maps $\dij$ on $\widetilde{X}$.  In particular, Lemmas \ref{lem:byone} and \ref{prop:quasi} still apply, and when $l = 0$, $\check \delta$ and $\grt$ coincide on $\hmhat(Y)$.  The proof of Theorem \ref{thm:anal} regarding invariance also readily adapts to a version for $\hmhat$, using a cobordism $W$ with $2^{l+2}$ hypersurfaces (that is, two copies of each hypersurfaces in Figure \ref{fig:intersectionsu}), where we maintain the round metric on the boundary of the ball throughout.  The $l = 0$ case is pictured in (b) of Figure \ref{fig:uchain}, and implies that $\hmhat(Y)$ is well-defined.  Similarly, the $l=1$ case for a general cobordism implies that $\hmhat(W)$ is well-defined.  The composition law \eqref{eqn:complaw} follows from the $l=2$ case by stretching from $W_2 \circ W_1$ to $W \circ (Y_0 \times [0,1])$.

\subsection{Realizations of graph associahedra}

\begin{figure}[htp]
\centering
\includegraphics[width=150mm]{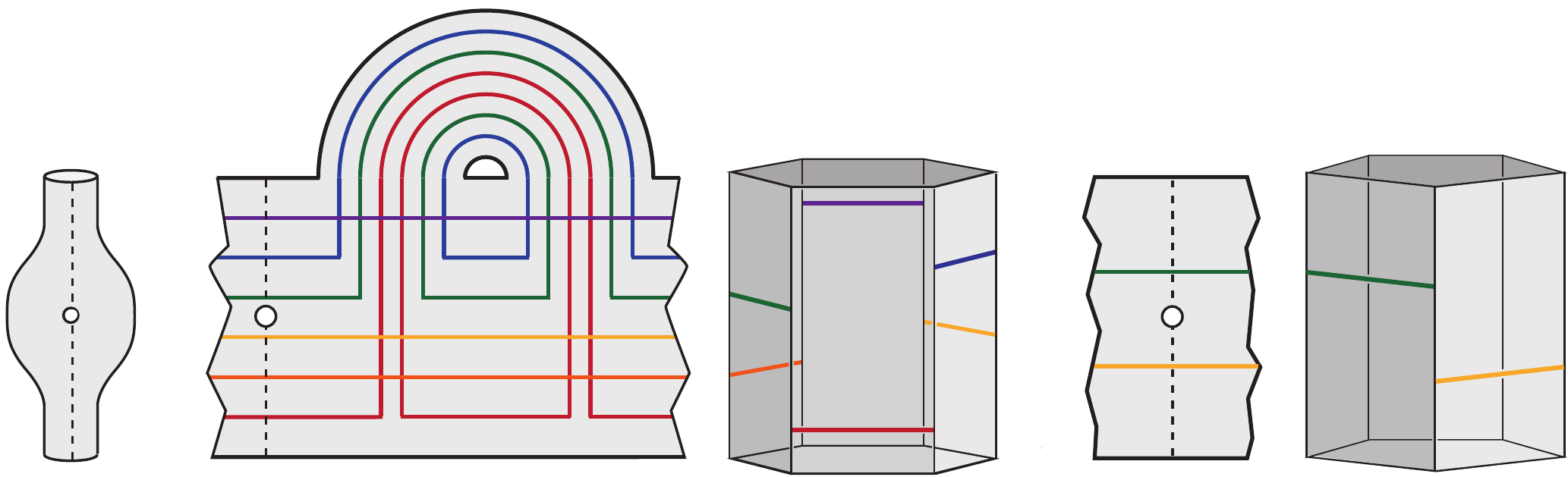}
\caption[Sliding the sphere to construct product-like realizations of graph associahedra.]{At left, we slide the sphere through the full cobordism from  Figures \ref{fig:intersections1} and \ref{fig:hexagon}.  Each time the sphere crosses an embedded hypersurface, we add a ridge to the corresponding lateral facet of $P_3 \times [0,1]$.  Once the sphere has completed its journey, we have a realization of $P_4$.  At right, we similarly slide the sphere through the full cobordism from Figure \ref{fig:intersections2}, adding ridges to $K_4 \times [0,1]$.  Once the sphere has completed its journey, we have $K_5$.}
\label{fig:hexvert}
\end{figure}

Recall the pentagonal realization of Stasheff's polytope at center in Figure \ref{fig:assocu}, and the inductive realization of $P_{l+1}$ as a refinement of $P_l \times [0,1]$ in the proof of Theorem \ref{thm:anal}.  Both of these realizations are motivated by the ``sliding-the-point'' proof of the naturality of the $U_\dagger$ action in Floer theory (though from our perspective we are sliding the sphere).  To see why, recall that to any product lattice $\Lambda$ we may associate a map $\dd$ whose longest component counts monopoles on $W$ over a family of metrics parameterized by a polytope $P_G$ (see Remark \ref{rem:genlat}).  We then expect the longest component of the homotopy which expresses the naturality of the $U_\dagger$ action with respect to $\dd$ to count monopoles over a family of metrics parameterized by $P_G \times [0,1]$, where the latter coordinate slides the sphere through $W$ (see Figure \ref{fig:hexvert}).  In Section \ref{sec:umap}, we instead embedded two copies of each hypersurface in $W$, because this approach fit more cleanly with our earlier constructions.

\newpage

More generally, consider any realization of the graph associahedron $P_G$ of dimension $n-1$ associated to a lattice $\Lambda$.  We may realize the graph associahedron associated to $\Lambda \times \{0,1\}$ as a refinement of $P_G \times [0,n]$.  Namely, for each internal vertex $I$ of $\Lambda$, we add a ridge at height $w(I)$ to the corresponding lateral facet of $P_G \times [0,n]$.  We ignore auxiliary hypersurfaces, as the sphere never passes through them.  In the two cubical realizations below, we have applied this construction twice, starting from a graph associahedron which is geometrically an interval.

\begin{figure}[htp]
\centering
\includegraphics[width=140mm]{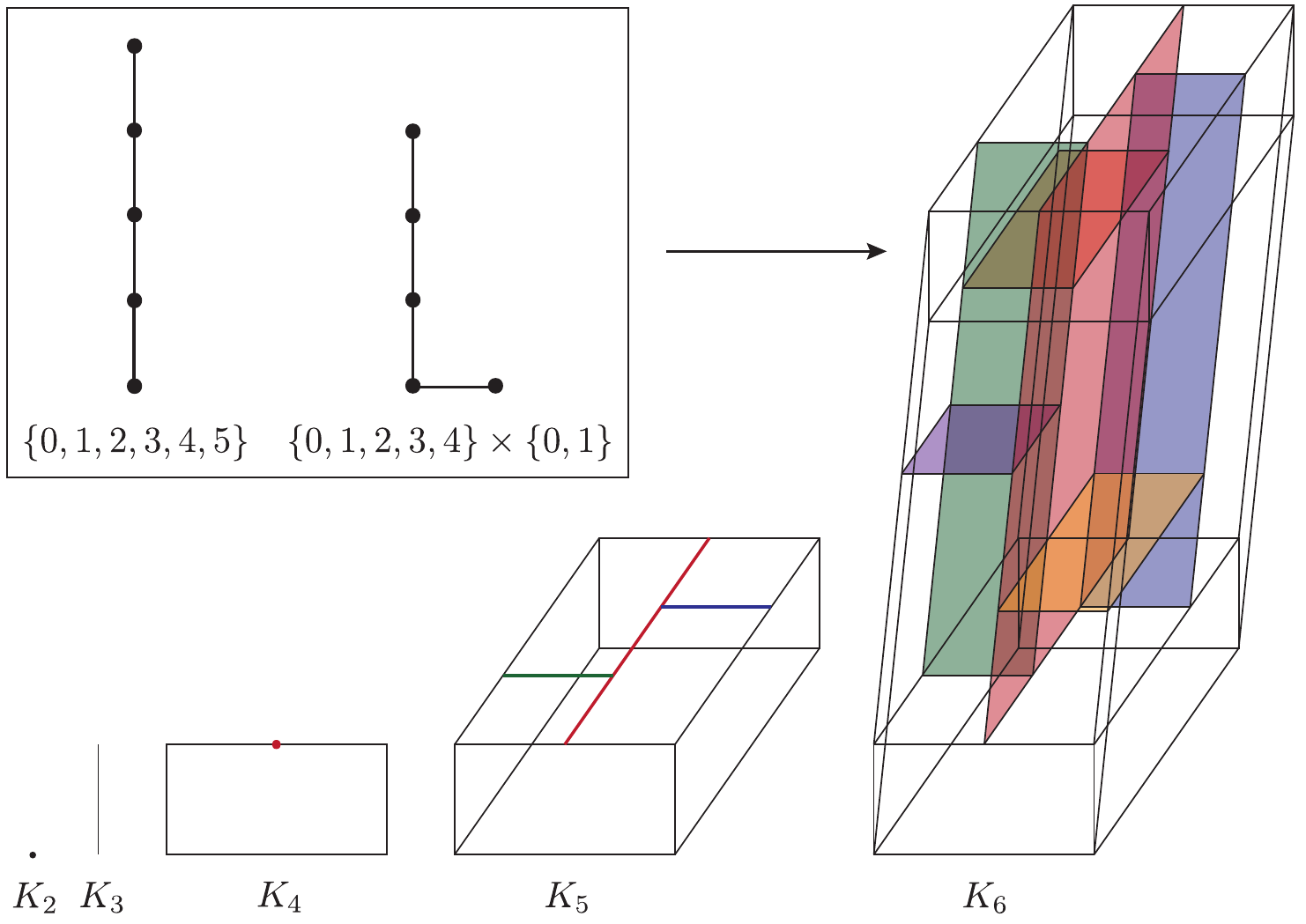}
\caption[Realizations of the associahedra.]{We realize $K_{n+2}$ by adding $\left( {}^n _2 \right)$ ridges with integral vertices to the facet $\{1\} \times [0,2] \times \cdots \times [0,n]$ of the hypercube $[0,1] \times [0,2] \times \cdots \times [0,n]$.  Similarly, $P_{n+1}$ is obtained by adding $2^{n+1} - 2(n+1)$ ridges to the hypercube.}
\label{fig:assocr}
\end{figure}

Note that the lattices $\{0,1,\dots,n\}$ and $\{0,1,\dots,n-1\} \times \{0,1\}$ both correspond to the graph consisting of a path of length $n-1$.  By the above remarks, this fact yields a simple realization of the associahedron $K_{n+2}$ by adding ridges to $K_{n+1} \times [0,n]$.  We obtain a similar realization of each permutohedron $P_n$ from the fact that $\{0,1\}^n = \{0,1\}^{n-1} \times \{0,1\}$.  If we build these realizations inductively, starting from the point, then each is naturally a refinement of the hypercube (see Figure \ref{fig:assocr}), and we may arrange that $P_{n+1}$ refines $K_{n+2}$ as well.  Upon sharing these realizations with experts, we were directed to similar ones in \cite{umb}, which also refine hypercubes but have vastly different origins.  For $K_n$ as a convex hull, see \cite{lod} and \cite{d}.

\newpage


\section{Khovanov homology and branched double covers}
\label{sec:khovanov}

We now turn to the construction of the spectral sequence from  $\kh(L)$ to $\hmhat(-\Sigma(L))$ in Theorem $\ref{thm:1}$.  We would like to interpret the $(E^1, d^1)$ page as the result of applying a much simpler functor than $\hmhat$ to the branched double cover of the hypercube of resolutions.  To this end, we prove the $\hmb$ and $\hmhat$ analogs of Proposition 6.2 in \cite{osz12}.

\begin{proposition}
\label{prop:s1s2}
Let $Y \cong \#^k(S^1 \times S^2)$.  Then, $\hmb(Y)$ is a rank one, free module over the ring $$\Lambda^*H_1(Y) \textstyle\bigotimes \ftwo[U_\dagger^{-1}, U_\dagger] / \ftwo[U_\dagger],$$
generated by some class $\Theta \in \cm(Y)$, and entirely supported over the torsion spin$^c$-structure. Moreover, if $K \subset Y$ is a curve which represents one of the circles in one of the $S^1 \times S^2$ summands, then the three-manifold $Y' = Y_0(K)$ is diffeomorphic to $\#^{k-1}(S^1 \times S^2)$, with a natural identification
$$ \pi: H_1(Y) / [K] \longrightarrow H_1(Y').$$
Under the cobordism $W$ induced by the two-handle, the map
$$\hmb(W) : \hmb(Y) \longrightarrow \hmb(Y')$$
is specified by 
\begin{align}
\label{eqn:merge}
\hmb(W)(\xi \otimes U_\dagger^n \cdot \Theta) = \pi(\xi) \otimes U_\dagger^n \cdot \Theta',
\end{align}
where here $\Theta'$ is some fixed generator of $\hmb(Y')$, and $\xi$ is any element of $\Lambda^*H_1(Y)$.  Dually, if $K \subset Y$ is an unknot, then $Y'' = Y_0(K) \cong \#^{k+1}(S^1 \times S^2)$, with a natural inclusion
$$i : H_1(Y) \longrightarrow H_1(Y'').$$
Under the cobordism $W'$ induced by the two-handle, the map
$$\hmb(W) : \hmb(Y) \longrightarrow \hmb(Y'')$$
is specified by
\begin{align}
\label{eqn:split}
\hmb(W')(\xi \otimes U_\dagger^n \cdot \Theta) = (\xi \wedge [K'']) \otimes U_\dagger^n \cdot \Theta'',
\end{align}
where here $[K''] \in H_1(Y'')$ is a generator in the kernel of the map $H_1(Y'') \to H_1(W')$.
\end{proposition}

\begin{proof}
The description of $\hmb(Y)$ holds whenever $Y$ has a metric of strictly positive scalar curvature, as shown in Proposition 36.1.3 of \cite{km}.  Note also that, by Theorem 3.4.4 of \cite{km}, the maps $\hmb(W)$ and $\hmb(W')$ behave naturally with respect to the above module structures.

We first consider the case where $K$ represents a circle in one of the $S^1 \times S^2$ factors.  As $[K]$ is null-homologous in $W$, we have $$\hmb(W)(([K] \wedge \xi) \otimes U_\dagger^n \cdot \Theta) = 0 \cdot \hmb(W)(\xi \otimes U_\dagger^n \cdot \Theta) = 0.$$  On the other hand,  since $Y_1(K) \cong Y'$, the corresponding surgery exact triangle splits (by comparing the ranks) as
$$0 \longrightarrow \hmb(Y') \longrightarrow \hmb(Y) \xrightarrow{\hmb(W)} \hmb(Y') \longrightarrow 0.$$
Therefore,  $\hmb(W)$ is surjective, which together with naturality, forces \eqref{eqn:merge}.

Now suppose $K$ is an unknot.  As $[K'']$ is null-homologous in $W'$, we have $$[K''] \cdot \hmb(W')(\xi \otimes U_\dagger^n \cdot \Theta) = \hmb(W')((0 \wedge \xi) \otimes U_\dagger^n \cdot \Theta) = 0.$$  So the image of $\hmb(W')$ is contained in $[K''] \wedge \hmb(Y'')$.  Now $Y_1(K) \cong Y$ and the short exact sequence reads 
$$0 \longrightarrow \hmb(Y) \xrightarrow{\hmb(W')} \hmb(Y'') \longrightarrow \hmb(Y) \longrightarrow 0.$$
Therefore,  $\hmb(W')$ is injective, which together with naturality, forces \eqref{eqn:split}.
\end{proof}

\begin{corollary}
\label{cor:s1s2}
Let $Y \cong \#^k(S^1 \times S^2)$.  Then, $\hmhat(Y)$ is a rank one, free module over the ring $\Lambda^*H_1(Y)$.  After removing $U_\dagger^n$ from \eqref{eqn:merge} and \eqref{eqn:split},  Proposition \ref{prop:s1s2} holds for  $\hmhat(Y)$.
\end{corollary}
\begin{proof}
The $U_\dagger$ map is surjective on the level of homology, so we have the short exact sequence
$$0 \longrightarrow \hmhat(Y) \longrightarrow \hmb(Y) \xrightarrow{U_\dagger} \hmb(Y) \longrightarrow 0.$$  In particular, $\hmhat(Y)$ is identified with $\text{Ker}(U_\dagger)$, which is the image of $\Theta$ under the action of $\Lambda^*H_1(Y) \otimes 1$.  Under this identification, the map $\hmhat(W)$ coincides with the restriction of $\hmb(W)$ to $\text{Ker}(U_\dagger)$.  The same holds in the case of $W'$.
\end{proof}

\begin{figure}[htp]
\centering
\includegraphics[width=100mm]{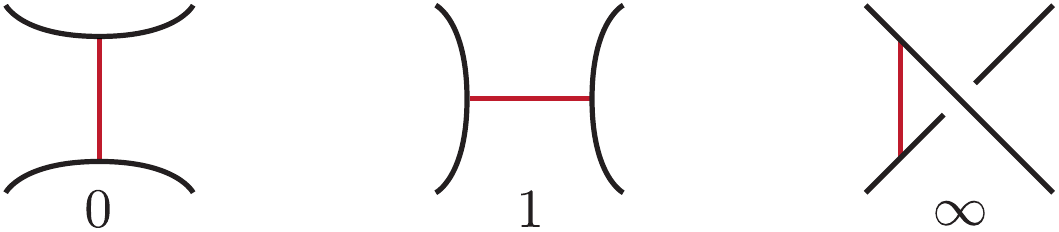}
\caption[The short arc between the strands of a resolved crossing.]{We have one short arc between the two strands near each crossing, in both the original diagram and its resolutions.}
\label{fig:conv}
\end{figure}

We now construct the $\hmhat$ spectral sequence associated to a link $L \subset S^3$.  We first fix a diagram $\mathcal{D}$ with $l$ crossings.  Following Section 2 of \cite{osz12}, we associate to $\mathcal{D}$ a framed link $\mathbb{L} \subset -\Sigma(L)$ to which we will apply the link surgery spectral sequence.  First, in a small ball $B_i$ about the crossing $c_i$, place an arc with an end on each strand as shown in the $\infty$ resolution of Figure \ref{fig:conv}.  Each of these arcs lifts to a closed loop $\mathbb{K}_i$ in the branched double cover $\Sigma(L)$, giving the components of a link $\mathbb{L} \subset \Sigma(L)$.  Note that all of the resolutions of $\mathcal{D}$ agree outside of the union of the $B_i$.  Furthermore, the branched double cover of $B_i$ over the two unknotted strands of $\mathcal{D}(I) \bigcap B_i$ is a solid torus, with meridian given by the preimage of either of the two strands pushed out to the boundary of $B_i$.  So for each $I \in \{0,1,\infty\}^l$, we may identify $\Sigma(\mathcal{D}(I)) - \nu(\mathcal{D}(I))$ with $\Sigma(L) - \nu(\mathbb{L})$.

In this way, for each crossing $c_i$, we obtain a triple of curves ($\lambda_i, \lambda_i + \mu_i, \mu_i)$ in the corresponding boundary component of $\Sigma(L) - \nu(\mathbb{L})$, which represent meridians of the fillings giving the branched covers of the 0-, 1-, and $\infty$-resolutions at $c_i$, respectively.  In this cyclic order, the curves may be oriented so that the algebraic intersection number of consecutive curves is +1.  We change this to $-1$ by flipping the orientation on the branched double cover (whereas in \cite{osz12} this is done by replacing $L$ with its mirror).  In the language of \cite{osz12}, each triple $(\lambda_i, \lambda_i + \mu_i, \mu_i)$ forms a triad.  From our 4-manifold perspective, this is precisely the condition that each 2-handle in a stack is attached to the previous 2-handle using the -1 framing with respect to the cocore (see the discussion preceding Theorem \ref{thm:linksurg}). From either perspective, framing $\mathbb{K}_i$ by $\lambda_i$, we are in precisely the setup of the link surgery spectral sequence, with $\yi = -\Sigma(\mathcal{D}(I))$ for all $I \in \{0,1,\infty\}^l$.  Now, using Corollary \ref{cor:s1s2}, the argument in \cite{osz12} may be repeated verbatim to show that the complexes $(E^1, d^1)$ and $\ckh(L)$ are isomorphic, and therefore that $E^2$ is isomorphic to $\kh(L)$.

\begin{figure}[htp]
\centering
\includegraphics[width=157mm]{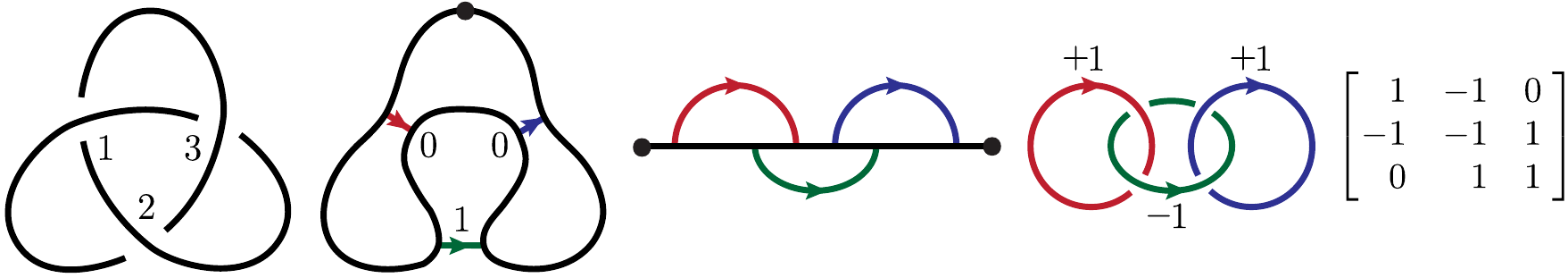}
\caption[A surgery diagram for the branched double cover of the trefoil.]{At left, we number the crossings in a diagram of the right-handed trefoil $\mathcal{T}$.  The resolution $\mathcal{D}(010)$ has one circle, and one (arbitrarily-oriented) arc for each crossing.  We then cut the circle open at the dot and stretch it out to a line, dragging the arcs along for the ride.  Reflecting each arc under the line yields the framed link $\mathbb{L}' \subset S^3$ and linking matrix at right.  Surgery on $\mathbb{L}'$ gives $-\Sigma(\mathcal{T})$, which is the lens space $-L(3,1)$.}
\label{fig:trefoil}
\end{figure}

As an alternative to the argument in \cite{osz12}, we now present a more global description of the isomorphism $E^1 \cong \ckh(L)$, taking advantage of the construction in \cite{b1} of a framed link $\mathbb{L}' \subset S^3$ (there denoted $\mathbb{L}$) which gives a surgery diagram for the branched double cover $-\Sigma(L)$.  Note that our convention in \cite{b1} on the orientation of surgery diagrams is precisely the opposite of the more sensible one we are using now (in this paper, the branched double cover of the right-handed trefoil is given by +3 surgery on an unknot, and we plan to revise \cite{b1} accordingly).  This construction of $\mathbb{L}'$ is illustrated in Figure \ref{fig:trefoil} for the standard diagram of the right-handed trefoil $\mathcal{T}$.  We first fix a vertex $I^* = (m^*_1,\dots,m^*_n) \in \{0,1\}^l$ for which the resolution $\mathcal{D}(I^*)$ consists of one circle.  The link $\mathbb{L}'$ is obtained as the preimage of the corresponding red arcs in Figure \ref{fig:conv}, where the component $\mathbb{K}'_i$ is given the framing $\lambda_i = (-1)^{m_i^*}$.  Note that the link $\mathbb{L} \subset -\Sigma(L)$ is represented in this surgery diagram by the framed push-off of $\mathbb{L}'$.

\begin{figure}[htp]
\centering
\includegraphics[width=150mm]{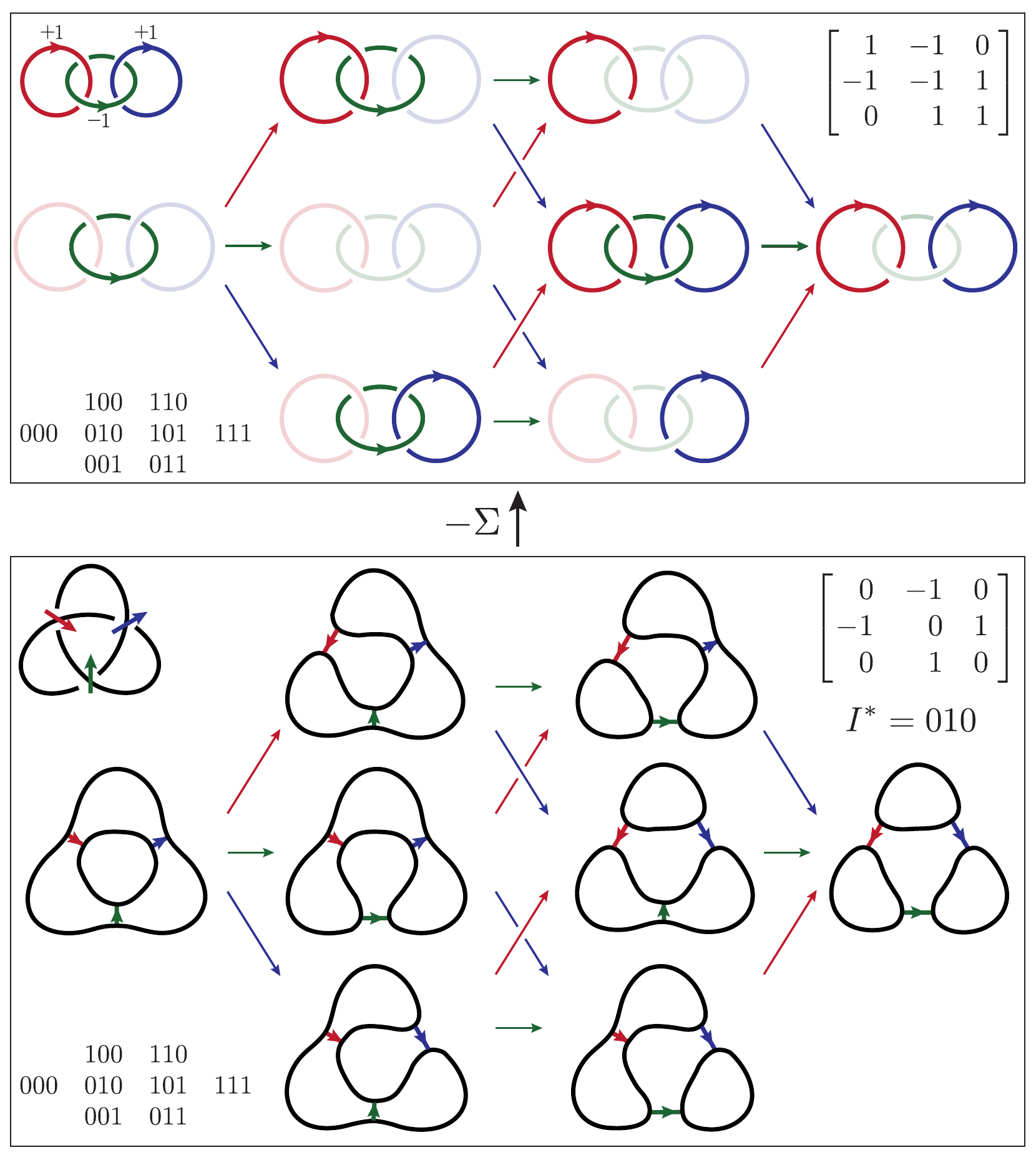}
\caption[The branched double cover of the Khovanov cube of resolutions of the trefoil.]{Continuing from Figure \ref{fig:trefoil}, above we obtain the cube of surgery diagrams $\mathbb{L}'_I$ for the right-handed trefoil $\mathcal{T}$.  All solid components are 0-framed, while all faded components are $\infty$-framed.
The link surgery cube of 3-manifolds $\yi$ and 4-dimensional 2-handle cobordisms $\wij$ above is the branched double cover of the Khovanov cube of 1-manifolds $\mathcal{D}(I) \subset S^3$ and 2-dimensional 1-handle cobordisms $F_{IJ} \subset S^3 \times [0,1]$ below.  At the upper right of each box, we record the associated numerical data, which determines $\mathbb{L}'$ as described in \cite{b1} and Remark \ref{rem:encode} below.  The orientation conventions for odd Khovanov homology are also described in \cite{b1}, but over $\ftwo$ these orientations are irrelevant.}
\label{fig:trefcube}
\end{figure}

For each $I = (m_1,\dots,m_n) \in \{0,1\}^l$, let $\mathbb{L}'_I$ be the link $\mathbb{L}'$ with framing modified to
\begin{align*}
\lambda_i &= \left\{
    \begin{array}{l}
        \infty \quad \text{if } m_i = m_i^*\\
        \hspace{1mm} 0 \hspace{.9mm} \quad \text{if } m_i \neq m_i^*
    \end{array} \right.
\end{align*}
on $\mathbb{K}'_i$.  Then $\mathbb{L}'_I$ gives a surgery diagram for $\yi = -\Sigma(\mathcal{D}(I)) \cong \#^k (S^1 \times S^2)$, where the resolution $\mathcal{D}(I)$ consists of $k+1$ circles.  This is illustrated in Figure \ref{fig:trefcube} for the trefoil $\mathcal{T}$.  Furthermore, each elementary 2-handle cobordism $\wij = -\Sigma(F_{IJ}) :-\Sigma(\mathcal{D}(I)) \to -\Sigma(\mathcal{D}(J))$ is given explicitly by 0-surgery on either $\mathbb{K}_i$ or its meridian $x_i$, in the case where $m_i^* = 0$ or $1$, respectively.  Here $F_{IJ} \subset S^3 \times [0,1]$ is the saddle-like surface cobordism from $\mathcal{D}(I)$ to $\mathcal{D}(J)$.

Motivated by Corollary \ref{cor:s1s2} and Figure \ref{fig:trefcube}, we define a complex with underlying $\ftwo$-vector space
$$\widehat C(\mathcal{D}) = \bigoplus_{I \in \{0,1\}^l} \Lambda^*H_1(-\Sigma(\mathcal{D}(I)))$$
and differential $\hat\partial$ given by the sum of maps $\hat\partial^I_J$ over all immediate successors $I < J$ in $\{0,1\}^l$.  These in turn are defined by
\begin{align*}
\hat\partial^I_J(\xi) &= \left\{
    \begin{array}{l}
        \pi(\xi) \hspace{15mm} \text{if $K$ represents a circle factor,} \\
         \text{[}K''\text{]} \wedge i(\xi) \quad \text{if $K$ is an unknot,}
     \end{array} \right.
\end{align*}
where our notation is from Proposition \ref{prop:s1s2}.

The central point of \cite{b1}, reduced mod 2, is that we can take the complex $(\widehat C(\mathcal{D}), \hat \partial)$ as the definition of $\ckh(D)$.  Namely, $(\widehat C(\mathcal{D}), \hat \partial)$ is precisely the version of $\ckh(D)$ defined at the end of Section 3, using the identification of $H_1(-\Sigma(\mathcal{D}(I)))$ with $\widehat V (\mathcal{D}(I))$ at the end of Section 4.
There, the main tool is to present each group $H_1(-\Sigma(\mathcal{D}(I)))$ by the linking matrix of $\mathbb{L}'_I$, using the common basis given by the meridians $x_i$.
On the other hand, Corollary \ref{cor:s1s2} immediately implies that $(\widehat C(\mathcal{D}), \hat \partial)$ is isomorphic to the page $(E^1, d^1)$ of the link surgery spectral sequence for $\mathbb{L} \subset -\Sigma(L)$.  We conclude that $(E^1, d^1)$ and $\ckh(L)$ are isomorphic complexes, and therefore that $E^2$ is isomorphic to $\kh(L)$ as an $\ftwo$-vector space.  This establishes Theorem \ref{thm:1}.

\begin{remark}
\label{rem:encode}
In \cite{b1}, we show that the framed isotopy type of $\mathbb{L} \subset -\Sigma(L)$ is completely determined by the linking matrix $A$ of $\mathbb{L'} \subset S^3$.  If follows that the pages $E^i$ for $i \geq 1$ are determined by $A$ as well (up to an overall shift in bigrading that depends on $n_\pm$).  In fact, since we are working with $\ftwo$ coefficients, the orientations of the arcs (and corresponding components of $\mathbb{L}'$) are extraneous as well.  We need only record which pairs of arcs in $\mathcal{D}(I^*)$ are linked, as well as $I^*$ itself.  On the other hand, the matrix $A$ with signs fully encodes the odd Khovanov homology of $L$ with $\zz$ coefficients (again, up to an overall shift in bigrading that depends on $n_\pm$), and should also encode a lift of the spectral sequence to $\zz$.
\end{remark}

\subsection{Grading}
\label{sec:khgrade}

For the duration of this paragraph, we return to the notation $\underline E^i$ to distinguish the $\hmhat$ version of the spectral sequence from the $\hmb$ version.  Our goal is to relate the mod 2 grading $\check \delta$ on $\underline E^1$ to the integer grading $\delta$ on $\ckh(L)$.  Since $U_\dagger$ is surjective on $\hm(\#^k S^1 \times S^2)$,  $\underline E^1$ may be identified as a $\check \delta$-graded vector space with the kernal of the map $$\sum_{I \in \{0,1\}^l} \hmb(U_\dagger \ | \ \yi \times [0,1]): E^1 \to E^1.$$  This permits us to work with the $\check \delta$ grading on $E^1$ instead.

\begin{proof}[Proof of Theorem \ref{thm:grading}]
Let $L \subset S^3$ be an oriented link and fix a diagram $\mathcal{D}$ with $n$ crossings.  Let $n_+$ and $n_-$ denote the number of positive and negative crossings, respectively.  Consider the hypercube complex $(X, \dd)$ given by surgeries on $\mathbb{L}' \subset -\Sigma(L)$.  Recall that $\yi = -\Sigma(\mathcal{D}(I)) \cong \#^k (S^1 \times S^2)$ when the resolution $\mathcal{D}(I)$ consists of $k+1$ circles.

We may think of a generator $x \in \Lambda^r(H_1(\yi))$ as an element of either $\underline E^1$ or $\ckh(\mathcal{D})$.  The group $\hm(\yi)$ is supported over the torsion spin$^c$ structure, and a short calculation shows that $\grq(x) = -r$, where $\grq(x)$ is the rational grading over the torsion spin$^c$-structure, defined in Section 28.3 of \cite{km}.  Moreover, on $\yi$ we have 
\begin{align}
\label{eqn:mod2}
\grq(x) \equiv \grt(x)  \mod 2.
\end{align}

Recall that $\ckh(\mathcal{D})$ has a quantum grading $q$ and a homological grading $t$.  The $\delta$-grading is defined as the linear combination $\delta = \frac{1}{2}q - t$.  Translating from the definitions in \cite{orsz}, we may express these gradings as
\begin{align*}
q(x) &= b_1(\yi) + 2\grq(x) + w(I) + n_+ - 2n_-\\
t(x) &= w(I) - n_- \\
\delta(x) &= \grq(x) - \frac{1}{2}w(I)+\frac{1}{2}b_1(\yi)+\frac{1}{2}n_+
\end{align*}
Here $n_+$ and $n_-$ denote the number of positive and negative crossings in $\mathcal{D}$.  The final formula defines a function $\delta : X \to \mathbb{Q}$.

We next define a function $\check{\delta}^\mathbb{Q}: X \to \mathbb{Z}$ which lifts the mod 2 grading $\check \delta$ on $X$.  Note that all cobordisms $\wij$ over the hypercube satisfy $\sigma(\wij) = 0$.  This is easily seen for an elementary cobordism, and follows in general from signature additivity. Using \eqref{eqn:iota} and \eqref{eqn:form1}, we have
$$\check{\delta}(x) = \grt(x) - \frac{1}{2}w(I) + \frac{1}{2} b_1(\yi) - \frac{1}{2}(\sigma(W_{\textit{0}\infty}) + b_1(\Sigma(L))) \mod 2.$$
By \eqref{eqn:mod2}, we may then define $\check{\delta}^\mathbb{Q}$ by
\begin{align*}
\check{\delta}^\mathbb{Q}(x) = \grq(x) - \frac{1}{2}w(I) + \frac{1}{2} b_1(\yi) - \frac{1}{2}(\sigma(W_{\textit{0}\infty}) + b_1(\Sigma(L))).
\end{align*}
Finally, we compare $\delta(x)$ with $\check\delta^\mathbb{Q}$:
\begin{align*}
\delta(x) - \check\delta^\mathbb{Q}(x) &= \frac{1}{2}(\sigma(W_{\textit{0}\infty}) + n_+ + b_1(\Sigma(L))). \\
&= \frac{1}{2}(\sigma(L) + \nu(L))
\end{align*}
The last line follows from Lemma \ref{lem:sig} below.  Reducing mod 2, we have the first claim of Theorem \ref{thm:1}.

For the remaining claim about the determinant, note that $$\chi_{\check\delta}(E^2) = (-1)^{(\sigma(L) + \nu(L))/2}\chi_\delta(\kh(L)) = (-1)^{(\sigma(L) + \nu(L))/2}V_L(-1) = \text{det}(L),$$
where $V_L(q)$ denotes the Jones polynomial of $L$ (when $\nu(L)>0$, everything vanishes).  Alternatively, one can show that the number of spin$^c$ structures on $-\Sigma(L)$ is $\det(L)$ and that $$\chi_{\grt} (\hmhat(-\Sigma(L),\mathfrak{s}))$$ is one when $\nu(L)  = 0$, and vanishes otherwise.
\end{proof}

\begin{lemma}
\label{lem:sig}
The signature and nullity of $L$ are given by
\begin{align*}
\sigma(L) &= \sigma(W_{\Ii \infty}) + n_+ \\
\nu(L) &= b_1(\Sigma(L))
\end{align*}
\end{lemma}
\begin{proof}
The nullity $\nu(L)$ is sometimes defined as the nullity of any symmetric Seifert matrix $S$ for $L$, and sometimes as $b_1(\Sigma(L))$.  These definitions are equivalent since $S$ presents $H_1(\Sigma(L))$.

We will prove the formula for $\sigma(L)$ by relating $\sigma(W_{\Ii \infty})$ to the signature of a certain 4-manifold $X_L$ bounding $\Sigma(L)$.  Recall that the diagram $\mathcal{D}$ has a checkerboard coloring with infinite region in white.  The black area forms a spanning surface $F$ for $L$ with one disk for each black region, and one half-twisted band for each crossing.  View $L$ as in the boundary of $D^4$, and push $F$ into the interior.   We then define $X_L$ as the branched double cover of $D^4$ over $F$.  In \cite{gl}, Gordon and Litherland show that the intersection form of $X_L$ is the Goeritz form $G$ associated to $\mathcal{D}$, and that
\begin{align*}
-\sigma(L) = \sigma(G) - \mu(\mathcal{D}),
\end{align*}
where $\mu(\mathcal{D}) = c- d$ (see Figure \ref{fig:crossings}).  The minus sign in front of $\sigma(L)$ is due to the fact that the signature convention in \cite{gl} is the opposite of ours.  Using the relations
\begin{align*}
w(B) &= b + c \\
n_- &= b + d
\end{align*}
we can also express $\mu(\mathcal{D})$ as
\begin{align*}
\mu(\mathcal{D}) = w(B) - n_-.
\end{align*}
Therefore, 
\begin{align}
\label{eqn:sig}
\sigma(L) = -\sigma(X_L) + w(B) - n_-.
\end{align}

\begin{figure}[htp]
\centering
\includegraphics[width=140mm]{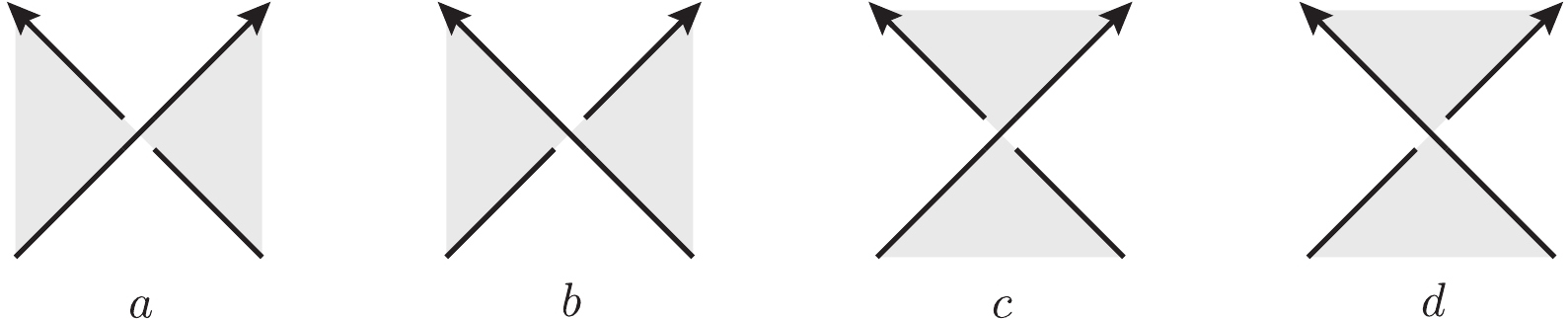}
\caption[Four types of crossings in an oriented link diagram with checkerboard coloring.]{Four types of crossings in an oriented diagram with checkerboard coloring.  The letters $a$, $b$, $c$, and $d$ denote the number of crossings of each type.}
\label{fig:crossings}
\end{figure}

We now construct a Kirby diagram for $X_L$ (see Section 3 of \cite{osz12} and Section 3 of \cite{gr} for similar constructions).  First, form the black graph resolution $\mathcal{D}(B)$ by resolving each crossing so as to separate the black regions into islands (that is, 1-resolve a crossing if and only if it is of type $b$ or $c$ in Figure \ref{fig:crossings}).  Draw a 1-handle in dotted circle notation along the boundary of each black region in $\mathcal{D}(B)$.  Next, add a 2-handle clasp at each crossing, with framing $+1$ if the crossing is of type $b$ or $c$, and $-1$ otherwise.  Finally, delete one of the 1-handles.

Since $\sigma(W_{0B})$ vanishes, signature additivity implies $\sigma(W_{0\infty}) = \sigma(W_{B\infty})$.
Next we construct a relative Kirby diagram for the cobordism $W_{B\infty}$.  First turn all but one of the circles in $\mathcal{D}(B)$ into 1-handles to get a surgery diagram for $Y_B = \Sigma(\mathcal{D}(B))$, regarded as the incoming end of $W_{B\infty}$.  Next, introduce a 0-framed clasp at each of the  $n - w(B)$ crossings corresponding to 0 digits of $B$.  This gives a relative Kirby diagram for the cobordism $W_{B\If}$.  Finally, introduce -1 framed clasps at each of the remaining crossings, and -1 framed meridians on each of the 0-framed clasps.  This gives a relative Kirby diagram for the cobordism $W_{B\infty}$.   After pulling off and blowing down all $n - w(B)$ of the -1 framed meridians, and filling in the incoming end with a boundary connect sum of copies of $S^1 \times D^3$, we recover the Kirby diagram for $-X_L$.  Therefore,
\begin{align*}
\sigma(W_{B\infty}) = -\sigma(X_L) - (n - w(B)).
\end{align*}
Combined with \eqref{eqn:sig}, we conclude
\begin{align*}
\sigma(L) &= -\sigma(X_L) + w(B) - n_- \\
&= \sigma(W_{B\infty}) + (n - w(B)) + w(B) - n_- \\
&= \sigma(W_{0\infty}) + n_+.
\end{align*}
\end{proof}

Modifying the above proof, we obtain the signature formula described in the introduction:

\begin{proof}[Proof of Proposition \ref{prop:newsig}]
Recall the construction of the surgery diagram $\mathbb{L}'$ for $-\Sigma(L)$, as in Figure \ref{fig:trefoil}.  Let $Z_L$ be the 4-manifold obtained by attaching 2-handles along $\mathbb{L}'$.  By construction, $Z_L$ bounds $-\Sigma(L)$.  Just as in the above proof, a Kirby diagram argument shows that
$$\sigma(W_{0\infty}) = \sigma(W_{I^*\infty}) = \sigma(Z_L) - (n - w(I^*)) = \sigma(A) - (n - w(I^*)),$$ where $A$ is the linking matrix of $\mathbb{L}$, which is congruent to the linking matrix of the arcs in Proposition \ref{prop:newsig}.  From Lemma \ref{lem:sig}, we arrive at the formula
\begin{align*}
\sigma(L) &= \sigma(W_{0\infty}) + n_+ \\
&= \sigma(A) + w(I^*) - n_-.
\end{align*}
Furthermore, $\det(L) = \left|\det(A)\right|$ since $A$ presents $H_1(-\Sigma(L))$.
\end{proof}

\begin{example}
Consider the resolution $\mathcal{D}(010)$ of the right-handed trefoil $\mathcal{T}$ in Figure \ref{fig:trefoil}, with $A$ given by the linking matrix at right.  The signature formula gives
\begin{align*}
\sigma(\mathcal{T}) = \sigma(A) + w(010) - n_-(\mathcal{D}) = 1 + 1 - 0 = 2.
\end{align*}
For the mirrored diagram $\overline{\mathcal{D}}$ representing the left-handed trefoil $\overline{\mathcal{T}}$, consider the mirrored resolution $\mathcal{D}(101)$.  Now the signature formula gives
\begin{align*}
\sigma(\overline{\mathcal{T}}) = \sigma(-A) + w(101) - n_-(\overline{\mathcal{D}})= -1 + 2 - 3 = -2.
\end{align*}
\end{example}

\subsection{Invariance}
\label{sec:linkinv}
We now turn to the proof of Theorem \ref{thm:linkinv}, which describes the extent to which the spectral sequence depends on the choice of diagram for the link L.

\begin{proof}[Proof of Theorem \ref{thm:linkinv}]
Let $\mathcal{D}_1$ and $\mathcal{D}_2$ be two diagrams of the link $L$.  Let $\underline X(\mathcal{D}_i)$ represent the hypercube complex associated to diagram $\mathcal{D}_i$, for some choice of analytic data (which we may suppress by Theorem \ref{thm:inv}).  The goal is to construct a filtered chain map
$$\phi: \underline X(\mathcal{D}_1) \to \underline X(\mathcal{D}_2)$$
which induces an isomorphism on the $\underline E^2$ page, and therefore on all higher pages as well.  It suffices to consider the case where $\mathcal{D}_1$ and $\mathcal{D}_2$ differ by a single Reidemeister move.

In \cite{b}, Baldwin defines such a map $\phi$ for each of the three Reidemeister moves.  While he was considering the Heegaard Floer version of the spectral sequence, his maps have direct analogues in the monopole Floer case.  The difficult part is proving that $\phi$ induces a homotopy equivalence from $\ckh(\mathcal{D}_1)$ to $\ckh(\mathcal{D}_2)$ on the $\underline E^1$ page.  However, this argument only involves properties of the Khovanov differential, drawing heavily on the proof that Khovanov homology is a bigraded link invariant (see \cite{orsz}).  It is also clear from the construction that $\phi$ preserves the bigrading on Khovanov homology, and therefore $\check \delta$.

Now, suppose links $L_1$ and $L_2$ are related by a mutation.  Fix diagrams $\mathcal{D}_1$ and $\mathcal{D}_2$ for $L_1$ and $L_2$ which exhibit the mutation.  Let $\mathbb{L}'_1$ and $\mathbb{L}'_2$ be the associated framed links in $S^3$, and let $\mathbb{L}_1$ and $\mathbb{L}_2$ be the associated framed links in $-\Sigma(L_1)$ and $-\Sigma(L_2)$.  In Section 4 of \cite{b1}, we prove that there is an orientation-preserving diffeomorphism
$$\psi':S^3 \to S^3$$
for which $\psi(\mathbb{L}'_1)$ and $\mathbb{L}'_2$ are isotopic as framed links.  This implies that there is an orientation-preserving diffeomorphism $$\psi:-\Sigma(L_1) \to -\Sigma(L_2)$$ for which $\psi(\mathbb{L}_1)$ and $\mathbb{L}_2$ are isotopic as framed links.  Appealing to Theorem \ref{thm:inv}, we conclude that the $\underline E^i$ pages of the spectral sequences associated to $\mathcal{D}_1$ and $\mathcal{D}_2$ agree for all $i \geq 1$.
\end{proof}


\subsection{The spectral sequence for a family of torus knots?}
\label{sec:torus}
\begin{figure}[htp]
\centering
\includegraphics[width=135mm]{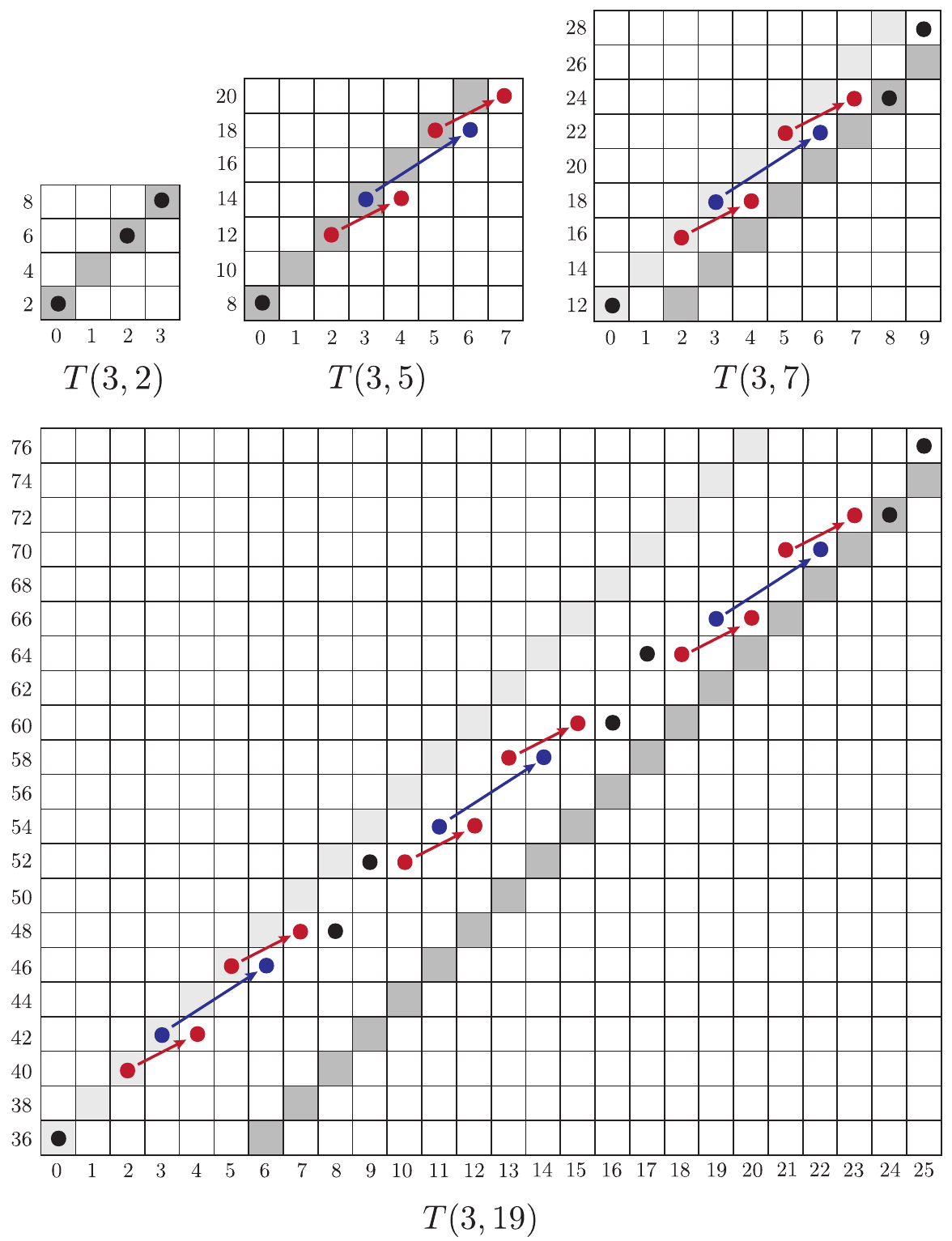}
\caption[Conjectural spectral sequence for torus knots of the form $T(3,6n \pm 1)$.]{Each dot represents an $\ftwo$ summand of $\kh(T(3,6n \pm 1))$ in the $(t,q)$-plane.  The diagonal $\delta = \sigma/2$ is heavily shaded and the diagonal $\delta = s/2$ is lightly shaded (unless $s = \sigma$).  The $d^2$ and $d^3$ differentials are in red and blue, respectively, as are their victims.  The surviving (black) dots generate $\hmhat(-\Sigma(2,3,6n \pm 1))$.  The shaded diagonals also correspond to $\check\delta^\mathbb{Q} = 0$ and $\check\delta^\mathbb{Q} = (2n-1)\pm 1$.  For $n \geq 1$, there is precisely one black dot on each diagonal in this range, giving the expected rank of $E^\infty$ in each $\grt$ grading.}
\label{fig:torus}
\end{figure}

In order to illustrate the spectral sequence in action, we now present an example which is both highly speculative and, we hope, compelling.  Consider the family of torus knots given by
$$\{T(3, 6n \pm 1) \, | \, n \geq 1\}.$$
For this family, the unreduced Khovanov homology with coefficients in $\mathbb{Q}$ takes the form of repeating blocks, and is stable, up to a shift in quantum grading, as $n$ grows \cite{stos}, \cite{turner}.  Computing $\kh(T(3,6n \pm 1))$ explicitly for several values of $n$ using Bar-Natan's {\em KnotTheory} package \cite{kat}, a similar pattern of repeating blocks emerges, as shown in Figure \ref{fig:torus}.  Surely this pattern persists and can be deduced by similar techniques.  Up to a shift in quantum grading, we have inclusions
$$\kh(T(3,5)) \subset \kh(T(3,7)) \subset \kh(T(3,11)) \subset \kh(T(3,13)) \subset \cdots.$$

\begin{conjecture}
\label{conj:torus}
For each such torus knot, and some choice of analytic data and diagram, the higher differentials are as shown in Figure \ref{fig:torus}.  In particular, the spectral sequence converges at the $E^4$ page, and the above inclusions on the $E^2$ page extend to the $E^3$ and $E^4$ pages.
\end{conjecture}

One intriguing way to frame this conjecture is as follows.  Using the $(t,q)$-bigrading in Figure \ref{fig:torus}, we may define {\em higher $\delta$-polynomials} by
\begin{align*}
U^k_{T(3,6n \pm 1)}(\delta) = \sum_{i,j} (-1)^i \, \text{rk} \, E^k_{i,j}(T(3,6n \pm 1)) \, \delta^{j/2-i}
\end{align*}
for each $k \geq 2$.  Then Conjecture \ref{conj:torus} implies that the $\delta$-polynomials on the $E^2$ and $E^3$ pages are monic monomials, while
\begin{align}
\label{eqn:u4}
U^4_{T(3, 6n \pm 1)} (\delta) = \delta^{\sigma/2} - \delta^{\sigma/2 +1}  + \delta^{\sigma/2 + 2} - \cdots + \delta^{s/2}.
\end{align}
Here $s$ denotes Rasmussen's $s$-invariant.  Based on Theorem \ref{thm:grading}, we suspect that these polynomials are highly relevant to the connection between Khovanov homology and Floer homology.  It would be very interesting to compare \eqref{eqn:u4} with the polynomials arising from Greene's conjectured $\delta$-grading on $\hf(-\Sigma(T(3, 6n \pm 1))$, defined in Section 8 of \cite{gr}.

Our primary evidence for Conjecture \ref{conj:torus} comes from \cite{b}, where Baldwin deduces that the Heegaard Floer spectral sequence for $T(3,5)$ is as shown in Figure \ref{fig:torus}.  His argument uses the Khovanov and Heegaard Floer contact invariants to show that the lower left generator survives to $E^\infty$ for every torus knot.  This is the only survivor in the case of $T(3,5)$, as the branched double cover is the Poincar\'e homology sphere.  Alas, we have not rigorously computed the monopole Floer spectral sequence even in this case, since for now we lack an analogous contact invariant.

As further evidence, we cite the compatibility of $\hmhat(-\Sigma(T(3,6n - 1)))$ with the $E^\infty$ page implied by our conjecture.  The branched double cover of $T(3, 6n \pm 1)$ is the Brieskorn integer homology sphere $$\Sigma(2,3,6n \pm 1),$$ which arises by $1/(6n \pm 1)$ Dehn surgery on a trefoil knot.  Using a surgery exact triangle, the Heegaard Floer groups $\mathit{HF}^+(\Sigma(2,3,6n \pm 1))$ are explicitly calculated in \cite{osz6}.  The same techniques should apply in the monopole case, and using \eqref{eqn:hatlong}, we expect that
\begin{align*}
\hmhat(-\Sigma(2,3,6n - 1)) &= \zz_{(-2)} \oplus \left(\zz_{(-2)} \oplus \zz_{(-1)} \right)^{n-1} \\
\hmhat(-\Sigma(2,3,6n + 1)) &= \zz_{(0)} \oplus \left(\zz_{(0)} \oplus \zz_{(1)} \right)^n
\end{align*}
where the subscript denotes $\grq$ grading (see also \cite{moyu}).  In particular, $$\text{rk} \, \hmhat(-\Sigma(2,3,6n \pm 1)) = 2n \pm 1.$$
We also expect that, for some choice of metric and perturbation, the first summand arises from the reducible generator of $\hmb(-\Sigma(T(3,6n - 1)))$ in lowest $\grq$ grading, while each pair of summands in parenthesis arises from a single irreducible generator.  Comparing with Figure \ref{fig:torus}, we imagine that the reducible is sitting in the lower left corner, whereas the pairs arising from each irreducible are in adjacent homological grading and $\delta$ grading.  In particular, the $\check \delta$ grading on $E^4$ is compatible with the $\grt$ grading on $\hmhat(-\Sigma(2,3,6n \pm 1))$, as required by Theorem \ref{thm:grading}.

For links such that the $(t,q)$ bigrading on the higher pages is well-defined, we may encode the higher pages of the spectral sequence in the form of a 2-variable {\em higher Khovanov polynomial}, given by
\begin{align*}
E^k_L(t,q) = \sum_{i,j} \text{rk} \, E^k_{i,j}(L) \, t^i q^j
\end{align*}
for each $k \geq 2$.  We then obtain {\em higher Jones polynomials}, given by
\begin{align*}
V^k_L(q) = E^k_L(-1,q^{1/2})
\end{align*}
for each $k \geq 2$.  The ordinary Jones polynomial $V_L(q)$ coincides with $V^2_L(q)$.  Furthermore, if $L$ is quasi-alternating, then $V^k_L(q) = V_L(q)$ for all $k \geq 2$.

We now record in full the various higher polynomials associated with the differentials in Figure \ref{fig:torus}, in order to provide a (conjectural) data-set to aid in the search for a combinatorial description.  For comparison, we include the polynomials on the $E^2$ page as well.  Note that the optimal input for an algorithm may not be a diagram of the link itself, but rather the arc-linking data which encodes the mutation equivalence class, as described in Remark \ref{rem:encode}.
\newpage

\begin{conjecture}
Let $\mathcal{S}^n = T(3,6n+1)$ and let $\mathcal{T}^n = T(3,6n-1)$.  Set
$$f_n(t,q) = \sum_{k = 0}^{n-1} t^8q^{12} \qquad \qquad f_n(q) = \sum_{k = 0}^{n-1} q^{6}$$
For each $n \geq 1$, the higher Khovanov polynomials are given by \\
\begin{align*}
q^{-s} E^2_{\mathcal{S}^n}(t,q) &= 1 +  \left((t^8 q^{12}+ t^5q^{16}) + (t^3 q^6 + t^6 q^{10})+  (t^2 q^4 + t^4 q^6 + t^5 q^{10} + t^7 q^{12})\right)f_n(t,q) \\
q^{-s} E^3_{\mathcal{S}^n}(t,q) &= 1 +  \left((t^8 q^{12}+ t^5q^{16}) + (t^3 q^6 + t^6 q^{10})\right)f_n(t,q) \\
q^{-s} E^4_{\mathcal{S}^n}(t,q) &= 1 +  (t^8 q^{12}+ t^5q^{16}) f_n(t,q) \\
q^{-s} E^2_{\mathcal{T}^n}(t,q) &= 1 +  (t^8 q^{12}+ t^5q^{16})f_{n-1}(t,q) + \left((t^3 q^6 + t^6 q^{10})+  (t^2 q^4 + t^4 q^6 + t^5 q^{10} + 
t^7 q^{12})\right)f_n(t,q) \\
q^{-s} E^3_{\mathcal{T}^n}(t,q) &= 1 +  (t^8 q^{12}+ t^5q^{16})f_{n-1}(t,q) + (t^3 q^6 + t^6 q^{10})f_n(t,q) \\
q^{-s} E^4_{\mathcal{T}^n}(t,q) &= 1 +  (t^8 q^{12}+ t^5q^{16}) f_{n-1}(t,q)
\end{align*}
\\
\noindent The higher Jones polynomials are given by
\\
\begin{align*}
q^{-s/2} \, V^2_{\mathcal{S}^n}(q) &= 1 +  q^2 - q^{6n+2} \hspace{109mm} \\
q^{-s/2} \, V^3_{\mathcal{S}^n}(q) &= 1 +  \left((q^6 - q^{8}) + (-q^3 + q^5)\right) f_n(q)\\
q^{-s/2} \, V^4_{\mathcal{S}^n}(q) &= 1 +  (q^{6} - q^{8}) f_n(q) \\
q^{-s/2} \, V^2_{\mathcal{T}^n}(q) &= 1 +  q^2 - q^{6n} \\
q^{-s/2} \, V^3_{\mathcal{T}^n}(q) &= 1 +  (q^{6}-q^{8})f_{n-1}(q) + (-q^3 + q^{5})f_n(q) \\
q^{-s/2} \, V^4_{\mathcal{T}^n}(q) &= 1 +  (q^{6}-q^{8}) f_{n-1}(q)
\end{align*}
\\
\noindent The higher $\delta$-polynomials are given by
\\
\begin{align*}
\delta^{-\sigma/2} \, U^2_{\mathcal{S}^n}(\delta) &= 1 \hspace{131mm} \\
\delta^{-\sigma/2} \, U^3_{\mathcal{S}^n}(\delta) &= 1\\
\delta^{-\sigma/2} \, U^4_{\mathcal{S}^n}(\delta) &= 1 - \delta + \delta^2 - \cdots + \delta^{2n}\\
\delta^{-\sigma/2} \, U^2_{\mathcal{T}^n}(\delta) &= \delta^{-1} \\
\delta^{-\sigma/2} \, U^3_{\mathcal{T}^n}(\delta) &= \delta^{-1} \\
\delta^{-\sigma/2} \, U^4_{\mathcal{T}^n}(\delta) &= 1 - \delta + \delta^2 - \cdots + \delta^{2n-2}
\end{align*}
\\
Here $s(\mathcal{S}_n) = 12n$, $s(\mathcal{T}_n) = 12n-4$, and $\sigma(\mathcal{S}_n) = \sigma(\mathcal{T}_n) = 8n$.  So both $U^4_{\mathcal{S}^n}(\delta)$ and $U^4_{\mathcal{T}^n}(\delta)$ may be expressed as
\begin{align*}
\delta^{\sigma/2} - \delta^{\sigma/2 +1}  + \delta^{\sigma/2 + 2} - \cdots + \delta^{s/2}.
\end{align*}
\end{conjecture}


\section{From pseudoholomorphic polygons to associahedra of metrics}
\label{sec:hfproof}

This section is intended for those readers interested in the relationship between monopole Floer homology and Heegaard Floer homology (or more broadly, gauge theory and Lagrangian intersection theory).  In the former, the identities arising from permutohedra encode the fact that the Floer chain maps associated to 2-handle additions {\em commute} up to homotopy.  In the latter, we have $A_\infty$ relations which encode {\em associativity} up to homotopy (say, for the Floer maps associated to a stack of 2-handles).  The $A_\infty$ relations arise from associahedra, which naturally compactify the space of conformal structures on polygons (see, for example, Section 2 of \cite{b}).  We now outline an alternative construction of the hypercube complex $X$ that is modeled closely on the version in Heegaard Floer homology \cite{osz12}.  In particular, the components of the differential will be built using only associahedra of metics.

Let $L \subset Y$ be an $l$-component, framed link.  By a standard construction starting with a Heegaard multi-diagram subordinate to a bouquet of the link $L$,  we fix a Heegaard diagram $(\Sigma, \boldsymbol\alpha, \boldsymbol\eta_I)$ for each $Y_I$ with $I \in \{0,1\}^l$ (see Section 4 in both \cite{osz12} and \cite{osz5}, keeping in mind Remark \ref{rem:bdcview}).  Fix vertices $I<J$ and let $\Gamma_{IJ}$ be the set of all paths from $I$ to $J$.  Fix a path $\gamma \in \Gamma_{IJ}$ given by $I = I_0 < I_1 < \cdots < I_n = J$.  We now give an alternative construction of the cobordism $\wij: \yi \to \yj$ with respect to the path $\gamma$.  

\begin{figure}[htp]
\centering
\includegraphics[width=155mm]{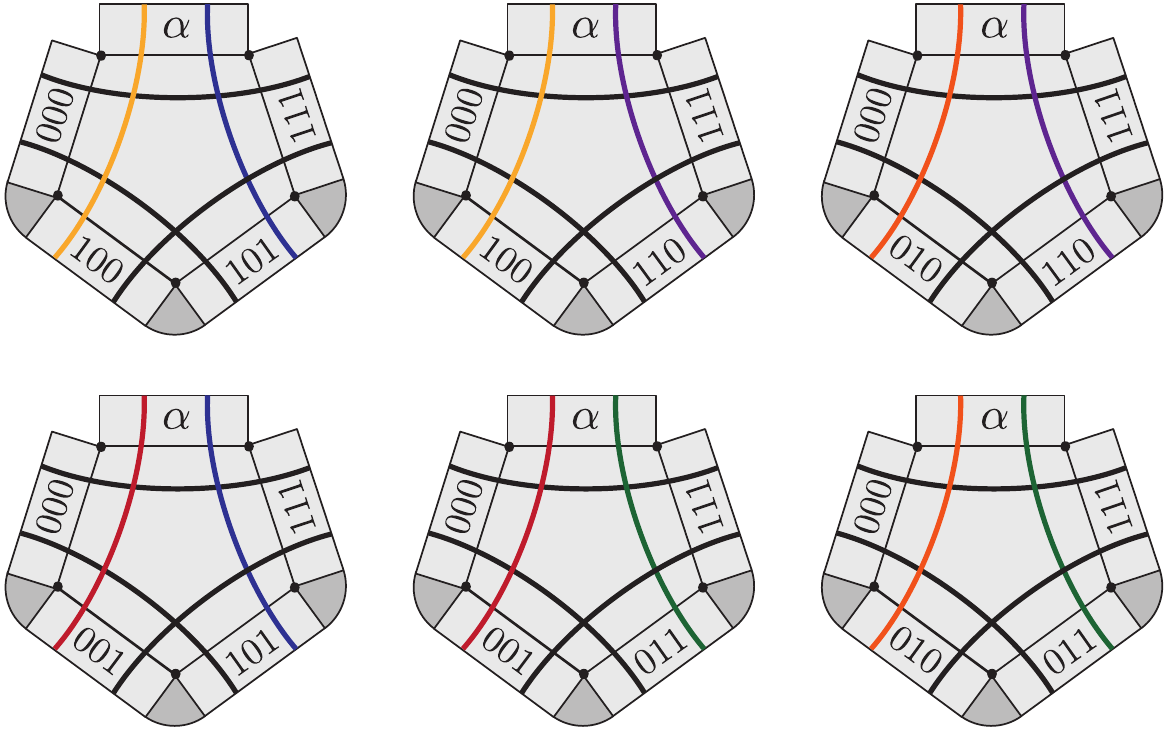}
\caption[Six pentagonal depictions of the cobordism for $\{0,1\}^3$.]{We have one depiction of $W$ for each of the six paths from $\Ii$ to $\If$ in $\{0,1\}^3$.  The regions that result from filling in the auxiliary boundary components are heavily shaded.  Conditions (i) and (ii) say that the (hypersurfaces associated to) arcs in distinct diagrams coincide in $W$ if and only if they have the same ends.  The internal hypersurfaces $\yi$ are colored consistently with Figures \ref{fig:intersections1} and \ref{fig:hexagon}.  The hexagon from Figure \ref{fig:hexagon} is built using the six squares of metrics that stretch on pairs of corresponding colored arcs in Figure \ref{fig:sixpaths}.  The remaining 24 squares cancel in pairs, in the sense of Proposition \ref{prop:familiar} below.}
\label{fig:sixpaths}
\end{figure}

Let $U_\alpha, U_0, \dots, U_n$ denote the $n+2$ handlebodies associated with $\gamma$ and let $\triangle_{n+2}$ denote the regular polygon with $n + 2$ edges $e_\alpha, e_0, \dots, e_n$.  We thicken these handlebodies, glue them onto $\triangle_{n+2} \times \Sigma$, and smooth the corners over the vertices of $\triangle_{n+2}$ to obtain a smooth, oriented cobordism $W_\gamma$ with $n+2$ boundary components (compare with Section 2.2 of \cite{osz5}):
\begin{align*}
W_\gamma = \frac{\triangle_{n+2} \times \Sigma \ \coprod \ e_\alpha \times U_\alpha \ \coprod \ e_0 \times U_0 \ \coprod \ \cdots \ \coprod \ e_n \times U_n}{e_\alpha \times \Sigma \sim e_\alpha \times \partial U_\alpha \ , \ e_0 \times \Sigma \sim e_0 \times \partial U_0 \ , \ \dots \ , \ e_n \times \Sigma \sim e_n \times \partial U_n}
\end{align*}
We have one embedded hypersurfaces in $W_\gamma$ for each pair of handlebodies.  The Heegaard diagram $(\Sigma, \boldsymbol\alpha, \boldsymbol\eta_{I_i})$ represents the internal hypersurface $Y_{I_i}$, which we also denote by $Y_i(\gamma)$.  The Heegaard diagram $(\Sigma, \boldsymbol\eta_{I_j}, \boldsymbol\eta_{I_k})$ represents an auxiliary hypersurface $S^j_k(\gamma)$ which is diffeomorphic to a connect sum of the form $\#^r (S^1 \times S^2)$.  Furthermore, $$\partial W_\gamma = Y_I \  \textstyle\coprod \ S^0_1(\gamma) \ \coprod \ S^1_2(\gamma) \ \coprod \ \cdots \ \coprod \ S^{n-1}_n(\gamma) \ \coprod \ Y_J.$$ 
After filling in each auxiliary boundary component of $W_\gamma$ with a connect sum of the form $\#^r (S^1 \times D^3)$, we recover precisely the cobordism $\wij$.  In Figure \ref{fig:sixpaths}, we have depicted this construction for the six paths in the hypercube $\{0,1\}^3$.

We embed the full set of internal hypersurfaces $Y_K$ in $\wij$ as before, as well as an auxiliary hypersurface $S^{K_1}_{K_2}$ for all $I \leq K_1 < K_2 \leq J$ with $w(K_2) - w(K_1) > 1$, so that:
\renewcommand{\labelenumi}{(\roman{enumi})}
\begin{enumerate}
\item $Y_i(\gamma)$ is identified with $Y_K$ whenever $\gamma = \{{I_0} < \cdots < I_n\}$ has $I_i = K$.
\item $S^{k_1}_{k_2}(\gamma)$ is identified with $S^{K_1}_{K_2}$ whenever the subpath $\gamma^{k_1}_{k_2} = \{ I_{k_1} < \cdots < I_{k_2} \}$ runs from $K_1$ to $K_2$.
\item The collection of hypersurfaces $$\{Y_i(\gamma)\} \bigcup \{S^j_k(\gamma) \ | \ k - j> 1\}$$ associated to any particular $\gamma$ intersect as suggested by Figure \ref{fig:sixpaths}.
\end{enumerate}
When $I = \Ii$ and $J = \If$, we have constructed the full cobordism $W$.  Fix an initial metric and perturbation on $W$, and on each $\wij$ by restriction.  The intersection data required by (iii) corresponds to the lattice $\{0,1,\dots, n\}$ in Theorem \ref{thm:graph}.  Degenerating the metric along a hypersurface mirrors degenerating the conformal structure on the polygon $\triangle_{n+2}$ along the corresponding arc.  So for each $\gamma \in \Gamma_{IJ}$, we have an associahedron of metrics $K_\gamma$ on $\wij$.  We define $\check{d}_\gamma: \cm(\yi) \to \cm(\yj)$ to be the map associated to $K_\gamma$, and then 
$\check{d}^I_J : \cm(\yi) \to \cm(\yj)$ by 
$$\check{d}^I_J = \sum_{\gamma \in \Gamma_{IJ}} \check{d}_\gamma.$$  Finally, the map $\check d : \tilde X \to \tilde X$ is given by the sum of the $\check{d}^I_J$ over all pairs $I \leq J$.

We now show that $\check{d}^2 = 0$ by considering the $A_\infty$ relations arising from the associahedra $K_\gamma$ (this is analogous to the approach in \cite{osz12}).  Note there will be no terms corresponding to pinching off bigons, since we have effectively placed the generating cycle $\Theta$ on each $S^j_{j+1}(\gamma)$ by filling it in.  We let $\check{S}^j_k(\gamma)$ be the operator associated to the facet $S^j_k(\gamma)$ in $K_\gamma$.  Using the appropriate end-counting map $\check{A}_\gamma$, we obtain the $A_\infty$ relation
\begin{align*}
\check{d}_\gamma \check{d}^I_I + \check{d}^J_J \check{d}_\gamma  = \sum_{k=1}^{n-1} \check{d}_{\gamma^k_n} \check{d}_{\gamma^0_k} + \sum_{0\leq j< j+1 < k \leq n} \check{S}^j_k(\gamma)
\end{align*}
which we rewrite as
\begin{align}
\label{eqn:pathid}
\sum_{k=0}^{n} \check{d}_{\gamma^k_n} \check{d}_{\gamma^0_k} = \sum_{0\leq j< j+1 < k \leq n} \check{S}^j_k(\gamma).
\end{align}
Summing \eqref{eqn:pathid} over all $\gamma \in \Gamma_{IJ}$, we have
\begin{align*}
\sum_{I \leq K \leq J} \check{d}^K_J \check{d}^I_K &= \sum_{\gamma \in \Gamma_{IJ}} \sum_{0\leq j< j+1<k \leq n} \check{S}^j_k(\gamma) \\
&= \sum_{\substack{I \leq K_1 < K_2 \leq J \\ w(K_2) - w(K_1) > 1}} \sum_{\gamma_1 \in \Gamma_{I K_1}} \sum_{\gamma_2 \in \Gamma_{K_2 J}} \sum_{\gamma \in \Gamma_{K_1 K_2}} \check{S}^{w(K_1) - w(I)}_{w(K_J) - w(K_2)} (\gamma_1 \cdot \gamma \cdot \gamma_2) \\
&= \sum_{\substack{I \leq K_1 < K_2 \leq J \\ w(K_2) - w(K_1) > 1}} \sum_{\gamma_1 \in \Gamma_{I K_1}} \sum_{\gamma_2 \in \Gamma_{K_2 J}} \left(\sum_{\gamma \in \Gamma_{K_1 K_2}} \check{S}^{K_1}_{K_2}\right) = 0
\end{align*}
since the innermost sum has an even number of terms, which are identical by (ii).  Here, $\gamma_1 \cdot \gamma \cdot \gamma_2$ denotes the concatenation of the three paths.  We conclude that $\check d^2 = 0$.

In fact, we have reproduced the maps associated to the permutohedra $\pij$:
\begin{proposition}
\label{prop:familiar}
Using the same choice of initial metric and perturbation on $W$, we have $$\check{d}^I_J = \dij.$$
\end{proposition}

\begin{proof}
Let $n = w(J) - w(I)$.  Each of the $n!$ associahedra $K_\gamma$ contains a single cube $C_\gamma$ which stretches on internal hypersurfaces $Y_{I_1} < \dots < Y_{I_{n-1}}$.  These cubes are exactly those we glued together to form $\pij$ in Section \ref{sec:surface}.  The remaining cubes in the $n!$ associahedra fall into equivalence classes of even size, where two cubes are equivalent if they parameterize identical families of metrics on $\wij$.  Thus, their total contribution to $\check{d}^I_J$ is zero.
\end{proof}

The map analogous to $\check{d}_\gamma$ in \cite{osz12} counts pseudo-holomorphic polygons $\triangle_{n+2}$ in the Heegaard multi-diagram $(\Sigma_g, \boldsymbol{\alpha}, \boldsymbol\eta(I_0), \dots, \boldsymbol\eta(I_n)),$ where the conformal structure on $\triangle_{n+2}$ is allowed to vary.  Figure \ref{fig:sixpaths} directly illustrates how the space of conformal structures on $\triangle_{n+2}$ corresponds to an associahedral family of metrics $K_\gamma$ on $W$.  Conversely, the maps $\dij$ associated to permutohedra of metrics do not have {\em direct} analogues in Heegaard Floer homology.  Of course, this did not prevent Ozsv\'ath and Szab\'o from deriving their version of the link surgery spectral sequence.  In order to prove that $\check d$ is a differential, they essentially used the above algebraic cancellation, in place of the geometric cancellation implicit in excising the duplicate cubes and gluing together what remains.  We hope that the reader familiar with \cite{osz12} has gleaned new intuition from access to a more expansive repository of polytopes.

\newpage
\section*{Appendix I: Morse homology via path algebras}

Monopole Floer homology may be viewed as an infinite dimensional version of Morse homology for manifolds with boundary.  For a beautiful treatment of the finite dimensional model, see Section 2 of \cite{km}.  We now give a brief presentation of its essential features, assuming familiarity with Morse homology for closed manifolds.   By recasting the combinatorics in terms of path algebras, we hope to illuminate the classification of ends in Lemma \ref{lem:ends} and the form of the matrices \eqref{eqn:dij}, \eqref{eqn:aij}, \eqref{eqn:exactd2}, \eqref{eqn:l}, and \eqref{eqn:exacta2} used to define $\dij$, $\check L$, and $\aij$.

Consider a manifold with boundary, equipped with a Morse function whose gradient is everywhere tangent to the boundary. The critical points in the boundary are classified as stable or unstable, according to whether the flow in the normal direction is toward or away from the boundary, respectively.  We denote interior, boundary-stable, and boundary-unstable critical points by $o$, $s$, and $u$, respectively.  Note that interior gradient trajectories always flow from $o$ or $u$ to $o$ or $s$, whereas boundary trajectories flow from $s$ or $u$ to $s$ or $u$.  We distinguish between interior and boundary trajectories from $u$ to $s$, so there are eight types in all.

\begin{figure}[htp]
\centering
\includegraphics[width=150mm]{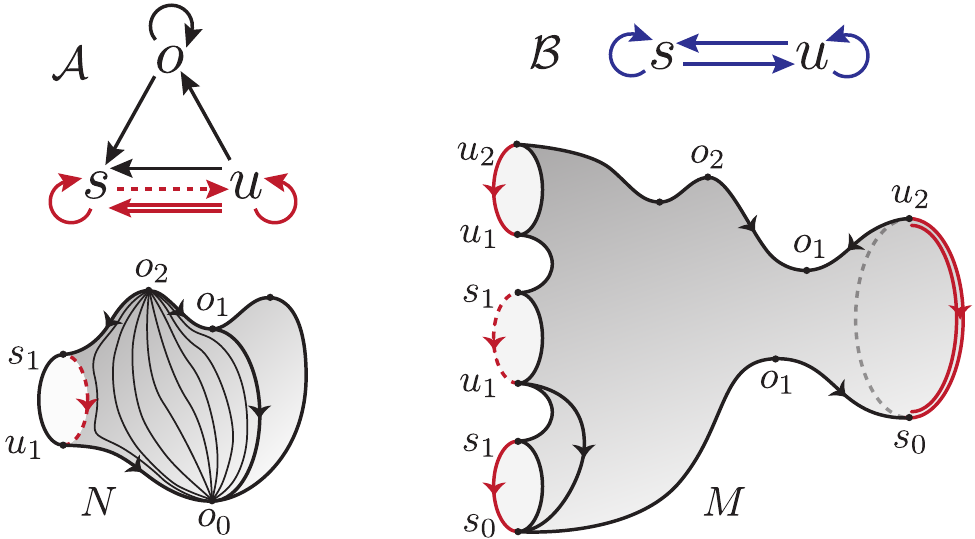}
\caption[Path algebras and Morse homology for manifolds with boundary.]{Path algebras and Morse homology for manifolds with boundary.}
\label{fig:morse}
\end{figure}

On the surface $M$ in Figure \ref{fig:morse}, we have marked one isolated gradient trajectory for each of these eight types, where those in $\partial M$ (in red) are isolated with respect to $\partial M$.  The subscripts on the critical points denote Morse index with respect to $M$.   While most of these isolated trajectories lower Morse index by 1, there are two exceptions.  The doubled trajectory from $u$ to $s$ in $\partial M$ lowers Morse index by 2, while the dashed trajectory from $s$ to $u$ in $\partial M$ fixes Morse index.  This last type is called {\em boundary-obstructed}.

All of this information may be neatly encoded in a weighted path algebra (or quiver) over $\ftwo$, denoted $\mathcal{A}$ (this was first pointed out to the author by Dylan Thurston).  As an $\ftwo$-vector space, $\mathcal{A}$ has a basis given by the set of all paths in the directed graph at left in Figure \ref{fig:morse}.  The product of two paths is given by concatenation if the end of the first coincides with the beginning of the second, and is zero otherwise.  The weight of a path is the sum of the weights of its edges, where the dashed, single, and doubled edges have weights 0, 1, and 2, respectively.

If we consider $\partial M$ as a closed manifold in its own right, then the Morse index of each boundary-unstable critical point is one less.  So now all four types of isolated trajectories in $\partial M$ lower Morse index by 1, as encoded in the path algebra $\mathcal{B}$ in Figure \ref{fig:morse}.

\begin{figure}[htp]
\centering
\includegraphics[width=120mm]{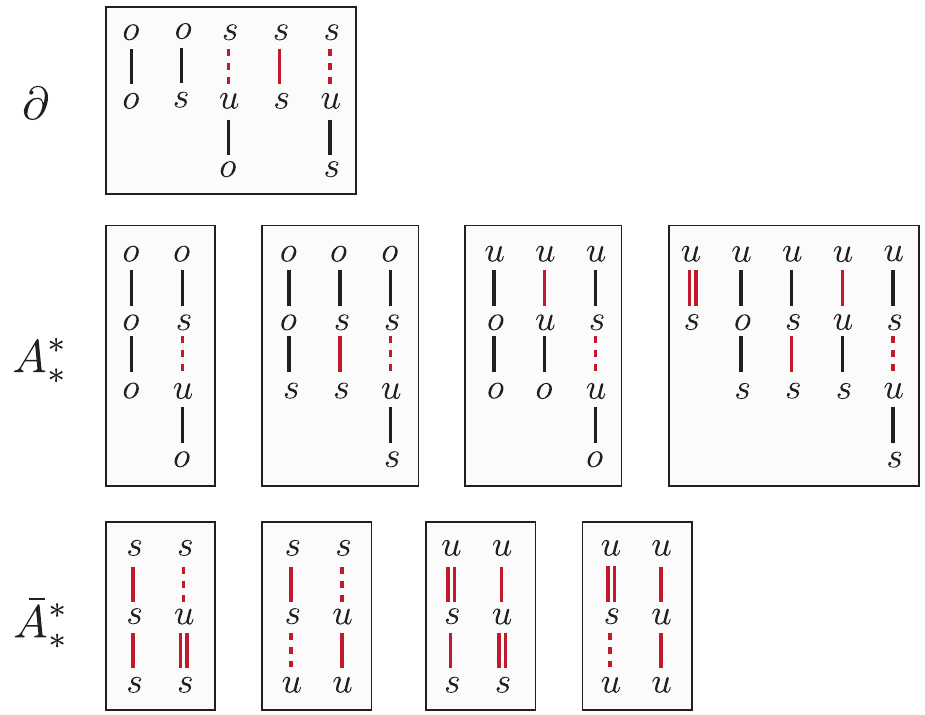}
\caption[The differential and identities as elements of a path algebra.]{The differential $\partial$ may be thought of as an element of the path algebra $\mathcal{A}$.  The eight elements in rows $A$ and $\bar A$ define an ideal of $\mathcal{A}$.}
\label{fig:morsepaths}
\end{figure}

The groups $H_*(\partial M)$, $H_*(M)$, and $H_*(M, \partial M)$ arise from the Morse complex generated by critical points of types $\{s, u\}$, $\{o,s\}$, and $\{o,u\}$, respectively.  The correspondence with the monopole Floer groups is reflected by the exact sequences
\begin{align*}
\cdots \longrightarrow H_*(\partial M) \longrightarrow H_*(M)  \longrightarrow H_*(M, \partial M) \longrightarrow \cdots  \\
\cdots \longrightarrow \hmbar(Y) \longrightarrow \hmb(Y) \longrightarrow \hmfrom(Y)  \longrightarrow \cdots.
\end{align*}
Since we are primarily concerned with $\hmb(Y)$, we focus on the absolute case $H_*(M)$.  The Morse complex then has the form
$$C(M) = C^o(M) \oplus C^s(M).$$
The differential $\partial$ may be thought of as an element of $\mathcal{A}$, given by the sum of all weight 1 paths from $\{o,s\}$ to $\{o,s\}$, as depicted in Figure \ref{fig:morsepaths}.  In matrix form, this becomes
\begin{align*}
\partial &= \left[\begin{array}{rr}
\doo & \duo\esu\\
\dos & \ess + \dus\esu
\end{array}\right].
\end{align*}

We introduce an ideal of $\mathcal{A}$, generated by the other eight elements in Figure \ref{fig:morsepaths}.  We have one relation for each interior (black) generator of $\mathcal{A}$, given by the sum all paths of weight 2 between its ends.  We similarly have one relation for each boundary (red) generator of $\mathcal{A}$, given by the sum all paths of weight 2 between the ends of the corresponding (blue) generator in $\mathcal{B}$.  These relations correspond precisely to maps counting the ends of 1-dimensional moduli spaces, and can be expressed in that form as
\begin{align*}
\aoo &=  \doo \doo + \duo\esu\dos & \abss &= \ess\ess + \eus\esu \\
\aos &= \dos\doo + \ess\dos + \dus\esu\dos & \absu &= \esu\ess + \euu\esu \\
\auo&= \doo\duo + \duo\euu + \duo\esu\dus & \abus &= \ess\eus + \eus\euu \\
\aus &=  \eus + \dos\duo + \ess\dus + \dus\euu + \dus\esu\dus & \abuu &= \esu\eus + \euu\euu
\end{align*}
We have illustrated two broken trajectories counted by the map $\aoo$ on the surface $N$ in Figure \ref{fig:morse}.  The 1-dimensional family of interior trajectories from $o$ to $o$ has one end with two components, and another end with three components, where the middle component is boundary-obstructed.  Note that the terms in the above relations correspond precisely to those described in Lemma \ref{lem:ends}.

As elements of $\mathcal{A}$, $\absu$ and $\abus$ have weights 1 and 3, respectively, while the other relations have weight 2.  Next, we define an element $A \in \mathcal{A}$ by extending each relation (if necessary and possible) to a weight 2 element from $\{o,s\}$ to $\{o,s\}$ and summing.  One observes cancellation of precisely those paths with no interior $o$ or $s$ (that is, no good break).  As a map, $A$ is given by
\begin{align*}
A &= \left[\begin{array}{rr}
\aoo & \auo\esu + \duo\absu \\
\aos & \abss + \aus\esu + \dus\absu
\end{array}\right] 
\end{align*}
It is now easy and instructive to check that $\partial^2$ and  $A$ coincide as elements of the (free) path algebra $\mathcal{A}$, so that $\partial^2$ is in the ideal generated by the relations.  Of course, this is the simplest case of the calculation done in Lemma $\ref{lem:aij}$, with the implication being that $\partial$ is a differential.

\begin{figure}[htp]
\centering
\includegraphics[width=150mm]{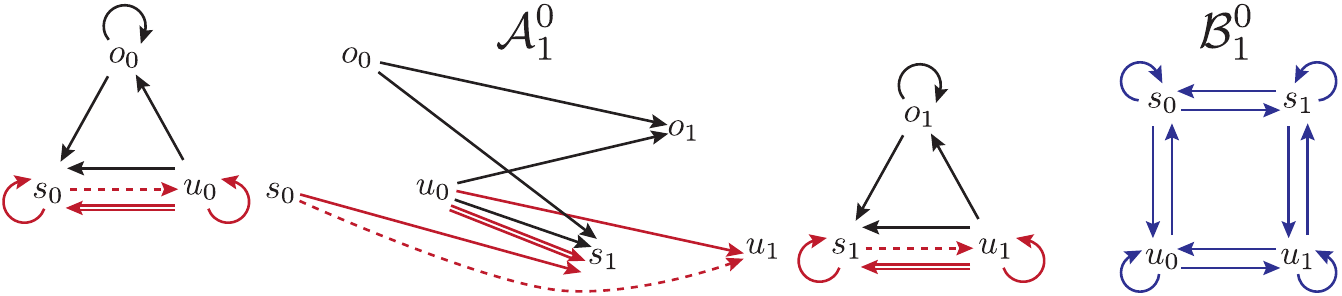}
\caption[The path algebra of a cobordism with fixed metric.]{The maps $\ddd^0_0$, $\mm^0_1$, and $\ddd^1_1$ may be thought of as elements of the weighted path algebra $\mathcal{A}^0_1$ over the red and black graph with 24 edges and 6 vertices.  We have one relation for each black edge in $\mathcal{A}^0_1$ and one relation for each blue edge in $\mathcal{B}^0_1$, each consisting of all paths of weight 2 between the corresponding ends.  We record the 24 relations above (the first 8 occur twice).}
\label{fig:pathalg}
\end{figure}

Kronheimer and Mrowka discuss functoriality in Morse homology in Section 2.8 of \cite{km}.  The above path algebra interpretation may be generalized to describe such maps, including those equipped with families of metrics and multiple incoming and outgoing ends, by including one copy of $o$, $s$, and $u$ for each end (see Figure \ref{fig:pathalg} for an example).  A map which counts unbroken trajectories on a cobordism with $n$ boundary components is represented by a star graph with $n$ labelled leaves.  The notion of boundary-obstructed trajectories generalizes in a natural manner to determine the weight of such a graph (see {\em boundary-obstructed of corank $c$} in Section 24.4 of \cite{km}).  In the case of 2-boundary components, everything is governed by the path algebras $\mathcal{A}$ and $\mathcal{B}$, as should be clear from the explicit examples of the next section.   Indeed, the author and Dave Bayer wrote a program in {\em Haskell}, which formally implements the path algebra associated to a cobordism equipped with a permutohedron of metrics (the hypercube case).  Sure enough, the program verifies Lemma \ref{lem:aij} in this language.

\newpage

\section*{Appendix II: Explicit maps and identities}

The path algebra formalism of Appendix I dictactes the form of the differential $\dd$ inducing the spectral sequence.  While the general case is given in \eqref{eqn:dij}, we now look at the maps associated to the lattices $\{0,1\}^l$ for $0 \leq l \leq 3$ in greater detail.   When $l = 1$, $2$, and $3$, the polytope of metrics $P_l$ will be a point, an interval, and a hexagon, respectively.  We record the component maps of
$$\ddd = \dd^0_0 \qquad \mm = \dd^0_1 \qquad \check{H} = \dd^{00}_{11} \qquad \check{G} = \dd^{000}_{111}$$
dropping superscripts and subscripts when it does not cause ambiguity:

\begin{align*}
\ddd &= \left[\begin{array}{rr}
\doo & \duo\esu\\
\dos & \ess + \dus\esu
\end{array}\right]\\
\mm &= \left[\begin{array}{rr}
\moo & \muo\esu + \duo\nsu\\
\mos & \nss + \mus\esu + \dus\nsu
\end{array}\right]\\
\hh &= \left[\begin{array}{rr}
\hoo & \huo\esu + \duo\ksu + \muob\nsua + \muod\nsuc\\
\hos & \kss + \hus\esu + \dus\ksu + \musb\nsua + \musd\nsuc
\end{array}\right]\\
\check{G} &= \left[\begin{array}{rr}
G^o_o & G^u_o\esu + \duo\bar{G}^s_u + H^u_o({}^{001}_{111})\bar{m}^s_u({}^{000}_{001}) + m^u_o({}^{011}_{111})\bar{H}^s_u({}^{000}_{011}) \\
G^o_s & \bar{G}^s_s + G^u_s\esu + \dus\bar{G}^s_u + H^u_s({}^{001}_{111})\bar{m}^s_u({}^{000}_{001}) + m^u_s({}^{011}_{111})\bar{H}^s_u({}^{000}_{011})
\end{array}\right] \\
&+ \left[\begin{array}{rr}
0 & H^u_o({}^{010}_{111})\bar{m}^s_u({}^{000}_{010}) + m^u_o({}^{101}_{111})\bar{H}^s_u({}^{000}_{101})+ H^u_o({}^{100}_{111})\bar{m}^s_u({}^{000}_{100}) + m^u_o({}^{110}_{111})\bar{H}^s_u({}^{000}_{110}) \\
0 & H^u_s({}^{010}_{111})\bar{m}^s_u({}^{000}_{010}) + m^u_s({}^{101}_{111})\bar{H}^s_u({}^{000}_{101})+ H^u_s({}^{100}_{111})\bar{m}^s_u({}^{000}_{100}) + m^u_s({}^{110}_{111})\bar{H}^s_u({}^{000}_{110})
\end{array}\right]
\end{align*}

Next, we record the component maps of
$$\Aa = \Aa^0_0 \qquad \check{B} = \Aa^0_1 \qquad \check{E} = \Aa^{00}_{11} \qquad \check{F} = \Aa^{000}_{111}$$
These count boundary points of 1-dimensional moduli spaces and thus all vanish identically.   We will omit those components of $\check F$ besides $F_o^o$ (shown below), as they run to three pages.

\begin{align*}
\Foo &=  \goo\doo + \doo\goo + \duo\esu\gos + \duo\gbsu\dos + \guo\esu\dos \\
&+ \hhood\mmooa + \hhuod\nnsua\dos + \hhuod\esu\mmosa + \duo\kksud\mmosa \\
&+ \hhooe\mmoob + \hhuoe\nnsub\dos + \hhuoe\esu\mmosb + \duo\kksue\mmosb \\
&+ \hhoof\mmooc + \hhuof\nnsub\dos + \hhuof\esu\mmosc + \duo\kksuf\mmosc \\
&+ \mmooj\hhooa + \mmuoj\kksua\dos + \mmuoj\esu\hhosa + \duo\nnsuj\hhosa \\
&+ \mmook\hhoob + \mmuok\kksub\dos + \mmuok\esu\hhosb + \duo\nnsuk\hhosb \\
&+ \mmool\hhooc + \mmuol\kksuc\dos + \mmuol\esu\hhosc + \duo\nnsul\hhosc \\
&+ \mmuoj\nnsud\mmosa + \mmuok\nnsue\mmosa + \mmuoj\nnsuf\mmosb \\
&+ \mmuol\nnsug\mmosb + \mmuok\nnsuh\mmosc + \mmuol\nnsui\mmosc \\
\end{align*}

\begin{align*}
\aoo &=  \doo \doo + \duo\esu\dos \\
\aos &= \dos\doo + \ess\dos + \dus\esu\dos \\
\auo&= \doo\duo + \duo\euu + \duo\esu\dus \\
\aus &=  \eus + \dos\duo + \ess\dus + \dus\euu + \dus\esu\dus \\
&\\
\abss &= \ess\ess + \eus\esu \\
\absu &= \esu\ess + \euu\esu \\
\abus &= \ess\eus + \eus\euu \\
\abuu &= \esu\eus + \euu\euu \\
&\\
\boo &= \moo\doo + \doo\moo + \duo\esu\mos + \duo\nsu\dos + \muo\esu\dos \\
\bos &= \mos\doo + \dos\moo + \mus\esu\dos + \dus\nsu\dos + \dus\esu\mos + \nss\dos + \ess\mos \\
\buo &= \moo\duo + \doo\muo + \muo\esu\dus  + \duo\nsu\dus + \duo\esu\mus + \muo\euu +\duo\nuu \\
\bus &= \mos\duo + \dos\muo + \mus\esu\dus  + \dus\nsu\dus + \dus\esu\mus + \mus\euu +\dus\nuu \\
&+ \nss\dus + \ess\mus + \nus \\
&\\
\bbss &= \nss\ess + \ess\nss + \nus\esu + \eus\nsu \\
\bbsu &= \nsu\ess + \esu\nss + \nuu\esu + \euu\nsu \\
\bbus &= \nss\eus + \ess\nus + \nus\euu + \eus\nuu \\
\bbuu &= \nsu\eus + \esu\nus + \nuu\euu + \euu\nuu \\
\end{align*}

\begin{align*}
\Eoo &= \hoo\doo + \doo\hoo + \duo\esu\hos + \duo\ksu\dos + \huo\esu\dos \\
&+ \muob\nsua\dos  + \muod\nsuc\dos  \\
&+ \moob\mooa +\muob\esu\mosa + \duo\nsub\mosa\\
&+ \mood\mooc + \muod\esu\mosc + \duo\nsud\mosc\\
& \\
\Eos &= \hos\doo + \dos\hoo + \hus\esu\dos + \dus\ksu\dos + \dus\esu\hos + \kss\dos + \ess\hos \\
&+ \musb\nsua\dos + \musd\nsuc\dos  \\
&+ \mosb\mooa + \musd\esu\mosc + \dus\nsud\mosc\\
&+ \mosd\mooc + \musd\esu\mosc + \dus\nsud\mosc\\
& \\
\Euo &= \hoo\duo + \doo\huo + \huo\esu\dus  + \duo\ksu\dus + \duo\esu\hus + \huo\euu +\duo\kuu \\
&+ \muob\nsua\dus + \muod\nsuc\dus \\
&+ \moob\muoa + \muob\esu\musa + \duo\nsub\musa + \muob\nuua\\
&+ \mood\muoc + \muod\esu\musc + \duo\nsud\musc + \muod\nuuc\\
& \\
\Eus &= \hos\duo + \dos\huo + \hus\esu\dus  + \dus\ksu\dus + \dus\esu\hus + \hus\euu +\dus\kuu + \kss\dus + \ess\hus + \kus \\
&+ \musb\nsua\dus + \musd\nsuc\dus\\
&+ \mosb\muoa + \musb\esu\musa + \dus\nsub\musa + \nssb\musa + \musb\nuua\\
&+ \mosd\muoc + \musd\esu\musc + \dus\nsud\musc + \nssd\musc + \musd\nuuc\\
& \\
& \\
\Ebss &= \kss\ess + \ess\kss + \kus\esu + \eus\ksu \\
&+ \nssb\nssa + \nusb\nsua + \nssd\nssc + \nusd\nsuc \\
& \\
\Ebsu &= \ksu\ess + \esu\kss + \kuu\esu + \euu\ksu \\
&+ \nsub\nssa + \nuub\nsua + \nsud\nssc + \nuud\nsuc \\
& \\
\Ebus &= \kss\eus + \ess\kus + \kus\euu + \eus\kuu \\
&+  \nssb\nusa + \nusb\nuua + \nssd\nusc + \nusd\nuuc\\
& \\
\Ebuu &= \ksu\eus + \esu\kus + \kuu\euu + \euu \kuu \\
&+ \nsub\nusa + \nuub\nuua + \nsud\nusc + \nuud\nuuc \\
\end{align*}

We bundle these identities into the vanishing maps \\
\begin{align*}
\check{A} &= \left[\begin{array}{rr}
\aoo & \auo\esu + \duo\absu \\
\aos & \abss + \aus\esu + \dus\absu
\end{array}\right] \\ \\
\check{B} &= \left[\begin{array}{rr}
\boo & \buo\esu + \duo\bbsu + \muo\absu + \auo\nsu\\
\bos & \bbss + \bus\esu + \dus\bbsu + \mus\absu + \aus\nsu
\end{array}\right] \\ \\
\check{E} &= \left[\begin{array}{ll}
\Eoo & \Euo\esu + \duo\Ebsu + \huo\absu + \auo\ksu \\
\Eos & \Ebss + \Eus\esu + \dus\Ebsu + \hus\absu + \aus\ksu
\end{array}\right] \\
&+ \left[\begin{array}{ll}
0 & \buob\nsua + \muob\bbsua + \buod\nsuc + \muod\bbsuc \\
0 & \busb\nsua + \musb\bbsua + \busd\nsuc + \musd\bbsuc
\end{array}\right] \\ \\
\check{F} &= \left[\begin{array}{rr}
\Foo & \Fuo\esu + \duo\Fbsu + \guo\absu + \auo\gbsu \\
\Fos & \Fbss + \Fus\esu + \dus\Fbsu + \gus\absu + \aus\gbsu
\end{array}\right] \\
&+ \left[\begin{array}{ll}
0 & \Euo({}^{001}_{111})\nsu({}^{000}_{001}) + \huo({}^{001}_{111})\bbsu({}^{000}_{001}) + \muo({}^{110}_{111})\Ebsu({}^{000}_{110}) + \buo({}^{110}_{111})\ksu({}^{000}_{110}) \\
0 & \Eus({}^{001}_{111})\nsu({}^{000}_{001}) + \hus({}^{001}_{111})\bbsu({}^{000}_{001}) + \mus({}^{110}_{111})\Ebsu({}^{000}_{110}) + \bus({}^{110}_{111})\ksu({}^{000}_{110})
\end{array}\right] \\
&+ \left[\begin{array}{ll}
0 & \Euo({}^{010}_{111})\nsu({}^{000}_{010}) + \huo({}^{010}_{111})\bbsu({}^{000}_{010}) + \muo({}^{101}_{111})\Ebsu({}^{000}_{101}) + \buo({}^{101}_{111})\ksu({}^{000}_{101}) \\
0 & \Eus({}^{010}_{111})\nsu({}^{000}_{010}) + \hus({}^{010}_{111})\bbsu({}^{000}_{010}) + \mus({}^{101}_{111})\Ebsu({}^{000}_{101}) + \bus({}^{101}_{111})\ksu({}^{000}_{101})
\end{array}\right] \\
&+ \left[\begin{array}{ll}
0 & \Euo({}^{100}_{111})\nsu({}^{000}_{100}) + \huo({}^{100}_{111})\bbsu({}^{000}_{100}) + \muo({}^{011}_{111})\Ebsu({}^{000}_{011}) + \buo({}^{011}_{111})\ksu({}^{000}_{011}) \\
0 & \Eus({}^{100}_{111})\nsu({}^{000}_{100}) + \hus({}^{100}_{111})\bbsu({}^{000}_{100}) + \mus({}^{011}_{111})\Ebsu({}^{000}_{011}) + \bus({}^{011}_{111})\ksu({}^{000}_{011})
\end{array}\right] \\
\end{align*}
and claim that
\begin{align*}
\check{A} &= \ddd^2 \\ \\
\check{B} &= \ddd\mm + \mm\ddd \\ \\
\check{E} &= \ddd\hh + \hh\ddd + \mmb\mma + \mmd\mmc \\ \\
\check{F} &= \ddd\gc + \gc\ddd + \hhhd\mmma + \hhhe\mmmb + \hhhf\mmmc \\
&+ \mmmj\hhha + \mmmk\hhhb + \mmml\hhhc
\end{align*} \\
In each case, the additional terms on the left are those with no good break, and these terms occur in pairs (see the proof of Lemma \ref{lem:aij}).  We conclude that $\dd^2 = 0$ for $l \leq 3$.

\newpage

The equation
\begin{align*}
\ddd\hh + \hh\ddd = \mmb\mma + \mmd\mmc
\end{align*}
should be viewed as a generalization of the equation
\begin{align*}
\ddd \mw_P + \mw_P \ddd = \mw_Q.
\end{align*}
The latter holds whenever $Q$ is the boundary of a family $P$ of {\em non-degenerate} metrics (see \eqref{eqn:qnull}).  However, in the above case, $\check{Q}$ counts isolated trajectories over the two {\em degenerate} metrics at $\{-\infty, \infty\}$, which split the cobordism $W$ as $W_{00,01} \coprod W_{01,11}$ and $W_{00,10} \coprod W_{10,11}$.  We have the (non-vanishing) identities

\begin{align*}
\qoo &= \moob\mooa +\muob\esu\mosa + \duo\nsub\mosa\\
&+ \mood\mooc + \muod\esu\mosc + \duo\nsud\mosc\\
\qos &= \mosb\mooa + \musd\esu\mosc + \dus\nsud\mosc\\
&+ \mosd\mooc + \musd\esu\mosc + \dus\nsud\mosc\\
\quo &= \moob\muoa + \muob\esu\musa + \duo\nsub\musa + \muob\nuua\\
&+ \mood\muoc + \muod\esu\musc + \duo\nsud\musc + \muod\nuuc\\
\qus &= \mosb\muoa + \musb\esu\musa + \dus\nsub\musa + \nssb\musa + \musb\nuua\\
&+ \mosd\muoc + \musd\esu\musc + \dus\nsud\musc + \nssd\musc + \musd\nuuc\\
&\\
\rss &= \nssb\nssa + \nusb\nsua + \nssd\nssc + \nusd\nsuc \\
\rsu &= \nsub\nssa + \nuub\nsua + \nsud\nssc + \nuud\nsuc \\
\rus &= \nssb\nusa + \nusb\nuua + \nssd\nusc + \nusd\nuuc \\
\ruu &= \nsub\nusa + \nuub\nuua + \nsud\nusc + \nuud\nuuc
\end{align*} \\
which we bundle into the map \\
\begin{align*}
\check{Q} &= \left[\begin{array}{rr}
\qoo & \quo\esu + \duo\rsu \\
\qos & \rss + \qus\esu + \dus\rsu
\end{array}\right] \\
&+ \left[\begin{array}{rr}
0 & \muob\esu\dus\nsua + \muob\esu\nssa + \muod\esu\dus\nsuc+ \muod\esu\nssc\\
0 & \musb\esu\dus\nsua + \musb\esu\nssa + \musd\esu\dus\nsuc+ \musd\esu\nssc
\end{array}\right].
\end{align*} \\
The reader may now check directly that indeed \\
\begin{align*}
\check{Q} = \mmb\mma + \mmd\mmc.
\end{align*}

\newpage

With both the complex $(\tilde{X}, \dd)$ and the surgery exact triangle in mind, we now turn to the lattice $\{0,1,\infty\}$.  The map $\hh = \dd^0_\infty$ is given by

\begin{align*}
\hh &= \left[\begin{array}{rr}
\hoo & \huo\esu + \duo\ksu + \muo({}^1_\infty)\nsu({^0_1})\\
\hos & \kss + \hus\esu + \dus\ksu + \mus({}^1_\infty)\nsu({^0_1})
\end{array}\right]
\end{align*} \\
and we have the identities \\
\begin{align*}
\Eoo &= \hoo\doo + \doo\hoo + \duo\esu\hos + \duo\ksu\dos + \huo\esu\dos \\
&+ \muol\nsuk\dos + \mool\mook +\muol\esu\mosk + \duo\nsul\mosk\\
\Eos &= \hos\doo + \dos\hoo + \hus\esu\dos + \dus\ksu\dos + \dus\esu\hos + \kss\dos + \ess\hos \\
&+ \musl\nsuk\dos + \mosl\mook + \musd\esu\mosc + \dus\nsud\mosc\\
\Euo &= \hoo\duo + \doo\huo + \huo\esu\dus  + \duo\ksu\dus + \duo\esu\hus + \huo\euu +\duo\kuu \\
&+ \muol\nsuk\dus + \mool\muok + \muol\esu\musk + \duo\nsul\musk + \muol\nuuk\\
\Eus &= \hos\duo + \dos\huo + \hus\esu\dus  + \dus\ksu\dus + \dus\esu\hus + \hus\euu +\dus\kuu + \kss\dus + \ess\hus + \kus \\
&+ \musl\nsuk\dus + \mosl\muok + \musl\esu\musk \\
&+ \dus\nsul\musk + \nssl\musk + \musl\nuuk \\
&\\  
\Ebss &= \kss\ess + \ess\kss + \kus\esu + \eus\ksu + \nssl\nssk + \nusl\nsuk \\
\Ebsu &= \ksu\ess + \esu\kss + \kuu\esu + \euu\ksu + \nsul\nssk + \nuul\nsuk \\
\Ebus &= \kss\eus + \ess\kus + \kus\euu + \eus\kuu +  \nssl\nusk + \nusl\nuuk \\
\Ebuu &= \ksu\eus + \esu\kus + \kuu\euu + \euu \kuu + \nsul\nusk + \nuul\nuuk \\
\end{align*}
which we bundle into the map \\
\begin{align*}
\check{E} &= \left[\begin{array}{ll}
\Eoo & \Euo\esu + \duo\Ebsu + \huo\absu + \auo\ksu \\
\Eos & \Ebss + \Eus\esu + \dus\Ebsu + \hus\absu + \aus\ksu
\end{array}\right] \\
&+ \left[\begin{array}{ll}
0 & \buo({}^1_\infty)\nsu({}^0_1) + \muo({}^1_\infty)\bbsu({}^0_1) \\
0 & \bus({}^1_\infty)\nsu({}^0_1) + \mus({}^1_\infty)\bbsu({}^0_1)
\end{array}\right]
\end{align*} \\
The lattice $\{1,\infty, 0^\prime\}$ is treated similarly.  The reader may now check directly that

$$\check{E} = \hh\ddd + \ddd\hh + \mm({}^1_\infty) \mm({}^0_1).$$

\newpage

Finally, consider the lattice $\{0, 1, \infty, 0^\prime\}$ used in the proof of the surgery exact triangle.  The map $\check{G} = \dd^0_{0^\prime}$ is given by \\
\begin{align*}
\check{G} &= \left[\begin{array}{rr}
G^o_o & G^u_o\esu + \duo \bar{G}^s_u + \huo\nsu + \muo\ksu  + \duo \bar{m}^{ss}_u(n_s \otimes \cdot)\\
G^o_s & \bar{G}^s_s + \bar{m}^{ss}_s(n_s \otimes \cdot) + G^u_s\esu + \dus \bar{G}^s_u + \hus\nsu + \mus\ksu + \dus \bar{m}^{ss}_u(n_s \otimes \cdot)
\end{array}\right]
\end{align*} \\
This map is also described Section 5 of \cite{kmos}, as are the first four of the following identities for $\Aa = \Aa{}^0_{0^\prime}$.  We only record those identities which are used in the proof of the surgery exact triangle, and only those terms which break over $R_1$.  This is exactly the information that is needed to verify the cancellation in parts (ii) and (iii) of Lemma \ref{lem:aij2}. \\
\begin{align*}
\aoo &= \duo \bar{m}^{ss}_u(n_s \otimes \dos(\cdot)) + m^{uo}_o(\bar{n}_u \otimes \cdot) + \cdots \\
\aos &= \dus \bar{m}^{ss}_u(n_s \otimes \dos(\cdot)) + m^{uo}_s(\bar{n}_u \otimes \cdot) + \bar{m}^{ss}_s(n_s \otimes \dos(\cdot)) +  \cdots \\
\auo &= \duo \bar{m}^{ss}_u(n_s \otimes \dus(\cdot)) + \duo \bar{m}^{su}_u(n_s \otimes \cdot) + m^{uu}_o(\bar{n}_u \otimes \cdot) + \cdots \\
\aus &= \dus \bar{m}^{ss}_u(n_s \otimes \dus(\cdot)) + \dus \bar{m}^{su}_u(n_s \otimes \cdot) + m^{uu}_s(\bar{n}_u \otimes \cdot) + \bar{m}^{su}_s(n_s \otimes \cdot) + \bar{m}^{ss}_s(n_s \otimes \dus(\cdot)) + \cdots \\
& \\ 
\abss &= \bar{m}^{us}_s(\bar{n}_u \otimes \cdot) + \cdots\\
\absu &= \bar{m}^{us}_u(\bar{n}_u \otimes \cdot) + \cdots\\
& \\
\bar{A}^{ss}_s &= \bar{m}^{ss}_s(n_s \otimes \ess(\cdot)) + \bar{m}^{su}_s(n_s \otimes \esu(\cdot)) + \ess \bar{m}^{ss}_s(n_s \otimes \cdot) + \eus \bar{m}^{ss}_u(n_s \otimes \cdot) \\
\bar{A}^{ss}_u &= \bar{m}^{ss}_u(n_s \otimes \ess(\cdot)) + \bar{m}^{su}_u(n_s \otimes \esu(\cdot)) + \esu \bar{m}^{ss}_s(n_s \otimes \cdot) + \euu \bar{m}^{ss}_u(n_s \otimes \cdot)
\end{align*} \\
We then bundle these identities into the map $\Aa$ as in \eqref{eqn:exacta2}.  We also have the map $\check{L}$ as defined in \eqref{eqn:l}.  The reader may now verify directly that \\
$$\check{A} =  \check{L} + \check{G}\ddd + \ddd \check{G} + \hh\mm + \mm\hh.$$ \\

\noindent Though perhaps, at this point, he or she would prefer to take our word for it.

\bibliographystyle{hplain}
\bibliography{linksurg}

\begin{thebibliography}{10}

\bibitem{kat}
The~Knot Atlas.
\newblock http://katlas.org/.

\bibitem{b}
J.~Baldwin.
\newblock On the spectral sequence from {K}hovanov homology to {H}eegaard
  {F}loer homology.
\newblock 2008, arXiv:0809.3293.

\bibitem{b1}
J.~Bloom.
\newblock Odd {K}hovanov homology is mutation invariant.
\newblock 2009, arXiv:0903.3746v2.

\bibitem{bd}
P.J. Braam and S.~K. Donaldson.
\newblock Floer's work on instanton homology, knots, and surgery.
\newblock In H.~Hofer, editor, {\em The {F}loer Memorial Volume}, volume 133 of
  {\em Progress in Math.}, pages 195--256. Birkh$\ddot{\mathrm{u}}$auser, 1995.

\bibitem{carr}
M.~Carr and S.~Devadoss.
\newblock Coxeter complexes and graph-associahedra.
\newblock {\em Topology Appl.}, 153:2155--2168, 2006.

\bibitem{cass}
W.~Casselman.
\newblock Strange associations.
\newblock www.ams.org/featurecolumn/archive/associahedra.html, 2006.

\bibitem{chm}
S.~V. Chmutov and S.~K. Lando.
\newblock Mutant knots and intersection graphs.
\newblock 2007, arXiv:0704.1313v1.

\bibitem{d}
S.~Devadoss.
\newblock A realization of graph-associahedra.
\newblock {\em Discrete Math.}, 309:271--276, 2009.

\bibitem{don}
S.~K. Donaldson.
\newblock Topological field theories and formulae of {Casson} and
  {M}eng-{T}aubes.
\newblock In {\em Proceedings of the {K}irbyfest}, volume~2 of {\em Geometry
  and Topology Monographs}, pages 87--102, 1999.

\bibitem{fl}
A.~Floer.
\newblock Instanton homology, surgery, and knots.
\newblock In S.~K. Donaldson and C.B. Thomas, editors, {\em Geometry of
  Low-Dimensional Manifolds: Gauge Theory and Algebraic Surface}, volume 150 of
  {\em London Mathematical Society Lecture Note Series}, pages 97--114.
  Cambridge University Press, 1990.

\bibitem{gs}
R.~Gompf and A.~Stipsicz.
\newblock {\em 4-Manifolds and Kirby Calculus}.
\newblock Graduate Studies in Math. AMS, 1999.

\bibitem{gl}
C.~McA. Gordon and R.~A. Litherland.
\newblock On the signature of a link.
\newblock {\em Invent. Math.}, 47:53--69, 1978.

\bibitem{gr}
J.~Greene.
\newblock A spanning tree model for the {H}eegaard {F}loer homology of a
  branched double-cover.
\newblock 2008, arxiv/0805.1381v1.

\bibitem{kh1}
M.~Khovanov.
\newblock A categorification of the {J}ones polynomial.
\newblock {\em Duke Math. J.}, 101(3):359--426, 2000.

\bibitem{km}
P.~Kronheimer and T.~Mrowka.
\newblock {\em Monopoles and Three-Manifolds}.
\newblock New Mathematical Monographs. Cambridge University Press, 2007.

\bibitem{km2}
P.~Kronheimer and T.~Mrowka.
\newblock Knots, sutures and excision.
\newblock 2008, arXiv:0807.4891v2.

\bibitem{kmos}
P.~Kronheimer, T.~Mrowka, P.~Ozsv{\'a}th, and Z.~Szab{\'o}.
\newblock Monopoles and lens space surgeries.
\newblock {\em Ann. of Math. (2)}, 165(2):457--546, 2007, math.GT/0310164.

\bibitem{lod}
J.~Loday.
\newblock Realization of the {S}tasheff polytope.
\newblock {\em Arch. Math. (Basel)}, 83(3):267--278, 2004.

\bibitem{mcc}
J.~McCleary.
\newblock {\em A User's Guide to Spectral Sequences}.
\newblock Cambridge Studies in Advanced Mathematics. Cambridge University
  Press, 2 edition, 2000.

\bibitem{mo}
C.~Monalescu and P.~Ozsv{\'a}th.
\newblock On the {K}hovanov and knot {F}loer homologies of quasi-alternating
  links.
\newblock In {\em Proceedings of the 14th {G}{\"o}kova Geometry-Topology
  Conference}, pages 60--81, 2007.

\bibitem{moyu}
T.~Mrowka, P.~Ozsv{\'a}th, and B.~Yu.
\newblock {S}eiberg-{W}itten monopoles on {S}eifert fibered spaces.
\newblock {\em MSRI 1996-093}, 1996, arXiv:math/9612221v1.

\bibitem{orsz}
P.~Ozsv{\'a}th, J.~Rasmussen, and Z.~Szab{\'o}.
\newblock Odd {K}hovanov homology.
\newblock 2007, math.QA/0710.4300.

\bibitem{osz6}
P.~Ozsv{\'a}th and Z.~Szab{\'o}.
\newblock Absolutely graded {F}loer homologies and intersection forms for
  four-manifolds with boundary.
\newblock {\em Adv. Math.}, 173:179--261, 2003, math.SG/0110170.

\bibitem{osz12}
P.~Ozsv{\'a}th and Z.~Szab{\'o}.
\newblock On the {H}eegaard {F}loer homology of branched double-covers.
\newblock {\em Adv. Math.}, 194(1):1--33, 2005, math.SG/0309170.

\bibitem{osz5}
P.~Ozsv{\'a}th and Z.~Szab{\'o}.
\newblock Holomorphic triangles and invariants for smooth four-manifolds.
\newblock {\em Adv. Math.}, 202:326--400, 2006, math.SG/0110169.

\bibitem{lrob3}
L.P. Roberts.
\newblock Notes on the {H}eegaard-{F}loer link surgery spectral sequence.
\newblock 2008, math.GT/0808.2817.

\bibitem{umb}
S.~Saneblidze and R.~Umble.
\newblock Diagonals on the permutahedra, multiplihedra, and associahedra.
\newblock {\em Homology, Homotopy Appl.}, 6(1):363--411, 2004.

\bibitem{seid}
P.~Seidel.
\newblock {\em {F}ukaya categories and {P}icard-{L}efschetz theory}.
\newblock Zurich Lectures in Advanced Mathematics. European Mathematical
  Society, 2008.

\bibitem{stos}
M.~Sto\v s\'ic.
\newblock Homology of torus links.
\newblock 2007, arXiv:math/0606656v2.

\bibitem{stash}
J.~Stasheff.
\newblock Homotopy associativity of {H}-spaces {I}.
\newblock {\em Trans. Amer. Math. Soc.}, 108:275--292, 1963.

\bibitem{turner}
P.~Turner.
\newblock A spectral sequence for {K}hovanov homology with an application to
  (3,q)-torus links.
\newblock 2006, arXiv:math/0606369v1.

\bibitem{w}
L.~Watson.
\newblock A remark on {K}hovanov homology and two-fold branched covers.
\newblock 2008, arXiv:0808.2797.

\end{thebibliography}

\end{document}